\documentclass[12pt]{iopart}
\pdfoutput=1

\usepackage{iopams}
\usepackage{bbm}
\usepackage{bm}
\usepackage{graphicx}

\expandafter\let\csname equation*\endcsname\relax
\expandafter\let\csname endequation*\endcsname\relax

\usepackage{amsmath}

 \usepackage[caption=false]{subfig}
 \usepackage{floatrow}
 \usepackage{longtable}
 \usepackage{multirow}
\usepackage{adjustbox}
\newcommand{\strike}[1]{\ifmmode\setbox0\hbox{$#1$}%
\else
\setbox0\hbox{#1}%
\fi
\makebox[\the\wd0][c]{%
\rule[0.48\ht0]{0.5\wd0}{0.25pt}}\hspace*{-\the\wd0}#1}

\usepackage{color}

\begin{document}
\title{
Heteroclinic units acting as pacemakers:
Entrained dynamics for cognitive processes}

\author{Bhumika Thakur and Hildegard Meyer-Ortmanns}

\address{School of Science, Jacobs University Bremen, Campus Ring 1, 28759 Bremen, Germany}
\ead{b.thakur@jacobs-university.de}
\ead{h.ortmanns@jacobs-university.de}
\vspace{10pt}
\begin{indented}
\item[]July 2022
\end{indented}

\begin{abstract}
Heteroclinic dynamics is a suitable framework for describing transient and reproducible dynamics such as cognitive processes in the brain. We demonstrate how heteroclinic units can act as pacemakers to entrain larger sets of units from a resting state to hierarchical heteroclinic motion that is able to describe fast oscillations modulated by slow oscillations. Such features are observed in brain dynamics. The entrainment range depends on the type of coupling, the spatial location of the pacemaker and the individual bifurcation parameters of the pacemaker and the driven units. Noise as well as a small back-coupling to the pacemaker facilitate synchronization. Units can be synchronously entrained to different temporal patterns encoding transiently excited neural populations, depending on the selected path in the heteroclinic network. Via entrainment, these temporal patterns, locally generated by  the pacemakers, can be communicated to the resting units in target waves over a spatial grid. For getting entrained there is no need of fine-tuning  the parameters of the resting units. Thus, entrainment provides one way of processing information  over the grid, when information is encoded in the generated spatiotemporal patterns.
\end{abstract}

\section{Introduction}
The dynamics of a number of systems is intrinsically transient even if it may look stationary. The transients may last long as compared to subsequent time intervals over which the state of the system then drastically changes. Examples are seen in ecological systems in extinction or switching events. These events refer to the extinction or switching of (part of) some species, or to the strong suppression of (sub)populations, in the latter case before another (sub)population becomes dominant. A second important branch of transient dynamics is brain dynamics. Brain dynamics proceeds via sequential segmentation of information that is manifested in sequences of electroencephalography (EEG) microstates \cite{EEG}. Cognitive processes are obviously transient, on the other hand they are well reproducible. Thus for a long time it was not clear which dynamical framework is able to capture the transient characteristics, but is reproducible at the same time and stable with respect to not too strong perturbations. As it turned out, a promising candidate for such a mathematical framework is heteroclinic dynamics \cite{rabi3,laurent,binding,rabi1}, to be defined later.

Hints on  underlying heteroclinic structures of the dynamics in ecological systems are seen in extinction and switching events \cite{valentin1,nowak}, for a review see, for example, \cite{szolnoki}. In relation to brain dynamics, in particular motor processing in the brain, heteroclinic cycles provide a possible explanation  for various gaits in animal and human motion \cite{friston}. Heteroclinic dynamics has been also discussed in relation to binding \cite{binding} and chunking dynamics \cite{rabi1} of the brain. Binding refers to the composition of information from different brain areas such as the visual and auditory system, possibly realized via heteroclinic connections. Chunking dynamics amounts to cutting long sequences of information into shorter ones (the chunks), for example  for better memorizing  long sequences. A possible realization in the brain is in terms of slow oscillations modulating fast oscillations, cutting the fast oscillations into ``chunks", corresponding to one period of the slow oscillations. Such motion can be reproduced in heteroclinic dynamics \cite{rabi1,max1,max2,max3}. Later we will consider this type of heteroclinic motion as a prototype of hierarchical  heteroclinic motion because of the two inherent time scales, the longer one referring to the modulation of fast oscillations.

It is well known that often a single input is sufficient to trigger a whole sequence of cognitive items from different areas of the brain which is then executed in a reproducible way. Assuming that this kind of transient dynamics corresponds to synchronized neural activities of sub-populations of neurons, the question arises as to what extent the involved neural entities must be in the ``right state''. Stated differently, how stable are these transient synchronization patterns against heterogenous parameters in the individual entities if they are described as heteroclinic units, assigned to a grid and obeying heteroclinic dynamics? Is it possible that a single heteroclinic unit (to be defined later) acts as a pacemaker and entrains other such units even if they are in a kind of resting state, a state, in which they would not perform oscillations on their own? Later we consider the maintenance and restoration of synchronization in larger sets of heteroclinic units  when  part of them is in a non-oscillatory (resting) state.
\footnote{A related question -though not in relation to brain dynamics- was addressed in \cite{daido1,daido2} where globally and diffusively coupled oscillators under the influence of ageing were studied. Ageing there was understood in the sense that a subset of oscillators fails to sustain their individual oscillations. Beyond a certain level of ageing, the collective oscillation stops in a so-called ageing transition with universal scaling behavior of the order parameter.  The models comprise Stuart-Landau oscillators and R\"ossler units in the oscillatory or chaotic regime.}

In relation to biological  applications, the action of pacemakers has been studied for coupled Kuramoto oscillators as paradigm for synchronization of limit cycles, for example in \cite{radicchi1,radicchi2}.
The entrained motion of coupled Kuramoto oscillators are simple limit cycles. In contrast, already the dynamics of a single heteroclinic unit can be quite versatile. The type we consider later either approaches an equilibrium with coexisting items, here called a resting state, or performs oscillations with one or more inherent time scales, imitating chunking dynamics. Thus  one may wonder whether a heteroclinic unit can act as a pacemaker and entrain other units to the same kind of intricate motion   when the other units are themselves in a resting state. The entrainment should be a result of the coupling.
\footnote{Heteroclinic units under weak periodic forcing have been considered in \cite{rabi2} and \cite{dawes1}. In our system, the periodic forcing is effectively provided by a pacemaker that itself is a heteroclinic unit under the action of noise (effectively acting like periodic forcing), but our driven units are chosen to be individually  in the regime of stable coexistence equilibria, differently from \cite{rabi2}.} What are the conditions for this entrainment and what determines its range?

One may wonder why such an entrainment to specific temporal sequences should be useful in view of brain dynamics.
The idea that information about the environment is represented in nervous systems through both the identity of an object (its spatial coordinate) and through temporally ordered sequences of excitation signals was proposed in \cite{rabi3}, based on experimentally observed features of olfactory processing networks \cite{laurent}. Even static signals (as from odor) are classified and coded by using time as an additional coordinate in storing information of the firing activity. Obviously, this way the storage capacity is highly increased \cite{rabi3}.

In view of this conjecture let us anticipate one of our main results as it illustrates that heteroclinic dynamics is able to realize temporal coding and process this information over a spatial grid via entrainment to the pacemaker dynamics.
We consider a few heteroclinic units in a small disc at the center of a two-dimensional grid, each of them performing a kind of chunking dynamics with a certain sequence of dominantly excited sub-populations, here labelled by integer numbers $\{1,2,...,9\}$. Say the parameters of the individual units are chosen such that the order of the temporarily dominant sub-populations   is $(1,2,3)\rightarrow(6,4,5)\rightarrow(8,9,7)$ with cyclic repetition of this sequence. All neighboring  units outside the disc are individually in a resting state. When these resting units are  coupled to the units of the disc, they get entrained to the same sequence of numbers with an amplitude that decays with the distance to the disc, but with high precision with respect to the sequence. This means that the information encoded in the temporal order of the labels becomes accessible to many units in the resting state, when the information spreads over the grid in a target wave (as demonstrated in the movie of \cite{supp}). \footnote{Emergence of complex spatiotemporal wave patterns  are observed, for example, in primates' cerebral cortex, and analyzed in \cite{townsend1,townsend2}.}

Not surprisingly, the range of synchronization with the pacemaker units is finite, and it is part of our analysis which parameters and coupling topology determine this range. Our results on the two-dimensional grid are based on numerical simulations.
To obtain a partial analytical understanding, we reduce the two-dimensional grid to a chain with  a single pacemaker and $N-1$ driven units. In more detail we consider two heteroclinic units without a hierarchy in time scales. The extension to $N>2$ reveals a proliferation of  saddle equilibria, so one  may wonder whether synchronization is still achievable. However, as it turns out, noise facilitates synchronization, precise entrainment works best for an intermediate level of noise. Our results may therefore explain why the precision of stochastic Turing patterns, as considered in \cite{barato}, is maximized for an intermediate thermodynamic cost \cite{barato} spent on suppressing thermodynamic fluctuations.

The paper is organized as follows. In section \ref{sec2} we present heteroclinic dynamics realized via generalized Lotka-Volterra equations and specify the chosen topologies as well as the implementation of the pacemaker. In section \ref{sec3} we consider heteroclinic units with one inherent time scale. We provide bifurcation analyses for two  coupled units to explain some qualitative features of entrainment for larger sets, coupled along chains  of N units. Here we study also the effect of noise. The focus of section \ref{sec4} is on heteroclinic units with an inherent hierarchy in  time scales. We study the entrainment to trajectories between different heteroclinic motion, in particular the dependence of the temporal order on the chosen path in the heteroclinic network. In section \ref{sec5} we address the question of maintenance of oscillations in the presence of randomly distributed units  in their resting state. Moreover, we discuss the entrainment of resting-state units via target waves. The target waves are emitted from pacemakers when the pacemakers are restricted to a small disc in the center of the grid. Section \ref{sec6} contains our summary and conclusions.

\section{The model}\label{sec2}
We first summarize some basic notions of heteroclinic dynamics. For more formal definitions we refer to the original literature (section \ref{sec21}). Next we consider generalized Lotka-Volterra (GLV) equations as one option of realizing winnerless competition in heteroclinic dynamics (section \ref{sec22}). The concept of a pacemaker is introduced in section \ref{sec23}, followed by a table of our nomenclature in section \ref{sec24}.

\subsection{Heteroclinic dynamics}\label{sec21}
When the unstable manifold of a saddle equilibrium intersects the stable manifold of another saddle, the intersection is called a heteroclinic orbit. A heteroclinic network  is a set of vertices, corresponding to saddles, connected by edges, which are heteroclinic orbits. Thus heteroclinic networks are networks in phase space. Saddles in these networks are not restricted to saddle equilibria, but refer to any invariant set with non-trivial stable and unstable manifolds.  Early seminal papers on heteroclinic networks are \cite{kirk,ashwin1999}, see, for example, also \cite{ashwin2016}. A simple special case is a heteroclinic cycle between saddle equilibria, in which the unstable direction of one saddle becomes the stable direction of the subsequent saddle. The trajectory associated with the heteroclinic cycle approaches the vicinity of the saddles in a cyclic way. Upon its revolutions, it comes closer and closer to the saddles, but slows down and reaches the heteroclinic cycle in the infinite-time limit. The slowing down may be weakened or entirely prevented by exposing the system to a small amount of noise.

In view of our applications we couple heteroclinic networks on a spatial grid, as considered in \cite{max3}. To disentangle heteroclinic networks (which are networks in phase space) from network aspects of the spatial grid, we term heteroclinic networks assigned to a single site of the grid a heteroclinic unit. Thus the coupling of heteroclinic units refers to their coupling in coordinate space. This construction allows to share typical features of brain dynamics: According to \cite{grossberg}, the overall organization of the brain is in terms of parallel processing streams with complementary properties, hierarchical interactions within each stream, and parallel interactions between the streams. Our assignment of heteroclinic units to different sites of the grid realizes such a performance that is ongoing in parallel at different locations, while the imposed hierarchy of individual heteroclinic networks in phase space implements processes with hierarchies  in time scales such as the modulation of fast oscillations via slow oscillations.

\subsection{Generalized Lotka-Volterra equations}\label{sec22}
The system of heteroclinic units is defined on regular grid topologies with nodes labelled by $k$. To each node, we assign a heteroclinic unit that obeys GLV-equations according to
\begin{equation}\label{eq:1}
 \partial_t s_{k,i} = \rho s_{k,i} - \gamma_{k} s_{k,i}^2 - \sum_{j\neq i} A_{i,j} s_{k,i} s_{k,j} + \sum_{l} K_{k,l} (s_{l,i}-s_{k,i}) + \sigma_k |\xi_i(t)|.
\end{equation}
The meaning of the variables $s_{k,i}$ depends on the intended application. In an ecological context, $s_{k,i}$ represent the concentration of the $i^{th}$ species of the $k^{th}$ GLV unit, $\rho$ is the net reproduction rate, chosen as $1$ to set the time scale, $\gamma_k$ alters the net reproduction rate of the $k^{th}$ unit in a density-dependent way. Throughout the paper we focus on possible applications to brain dynamics, in particular transient cognitive processes. Here we interpret $s_{k,i}$ as densities of neural sub-populations. During the winnerless competition, a sequence of saddles is approached. In the vicinity of individual saddles various specific sub-populations get temporarily dominantly excited. Information (for example from sensory input) is believed to be encoded  in the specific spatiotemporal excitation pattern. In the vicinity of a saddle, the variable $s_{k,i}(t)$ tells us which subpopulation (i) gets dominantly excited at what location (k) and at what time. Thus,  information is contained in the very selection of the specific saddle and its time coordinate within the ordered temporal sequence of approached saddles as well as in the sequence itself. In short, we term the variables $s_{k,i}$  the information items or just the ``items'' in the following, without specifying the neural populations any further.

The parameter $\rho$ is again the net reproduction rate and $\gamma_k$ the item-density dependent decay rate, $\sigma_k$ the strength of the additive noise to the $k^{th}$ unit, and $\xi$ is Gaussian white noise with zero mean. $A$ is called the rate matrix with $A_{i,j}$ being the competition rate with which the item $j$ of a unit acts  on  item $i$ of the same unit. The interaction between the units is determined by the coupling matrix $K$: $K_{k,l}$ is finite if unit $l$ influences unit $k$, otherwise it is zero.

\subsubsection{GLV-equations for three items}
When each unit comprises of three items and the items play a rock-paper-scissors game, the rates characterize intransitive (non-hierarchical) competition. For the brain dynamics they characterize the `competition' strength in a winnerless game of items, the winnerless competition among its items then results in a single hierarchy level (one time scale), and $A$ is a $3\times3$  matrix given by
  \begin{equation}\label{eq:2}
 A =
\begin{pmatrix}
0 & c & e\\
e & 0 & c\\
c & e & 0
\end{pmatrix}.
 \end{equation}
As studied previously \cite{max1,max2}, there is a critical value of the decay rate $\gamma:\gamma_{crit}=(c+e)/2$, at which a single unit undergoes a Hopf bifurcation. For $\gamma_{k}<\gamma_{crit}$, the dynamics of the $k^{th}$ unit is characterized by a stable heteroclinic cycle connecting the saddle equilibria $\xi_i=\{\vec{s_k}\in \mathbb{R}^3: s_{k,i}=\rho/\gamma_k,s_{k,j}=0 \; \forall j \neq i \}$. For $\gamma_k>\gamma_{crit}$, the coexistence equilibrium  $\xi^*=\{\vec{s_k}\in \mathbb{R}^3:s_{k,i}=\frac{\rho}{c+e+\gamma_k} \forall i\} $ is stable.
Throughout the manuscript, we fix $c=2.0$ and $e=0.2$, and therefore, $\gamma_{crit}=1.1$. At $\gamma=\gamma_{crit}$, due to an emerging conservation of the sum and product of three items for $\gamma_{crit}$, the system has an infinity of periodic orbits, depending on the choice of initial conditions, with typical signatures of critical behavior \cite{max4}.

\subsubsection{GLV-equations for nine items}
When each unit contains nine items, the rate matrix $A$ is chosen in such a way that the winnerless competition among its items has an attractor that is a (large) heteroclinic cycle (LHC) between three saddles that are themselves (small) heteroclinic cycles (SHCs) between three saddle equilibria. As we have shown in \cite{max1,max2}, this structure of the attractor induces a hierarchy in time scales of slow oscillations (due to LHCs) modulating fast oscillations (due to SHCs).  The long (short) time scale amounts to one revolution in the long (short) heteroclinic cycle, respectively, which is well defined as long as the slowing down is suppressed by the application of a small amount of noise, otherwise it is itself a function of time.
Such a modulation of oscillations is a typical feature that is also experimentally observed in brain dynamics (see, for example, \cite{moelle,saleh}).
The rate matrix $A$ for the units with two hierarchy levels is then given as the following block matrix
\begin{equation}\label{eq:3}
 A=\begin{pmatrix}
m_0 & m_d & m_f\\
m_f & m_0 & m_d\\
m_d & m_f & m_0
\end{pmatrix},
\text{where }
 m_0=
\begin{pmatrix}
0 & c & e\\
e & 0 & c\\
c & e & 0
\end{pmatrix}.
\end{equation}
The matrix $m_d$ is also $3\times 3$ and has rates $d$ on the diagonal elements and $r$ on the remaining elements. Similarly, the matrix $m_f$ has a rate $f$ on the diagonal with all off-diagonal elements being $r$.
Again, on tuning the decay rate $\gamma$, the system undergoes a sequence of Hopf bifurcations, whose order and values depend on the choice of the other parameters \cite{max1,max2}. Two types of Hopf bifurcations are relevant to us here. For low values of $\gamma$, the system has an intricate heteroclinic dynamics with two hierarchy levels. When $\gamma$ is increased towards a bifurcation point $\gamma_c$, the system undergoes a Hopf bifurcation, at which the highest hierarchy level is gone. What remains is an LHC between three 3-items saddle equilibria as a remnant from the SHCs. Further increasing $\gamma$, the system undergoes a second Hopf bifurcation which drives the unit into stable global coexistence of nine items with the same concentrations.
In extension of such a  single hierarchical heteroclinic network, the dynamics on a spatial network (a two-dimensional square lattice) of identical such heteroclinic units  was analyzed in \cite{max3}, where the units were coupled via diffusion.

The very construction of the adjacency matrix $A$ in phase space may appear rather peculiar. It serves to enforce hierarchical heteroclinic motion just by the choice of rates. The rates are chosen out of certain intervals, so that the eigenvalues of the corresponding Jacobians at the saddle equilibria lead to  a preference of the trajectory for  one (desired) direction out of two unstable directions when escaping from the saddles. This construction was based on results of \cite{ashwin2013}. A hierarchy in time scales can be implemented differently, see e.g. \cite{rabi1}. In the brain, external input may be responsible for selecting certain paths in the (counterparts of) heteroclinic networks. We do not consider external input here. Our choice is just for convenience to generate an intricate motion that should be entrained, independently of how it was generated.

\subsection{Implementation of the pacemaker}\label{sec23}
In chemical oscillatory systems,  the role of pacemakers seems to be played by defects and impurities \cite{kurabook}. Accordingly, in the simplest case we treat the heteroclinic unit which should play the role of a pacemaker on the same footing as all other units (the driven units), but assign a ``defect'' in the sense that we distinguish its individual parameters from those of all other units: for given parameters $c$ and $e$ we choose $\gamma$ in the range of the highest (here at most second) hierarchy level of motion, so that this unit performs either a heteroclinic cycle of heteroclinic cycles (nine items) or a heteroclinic cycle between saddle equilibria (three items), while the other units are chosen in the range of $\gamma$ for which the dynamics converges to a stable equilibrium point with coexisting items. In the following we call the corresponding parameter regimes the heteroclinic regime (HC) or the coexistence regime (CE).
\begin{figure}[ht!]
    \centering
\includegraphics[width=\textwidth]{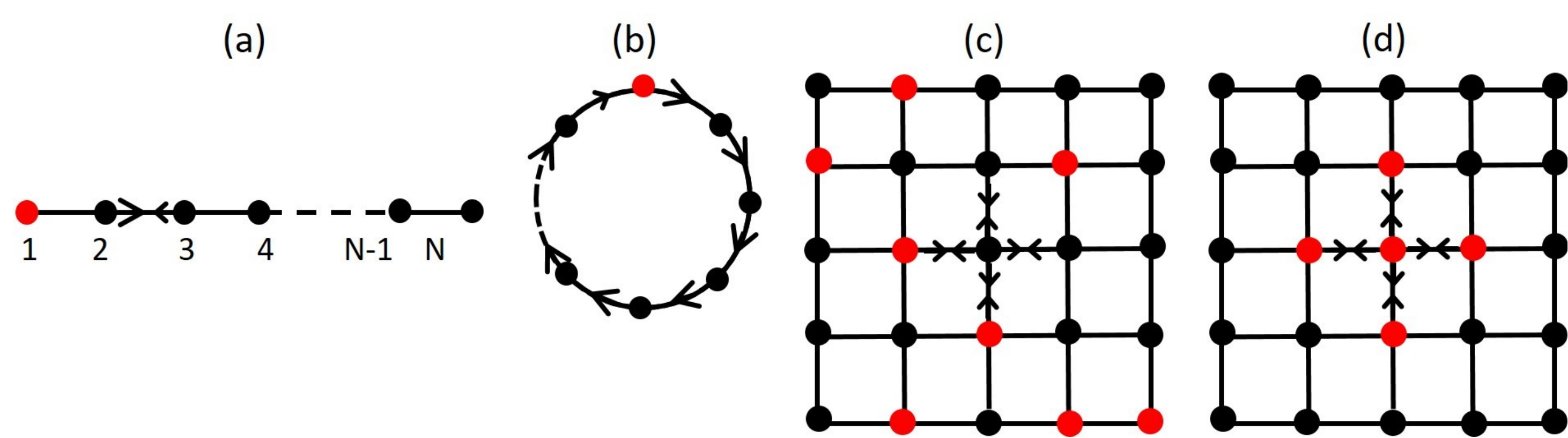}
\caption{Grid configurations (a) chain, (b) ring, (c) and (d) square lattices, in general with stronger forward coupling $\delta$ than back-coupling $\delta_b$. Pacemaker units are indicated  in red, driven units in black.
}
 \label{fig:1}
 \end{figure}
In addition, the pacemaker may be distinguished by its location and its coupling to the driven units. Along a chain, for local coupling, the pacemaker will be distinguished by its position at one of the edges of the chain and directed local coupling along all units with a smaller back-coupling than forward coupling. When we close the chain to a ring, we assign a smaller back-coupling to the pacemaker at least along the last link that leads to the closure.  In \cite{supp} we discuss  an implementation also with distance-dependent couplings of driven units to the pacemaker. The chains (and rings in \cite{supp}) are studied for varying size, in general with stronger forward than backward coupling. On the two-dimensional grid, the coupling is  bidirectional (diffusive), in order to compare with earlier results  of \cite{max3}, defects are either randomly assigned or confined to a region outside a small disc of ``non-defects'' at the center of the grid. The different options are displayed in Fig.1.

\subsection{Nomenclature}\label{sec24}
The nomenclature which we use in this paper is summarized in Table \ref{table:1}.
Some remarks to the table are in order. We call oscillations to be large (LHCs) or small (SHCs) heteroclinic cycles, even if the characteristic slowing down of the HCs is suppressed by a small amount of noise and it would be more precise to call the motion a quasi-limit cycle in the vicinity of the heteroclinic orbit. This is to distinguish such oscillations from limit cycles remote from the heteroclinic  orbit. Such cycles are also encountered  in this work, they are called  small limit cycles (SLCs). The attributes small or large do not reflect the distance in phase space, but are merely due to the visualization and correspond to hierarchy levels with slow (large cycle) or fast (small cycle) time scales.

\begin{table}[]
\begin{small}
\begin{tabular}{ll}
\multicolumn{2}{}{}                                                                                                                                                                                       \\ \hline
\multicolumn{1}{|l|}{HCU}        & \multicolumn{1}{l|}{heteroclinic unit, satisfying Eq.~1 with Eq.~2 or Eq.~3 at a single site}                                                                                                           \\ \hline
\multicolumn{1}{|l|}{$\rho$}     & \multicolumn{1}{l|}{net reproduction rate, set to 1 to set the timescale}                                                                                                                                \\ \hline
\multicolumn{1}{|l|}{$c$}        & \multicolumn{1}{l|}{ rate, set to 2.0}                                                                                                                                   \\ \hline
\multicolumn{1}{|l|}{$e$}        & \multicolumn{1}{l|}{ rate, set to 0.2}                                                                                                                                   \\ \hline
\multicolumn{1}{|l|}{$f$}        & \multicolumn{1}{l|}{ rate, taken as either 0.3 or 0.4}                                                                                                                                                               \\ \hline
\multicolumn{1}{|l|}{$r$}        & \multicolumn{1}{l|}{rate, set to 1.25}                                                                                                                                  \\ \hline
\multicolumn{1}{|l|}{$\gamma$}   & \multicolumn{1}{l|}{decay rate, used as a bifurcation parameter for an HCU}                                                                                                        \\ \hline
\multicolumn{1}{|l|}{$\delta$}   & \multicolumn{1}{l|}{forward coupling between HCUs on a grid}                                                                                                                        \\ \hline
\multicolumn{1}{|l|}{$\delta_b$} & \multicolumn{1}{l|}{backward coupling between HCUs on a grid}                                                                                                                       \\ \hline
\multicolumn{1}{|l|}{HC}         & \multicolumn{1}{l|}{\begin{tabular}[c]{@{}l@{}}heteroclinic cycle for 3 items: between 3 1-item saddle equilibria;\\ for 9 items: between 3 3-items saddle equilibria\end{tabular}} \\ \hline
\multicolumn{1}{|l|}{SHC}        & \multicolumn{1}{l|}{\begin{tabular}[c]{@{}l@{}}small HC (quasi-limit cycle) between 3 1-item equilibria if the trajectory \\ is close to the heteroclinic orbit\end{tabular}}       \\ \hline
\multicolumn{1}{|l|}{LHC}        & \multicolumn{1}{l|}{large HC between 3 SHCs}                                                                                                                                        \\ \hline
\multicolumn{1}{|l|}{SLC}        & \multicolumn{1}{l|}{\begin{tabular}[c]{@{}l@{}}small limit cycles (if the trajectory is remote from the corresponding \\ heteroclinic orbit)\end{tabular}}                          \\ \hline
\multicolumn{1}{|l|}{HC-regime}  & \multicolumn{1}{l|}{regime with HCs with 1 or 2 timescales}                                                                                                                         \\ \hline
\multicolumn{1}{|l|}{CE-regime}  & \multicolumn{1}{l|}{regime with coexistence equilibria (fixed-points)}                                                                                                              \\ \hline
\multicolumn{1}{|l|}{1L-unit}    & \multicolumn{1}{l|}{\begin{tabular}[c]{@{}l@{}}HCU with 3 or 9 items, each in the respective HC-regime with 1 timescale;\\ used for pacemakers\end{tabular}}                        \\ \hline
\multicolumn{1}{|l|}{2L-unit}    & \multicolumn{1}{l|}{HCU with 9 items in the HC-regime with 2 timescales; used for pacemakers}                                                                                       \\ \hline
\multicolumn{1}{|l|}{CE-unit}    & \multicolumn{1}{l|}{HCU with 3 or 9 items in the CE regime; used for driven units}                                                                                                  \\ \hline
\end{tabular}
\end{small}
\caption{Nomenclature}
\label{table:1}
\end{table}

For systems with noise and on a two-dimensional grid, we use XMDS2 software \cite{xmds2} for the numerical integration; everywhere else Matlab is used. For analytical calculations, we make use of the Mathematica software.

\section{Heteroclinic units with one level of hierarchy}\label{sec3}
For two units, a pacemaker and a driven unit, an analytical understanding is possible in terms of an effective model for the driven unit, exposed to stepwise constant forcing  by the pacemaker (\ref{subsec31} and Appendix A). The extension to $N-1$ driven units reveals a proliferation of  saddle equilibria which are in principle accessible to the driven units, the more, the larger the distance to the pacemaker. We discuss this proliferation (\ref{subsec32} and Appendix B.1), because it is the constructive role of noise that turns out to facilitate a precise entrainment in spite of the multitude of saddles.

\subsection{Bifurcation analyses of two coupled units}\label{subsec31}
We start with a system of two 1L-units and three items, with one unit acting as a pacemaker to the other (driven) unit. The pacemaker is at location $k=1$ with decay rate  $\gamma_P$, the driven unit at $k=2$
with decay rate $\gamma_D$.
The governing equation (Eq.~\ref{eq:1}) takes the following form:
\begin{eqnarray}\label{eq:4}
 \partial_t s_{1,i} &=& \rho s_{1,i} - \gamma_{P} s_{1,i}^2 -s_{1,i}(c s_{1,i+1} + e s_{1,i+2})+\delta_b (s_{2,i} - s_{1,i}),\nonumber\\
 \partial_ts_{2,i} &=& \rho s_{2,i} - \gamma_D s_{2,i}^2 - s_{2,i} (c s_{2,i+1} + e s_{2,i+2})+\delta (s_{1,i} - s_{2,i}),
\end{eqnarray}
where we denote the interaction strength from the pacemaker to the driven unit, that is $K_{2,1}$, by $\delta$ and the back-coupling from the driven unit to the pacemaker ($K_{1,2}$) by $\delta_b$.

We first consider $\delta_b=0$.
The $\gamma_P$ of the pacemaker is chosen such that it is characterized by a stable heteroclinic cycle
that connects three saddles of the pacemaker $(\rho/\gamma_P,0,0), (0,\rho/\gamma_P,0),(0,0,\rho/\gamma_P)$ via heteroclinic orbits. Due to the characteristic slowing-down of the heteroclinic cycle near the saddles, the trajectory spends most of the time in the vicinity of one of the saddles. During this dwell time,  only one of its items is dominating (this is winnerless competition with one temporary winner), while the concentration of the other two items approaches zero. Thus the driven unit   feels   a constant forcing to one of its items during that time.

This leads us to an effective modeling of the driven unit exposed to constant forcing that we present in detail in Appendix A. The result is that depending on the forcing, parameterized by the forward coupling $\delta$, the driven unit approaches a 1-item equilibrium (1S), a 2-items equilibrium (2S), or a 3-items equilibrium (3S), where the numbers (1,2,3) refer to the number of coexisting dominating items. In the original model, as the pacemaker approaches one of its saddles, the driven unit approaches one of these equilibria. If the equilibrium approached is a stable equilibrium of the effective model, it stays there until the pacemaker moves towards another of its saddles, thus changing forcing. Otherwise, if it approaches an unstable equilibrium, the driven unit moves towards the stable equilibrium of the effective model while the pacemaker continues to reside in the vicinity of its saddle.

\subsection{Chains of N heteroclinic units: Zooming into the proliferation of equilibria of a driven unit}\label{subsec32}
We next consider larger systems of coupled units in a chain  configuration as shown in Fig.~\ref{fig:1}(a). If we take one 1L-unit as the pacemaker  and all other $N-1$ units as CE-units, under what conditions will the pacemaker be able to entrain the other units to heteroclinic motion? We consider first unidirectional couplings,  bidirectional ones  with a strong forward and a weak back-coupling are discussed in \cite{supp}.
The unidirectional  nearest-neighbour couplings along the chain are chosen according to the coupling matrix $K$, given as

$K_{k,l} =
\begin{cases}
    \delta,& \text{if } k\in \{2,...,N \} \text{ and } l=k-1\\
        0,              & \text{otherwise}.
\end{cases}
$

We place a $1L$-unit as the pacemaker at the edge of the chain at $k=1$ with $\gamma_1=\gamma_P<1.1$.
For all  other CE-units we choose ($\gamma_k=\gamma_D>1.1 \;\forall k\in \{2,...,N \}$). Note that in this configuration, a unit at location $k-1$ acts effectively as a pacemaker to the $k^{th}$ unit which itself drives the unit at  $k+1$. Therefore, each unit directly or indirectly feels the effect of all the units to its left but is not affected by the dynamics of the units to its right side in the schematic Fig.~\ref{fig:1}(a).

In Appendix B.1 we zoom into the proliferation of equilibria in such a simple chain as it can be pursued to some extent what causes the proliferation.
The dynamics of the $n^{th}$ unit along the chain depends on the dynamics of all the $n-1$ units to its left.
Even though there are $3 n (n+1)/2$ possible equilibria in the phase space of the $n^{th}$ unit that it can visit in principle as shown in Table \ref{table:2}, it effectively gets to visit only a few of them. Which path in the branching tree of possible visits it actually visits depends on a number of factors such as
\begin{enumerate}
 \item which equilibria its (predecessor) pacemaker visits and for how long;
 \item whether the approached equilibrium  is a stable or unstable equilibrium of the corresponding system with constant forcing. If it gets kicked to an unstable equilibrium, it can move to a more stable one if its experienced forcing  does not change during that time interval.
 \item how closely it approaches an unstable equilibrium. The closer it comes, the longer it stays in the vicinity of the saddle equilibrium.
 \item the dwell time of the pacemaker at $k=1$ in the vicinity of its equilibria.
\end{enumerate}
All these factors in turn depend on the  bifurcation parameter of the pacemaker $\gamma_P$, the bifurcation parameter of the driven units $\gamma_D$, and the coupling along the chain $\delta$. This explains why a precise prediction of the path that is finally taken is practically not feasible for these chains.

Moreover, in Appendix B.2 we present the impact of $\delta$, $\gamma_P$ and $\gamma_D$ on the entrainment length in some detail. Here is the summary: The stronger the forcing $\delta$, the more units can be entrained, as expected. The more remote a given driven unit is from the original pacemaker at the end of the chain, the more the trajectory wiggles on its approach of the heteroclinic cycle due to the proliferated number of saddles, felt by the trajectory. Furthermore, the entrainment length increases with $\gamma_P$, which seems counterintuitive at a first view, and also decreases with an increase in $\gamma_D$ as expected. For further details we refer to Appendix B.

\begin{table}[]
 \begin{adjustbox}{width=\columnwidth,center}
\begin{tabular}{|l|lllllllllllllll|l|l|}
\hline
Unit-1 & \multicolumn{15}{l|}{1S}                                                                                                                                                                                                                                                                                                                                                       & 1        & 3         \\ \hline
Unit-2 & \multicolumn{10}{l|}{1S}                                                                                                                                                                                                                                          & \multicolumn{4}{l|}{2S}                                                                               & 3S & 3        & 9         \\ \hline
Unit-3 & \multicolumn{6}{l|}{1S}                                                                                                                                   & \multicolumn{3}{l|}{2S}                                                     & \multicolumn{1}{l|}{3S} & \multicolumn{3}{l|}{2S}                                                     & \multicolumn{1}{l|}{3S} & 3S & 6        & 18        \\ \hline
Unit-4 & \multicolumn{3}{l|}{1S}                                                     & \multicolumn{2}{l|}{2S}                           & \multicolumn{1}{l|}{3S} & \multicolumn{2}{l|}{2S}                           & \multicolumn{1}{l|}{3S} & \multicolumn{1}{l|}{3S} & \multicolumn{2}{l|}{2S}                           & \multicolumn{1}{l|}{3S} & \multicolumn{1}{l|}{3S} & 3S & 10       & 30        \\ \hline
Unit-5 & \multicolumn{1}{l|}{1S} & \multicolumn{1}{l|}{2S} & \multicolumn{1}{l|}{3S} & \multicolumn{1}{l|}{2S} & \multicolumn{1}{l|}{3S} & \multicolumn{1}{l|}{3S} & \multicolumn{1}{l|}{2S} & \multicolumn{1}{l|}{3S} & \multicolumn{1}{l|}{3S} & \multicolumn{1}{l|}{3S} & \multicolumn{1}{l|}{2S} & \multicolumn{1}{l|}{3S} & \multicolumn{1}{l|}{3S} & \multicolumn{1}{l|}{3S} & 3S & 15       & 45        \\ \hline
Unit-n & \multicolumn{15}{l|}{}                                                                                                                                                                                                                                                                                                                                                         & n(n+1)/2 & 3n(n+1)/2 \\ \hline
\end{tabular}
\end{adjustbox}
\caption{Proliferation of possible equilibria (stable or unstable) in the phase space of each unit for a chain of 1L-units with unit-1 being the driving pacemaker.
The last two columns denote the total number of possible equilibria when the pacemaker is near one of its equilibria, or when it has completed a whole cycle, respectively.}
\label{table:2}
\end{table}

\subsection{Facilitation of synchronization due to noise and stronger coupling}\label{subsec33}
\begin{figure}[ht!]
     \centering
 \includegraphics[width=1.1\textwidth]{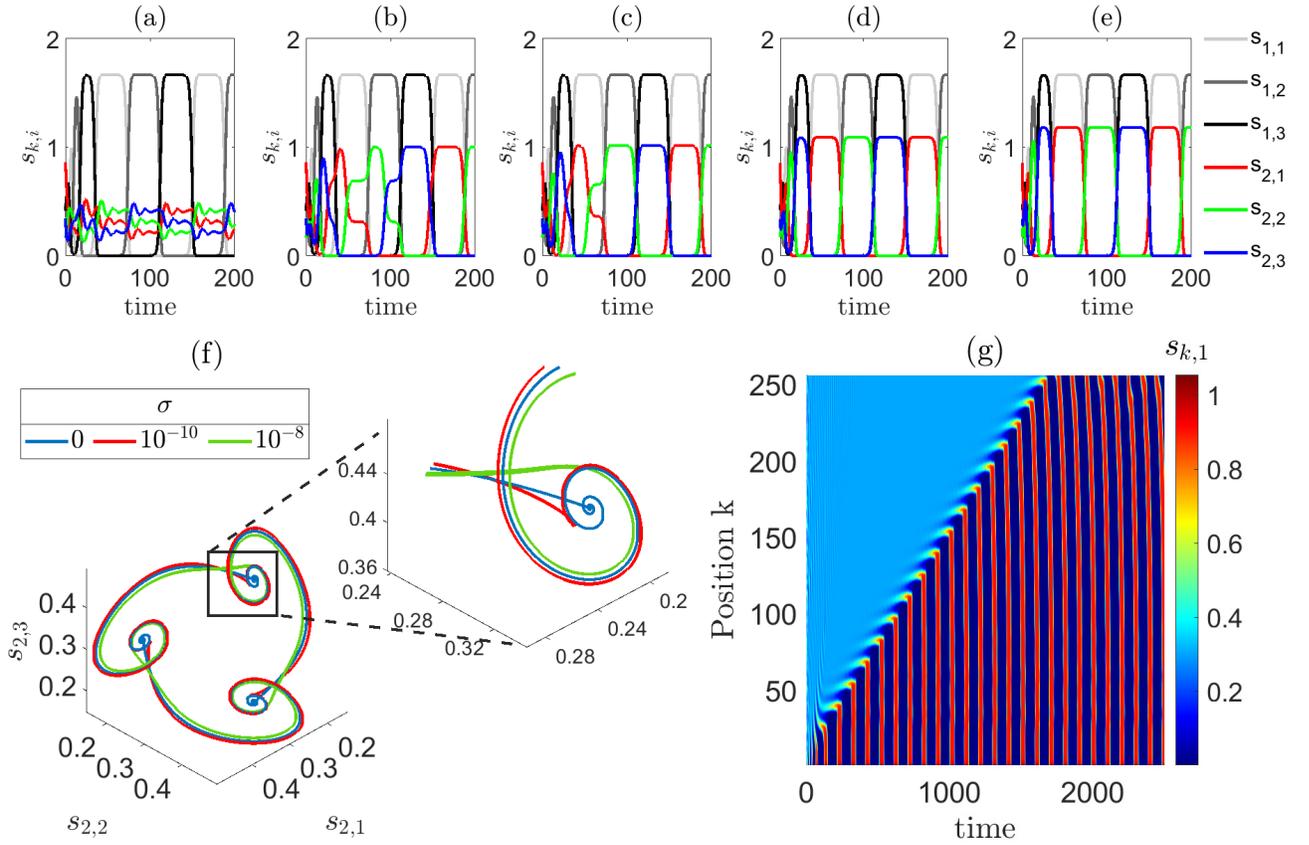}
\caption{Facilitation of synchronization due to noise and stronger coupling. Panels (a)-(e): faster entrainment with larger amplitudes due to increased forcing. Single driven unit  for $\sigma=10^{-12}$ and (a) $\delta=0.05$, (b) $\delta = 0.17$, (c) $\delta = 0.2$, (d) $\delta = 0.4$, and (e) $\delta = 0.8$. Panel (f): Zoom into the trajectory of the driven unit for $\delta=0.05$ for three values of noise $\sigma$.
 Panel (g): Spatiotemporal plot of the first item $s_{k,1}$ of each unit along a chain of unidirectionally coupled units in the presence of noise ($\sigma=10^{-12}$).
For (a)-(f)  $\gamma_P=0.6$, for  (g), $\gamma_P=0.95$ and $\delta=0.75$.
Other parameters are the same throughout: $\gamma_D=1.11$, $\rho=1$, $c=2.0$, and $e=0.2$. For further explanations see the text.}
 \label{fig:2}
  \end{figure}

The effect of noise on heteroclinic dynamics can be subtle \cite{armbruster,bakhtin,stone,ashwin2013}. However, for a certain range of noise strengths (for our other parameter choice this range is between $10^{-12}-10^{-8}$), its main effect  on the heteroclinic dynamics considered here is to prevent the slowing down of heteroclinic dynamics in the neighborhood of the saddles \footnote{For stronger noise of the order of $10^{-3}$ it can act as a control parameter as we shall see in section~\ref{subsec44}}. Without noise, we have analyzed the proliferation of possible equilibria in Appendix B.1 and Table~\ref{table:2}.
If we add noise to the dynamics of the pacemaker, the slowing down is prevented,  thus the forcing that is experienced by the driven units is continually changing.
These results are shown for two units in Fig.~\ref{fig:2}(a)-(f) for $\sigma=10^{-12}$.
For low values of the coupling strength, the driven unit continually switches between its three 3-items focus nodes (Fig.~\ref{fig:2}(a)), but the stronger the coupling (that is, the forcing), the faster the entrainment to the pacemaker's trajectory and the larger the amplitude, without any further visits to the 2- or 3-items equilibria after a transient. Finally in (d) and (e), after a very short transient, the driven unit only switches between its three 1S with the same frequency as the pacemaker, being  entrained to the heteroclinic dynamics of the pacemaker.

The zoom in Fig.~\ref{fig:2}(f) clearly shows how an increasing noise strength avoids to reach the vicinity of the saddle equilibrium. While the blue curve ($\sigma=0$) seems to spiral into the saddle before it escapes to the next one, the red curve ($\sigma=10^{-10}$) makes half a revolution of the saddle before it turns off, while the green one ($\sigma=10^{-8}$) surrounds the saddle at the largest distance, resolving the vicinity of the saddle the least.

A similar effect is seen for larger systems with $N > 2$ units. In Fig.~\ref{fig:2}(g), we show the results for a chain of unidirectionally coupled 1L-units in the presence of noise, applied to the pacemaker.
Without noise, the pacemaker would be  able to entrain around 100 units before its dynamics slows down, with noise it entrains all 256 units (Fig.~\ref{fig:2}(g)). Similar is the result for a small amount of back-coupling along the chain \cite{supp}, no de-entrainment of the units after a peak of entrainment is reached.

Thus, in summary, both stronger coupling and some amount of noise facilitate the entrainment to the pacemaker dynamics. Also a small back-coupling to the pacemaker  can have a similar effect as noise in supporting collective heteroclinic oscillations.\\

In view of brain dynamics, our result illustrates how a constructive role of noise comes about.  The large number of bifurcation points in the phase space of coupled heteroclinic networks may naively suggest that much fine-tuning of the parameters is needed to achieve a state of reproducible partial synchronization. When noise impedes  a detailed exploration of the phase space in all its possible bifurcation points, the driven units easily follow a pacemaker, and as ``noise'' is ubiquitous, this is a convenient way of facilitating synchronization by exploiting the noise.

In our modeling we do not take into account the thermodynamic cost for achieving precision of the entrainment. In \cite{barato} the relation between a stochastic Turing pattern in a Brusselator model and the thermodynamic cost have been addressed. According to the results of \cite{barato}, the precision of the pattern is maximized for an intermediate thermodynamic cost that is paid for suppressing fluctuations; this means that larger fluctuations seem to have a  positive effect on the precision of the pattern. A possible explanation for this observation may be the constructive role of noise on synchronization that we see here  for heteroclinic dynamics, while the underlying mechanism is not restricted to heteroclinic dynamics: A certain substructure of the attractor space is not accessible in the presence of noise so that the trajectories of different units in phase space become more similar and synchronize more easily. They cannot explore the close vicinities of the multitude of saddles.

\section{Heteroclinic units with two levels of hierarchy along a chain}\label{sec4}
In this section we consider heteroclinic units of which each admits two levels of hierarchy so that it provides an effective description of slow oscillations modulating fast oscillations.  As we know from previous work \cite{max1}, for an appropriate choice of  rates and other parameters, the trajectory performs an LHC during which it approaches subsequently three SHCs between three 1-item saddle equilibria.
In  section \ref{subsec41} we want to explore first with two, then with $N$ units, whether a 2L-unit can entrain one or more CE-units towards this intricate motion with two inherent time scales. As we shall see, the answer is positive.
Vice versa one may wonder whether a single unit in a resting state is also able to stop the intricate heteroclinic motion of $N-1$ $2L$-units. This is addressed in section \ref{subsec42}.

Moreover, we had emphasized that information may be encoded in the temporal sequence of visited saddles of the hierarchical heteroclinic network, corresponding to a specific path in our heteroclinic network. In our realization, this path is enforced by the choice of rates which determines the eigenvalues of the Jacobian in a way that one out of two directions for leaving the vicinity of the saddles is preferred. In section \ref{subsec43} we give an example for a different selected path  as a result of some back-coupling to the pacemaker and the choice of initial conditions under the same choice of rates.
In section \ref{subsec44} we indicate how a sufficiently large strength of noise can control the state of the system.

\subsection{Entrainment towards slow oscillations of fast oscillations}\label{subsec41}
\begin{figure}[ht!]
     \centering
 \includegraphics[width=\textwidth]{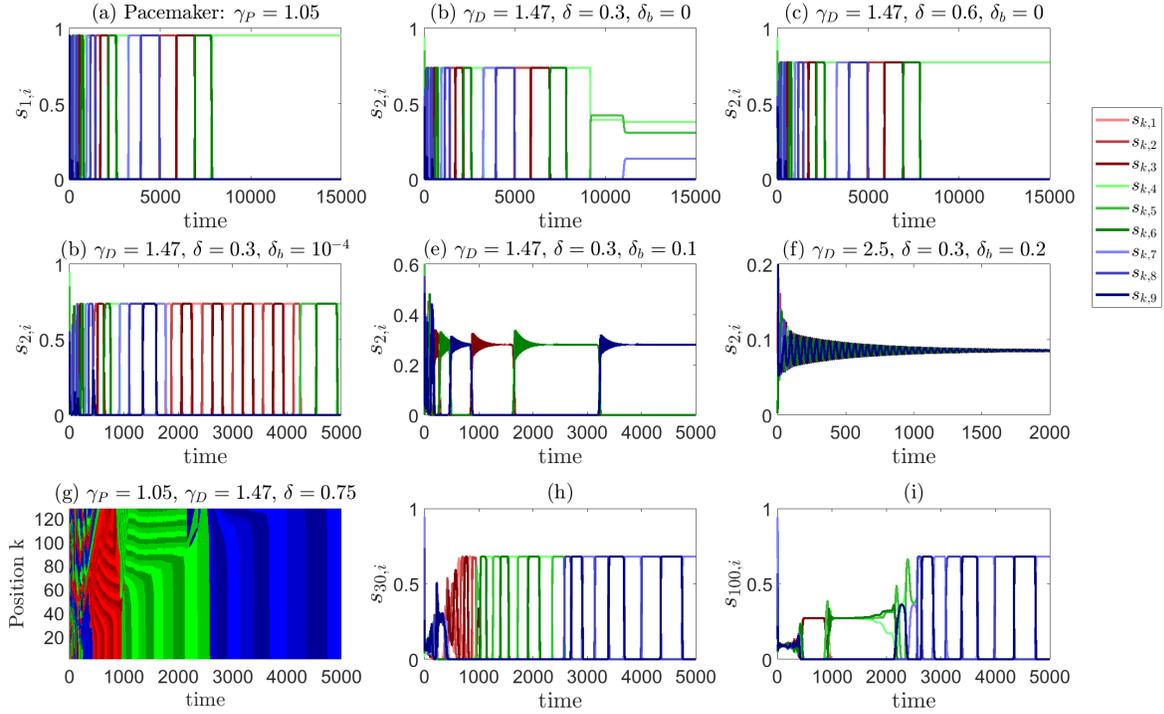}
\caption{Entrainment of the CE-units  towards two levels of hierarchy with slow oscillations of fast oscillations. Panels (a)-(f) show the result for  two units, a driven one and a pacemaker, (g)-(i) for a chain of 128 unidirectionally coupled units.  (a): time-series of the nine items of the pacemaker, (b) for the driven unit for forcing $\delta=0.3$. (c): driven unit for increased forcing $\delta=0.6$. (d): including back-coupling $\delta_b=10^{-4}$, while $\delta=0.3$. (e): for further increased back-coupling to $\delta_b=0.1$, the first level of hierarchy gets eliminated. For the CE-unit deeper in the regime of coexistence equilibria ($\gamma_D=2.5$), the second hierarchy level also breaks down as shown in (f) for $\delta_b=0.2$. In (g),  the spatiotemporal plot shows the dominant items at each site (color codes the items) of the chain with $\delta=0.75$, (h) time series of the $30^{th}$ unit, (i) for the $100^{th}$ unit. The fixed parameters for all figures are  $r=1.25$, $e=0.2$, $f=0.3$, $c=d=2$, and $\rho=1$. For further explanations see the text.
}
 \label{fig:3}
  \end{figure}

{\bf N=2 units: a pacemaker and a driven unit.} For the pacemaker we choose the 2L-unit with $\gamma_P<\gamma_c$, for the driven unit a CE-unit with $\gamma_D>\gamma_g$ ($\gamma_g$ the bifurcation point of the second Hopf bifurcation), unidirectional coupling and sufficiently strong forcing $\delta$, so that  the pacemaker (Fig.~\ref{fig:3}(a)) can entrain the driven unit towards slow oscillations of fast oscillations.  An example is shown in Fig.~\ref{fig:3}(b) where the driven unit entrains to the dynamics of the pacemaker till around ten thousand time units. The three items of the first (second, third) SHC of both the pacemaker and the driven unit are shown in shades of red (green, blue), respectively. Beyond $\text{time}=10000$, the trajectory of the pacemaker spends so much time near the first item saddle of the second SHC that the trajectory of the driven unit has enough time to go through a sequence of equilibria  upon its pursuit to reach the stable equilibrium. This stable equilibrium  would be approached in the equivalent system of a unit driven by a constant force. The dynamics is the same as we discussed for the systems with a single hierarchy level. However, if the pacemaker applies a stronger forcing to the driven unit, its trajectory  is pushed towards the corresponding saddle  and hence  spends a long time together with the pacemaker as shown in Fig.~\ref{fig:3}(c).

As an alternative to applying noise to prevent the slowing down, a small amount of back coupling ($\delta_b$) to the pacemaker has the same effect, as we saw already for the systems with a single hierarchy \cite{supp}. Both units then also entrain to heteroclinic dynamics with two levels of hierarchy. Also, the number of revolutions within an SHC increases with an increase in the back coupling.
We show an example in Fig.~\ref{fig:3}(d) for the driven unit. The pacemaker  has similar dynamics, with its amplitude being around 1.

If we increase $\delta_b$ beyond a certain value, the first level of hierarchy collapses, both units entrain to an HC between their corresponding 3-items saddle equilibria as shown in Fig.~\ref{fig:3}(e). If the driven unit is deeper in the CE-region, the second hierarchy level also breaks down as shown in Fig.~\ref{fig:3}(f).

{\bf $N > 2$ units without back-coupling.}  Also  a larger system of units with $N>2$ can be entrained towards slow oscillations of fast oscillations. We show the results for a unidirectionally coupled chain of 2L-units of size $N=128$ in  Fig.~\ref{fig:3}(g)-(i). The pacemaker is again placed at the edge position $k=1$. The units closer to the pacemaker get entrained first, for example, the $30^{th}$  at around 1000 time steps (Fig.~\ref{fig:3}(h)), the $100^{th}$ unit  at around 2000 time steps, again there is a pronounced transient time for the entrainment process to propagate along the chain as revolutions along the cycles of successive units take time. The different shades in panel (g) beyond the three main colors indicate the presence of two hierarchy levels.

\subsection{Inverting the role of pacemaker and driven units}\label{subsec42} So far our main interest was in the option of whether a single unit with oscillatory motion in the vicinity of heteroclinic orbits can entrain other units in a resting state, therefore acting as a pacemaker. The directed coupling, acting more strongly in forward than back-direction and the position of the pacemaker at the edge of a chain were obviously important when a single unit should entrain a whole set of driven units. (For bidirectional (diffusive) coupling, as we consider in section \ref{sec5}, either more pacemakers are needed, or again an asymmetry in the allocation within a whole disc.) Thus it is natural to ask what happens in the inverse case, a single CE-unit at the edge of a chain, directed coupling as before, but 1L or 2L-units as the driven units along the chain. Is a single CE-unit able to stop the heteroclinic oscillations of the entire set of driven units? The answer is positive. Thus it is the combination of the location with the directed coupling and a certain parameter set that enables a single heteroclinic unit to control a larger set of other such units and entrain this set to its own dynamics.

\subsection{Selection of different trajectories in the hierarchical attractor landscape and their synchronization via entrainment }\label{subsec43}
\begin{figure}
     \centering
      {\includegraphics[width = 0.4 \textwidth, height=2in]{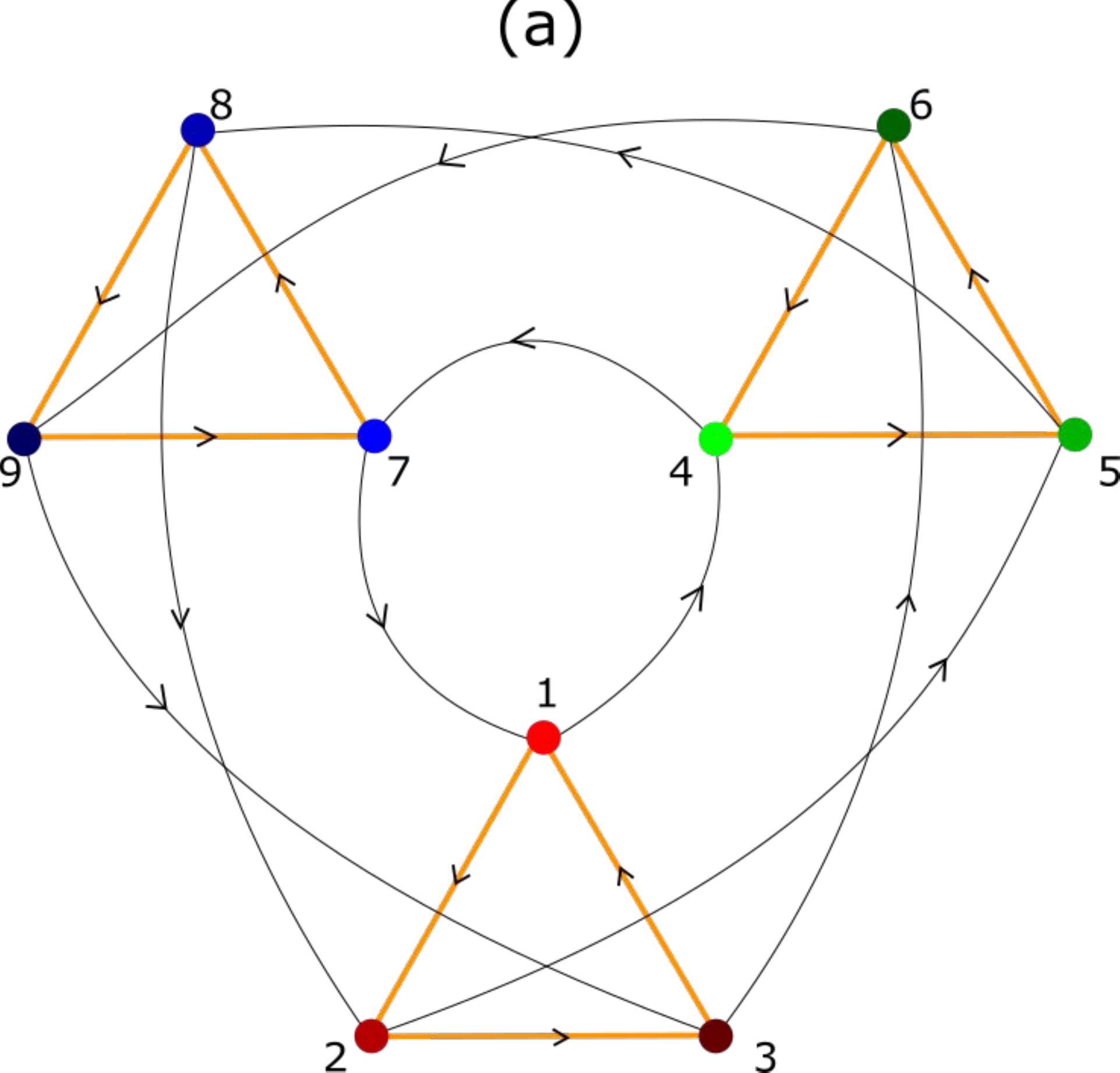}}\hfill
     {\includegraphics[width = 0.4\textwidth, height=2in]{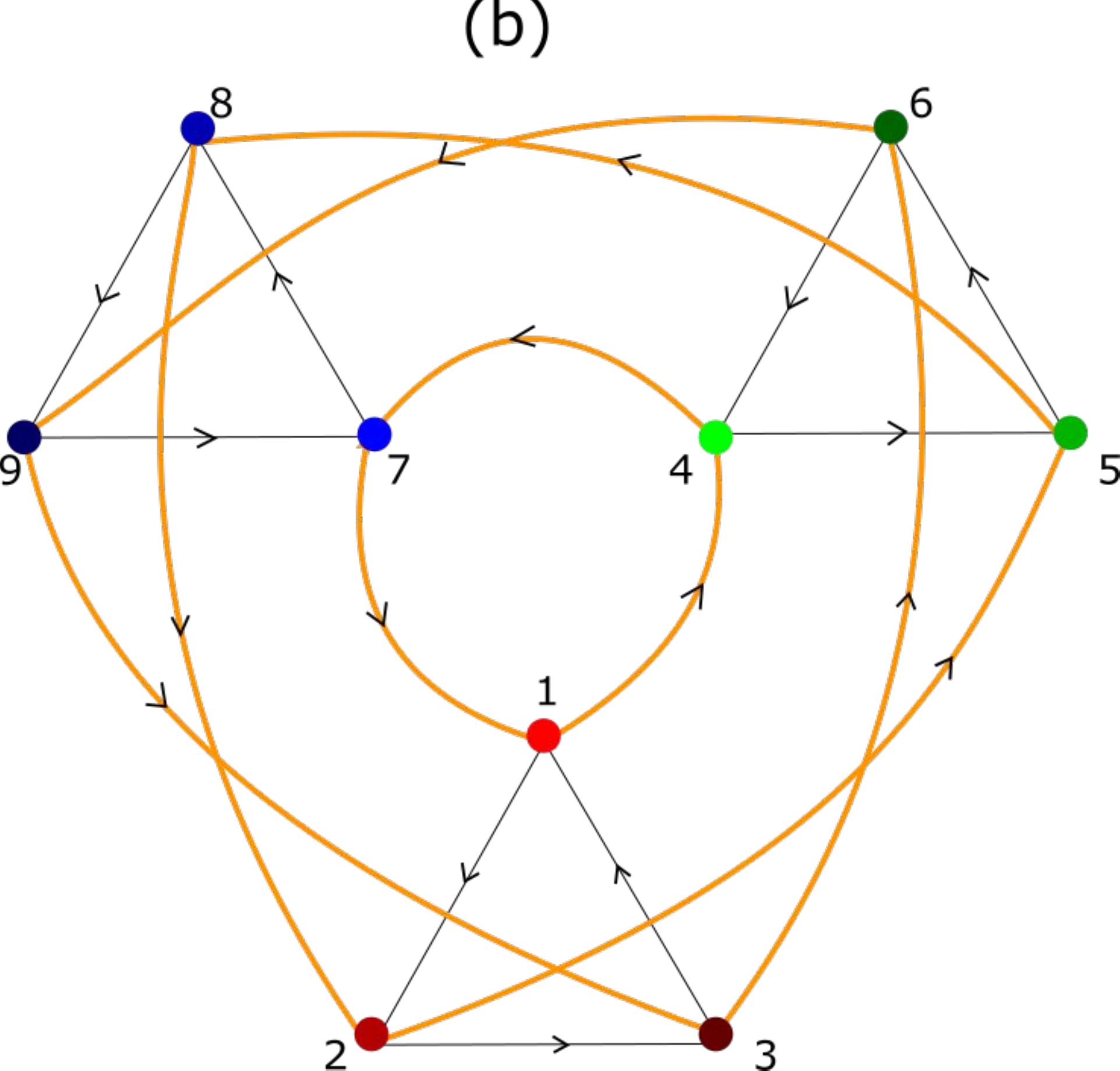}} \\
  {\includegraphics[width = 1.1\textwidth]{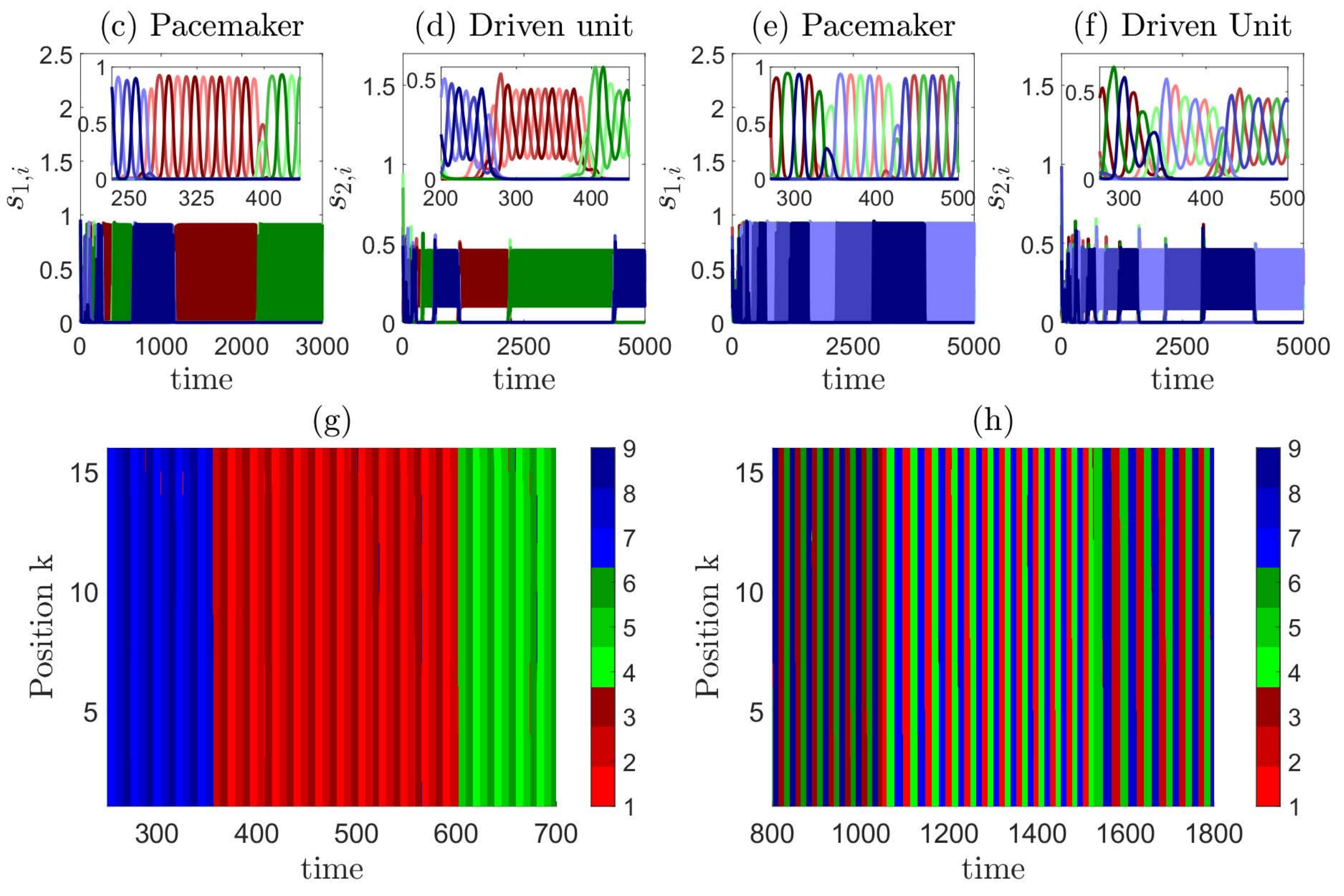}}
     \caption{Two coexisting  spatiotemporal patterns, generated from two different paths within the heteroclinic network. (a) The path approaches the saddles $(1,2,3)\rightarrow (6,4,5) \rightarrow (8,9,7)$, (b) $(1,4,7)\rightarrow (8,2,5) \rightarrow (6,9,3)$. In both cases orange lines denote SHCs, black lines connections between SHCs in phase space along possible LHCs. For the system of two units, the time series of both units for the dynamics along the first path are shown in (c)-(d), along the second path in (e) and (f). The values of the coupling strength are $\delta=0.1$, $\delta_b=0.001$.
     Panels (g) and (h): Differing spatiotemporal patterns of 16 synchronized units entrained by one pacemaker to hierarchical heteroclinic motion along either of the two paths as indicated in (a) and (b), respectively, color coded are the items which are dominant, $\delta=1.5$. Other parameters are $\gamma_P=1.05$, $\gamma_D=1.47$, $\delta_P=0.01$, $r=1.25$, $e=0.2$, $f=0.3$, $c=d=2$, and $\rho=1$.}
    \label{fig:4}
\end{figure}

As known from previous results \cite{max1,max2} and the last section, the predation rates of a 2L-unit can lead to a trajectory in the attractor space that is schematically indicated in Fig.~\ref{fig:4}(a). Corresponding to (a), a possible sequence of saddles, whose neighborhood the trajectory may visit subsequently, would be
$(1,2,3)\rightarrow (6,4,5)\rightarrow (8,9,7)$ with cyclic repetition; an alternative sequence according to Fig.~\ref{fig:4}(b) would be $(1,4,7)\rightarrow (8,2,5)\rightarrow (6,9,3)$, cyclically repeated. \footnote{Note that the visualization indicates that the terminology in terms of small or large heteroclinic cycles does not reflect the metric in phase space, but merely corresponds to the graphical representation, here adapted to (a), while the heteroclinic connections along the ``large'' cycle in (b), that is along $(7,8),(5,6),(3,1)$, are graphically represented as ``short'' connections.}

{\bf N=2 units.} We consider again two units, the pacemaker with $\gamma_P=1.05$, the driven unit with $\gamma_D=1.47$, $\delta=0.1$ and $\delta_b=0.001$.

For most initial conditions, the trajectories of both units follow  an LHC of small {\it limit} cycles (SLCs) along the first path as shown in  Fig.~\ref{fig:4}(c)-(d). \footnote{Note that here we call the cyclic behavior limit cycles rather than heteroclinic cycles, as the amplitudes are clearly off the heteroclinic contour.} However, for some initial conditions, the trajectories follow an LHC of SLCs along the second path as shown in Fig.~\ref{fig:4}(e)-(f).
The LHC of SLCs along the second path is  seen for a range of back-couplings which also depends on the choice of the initial conditions. For example, for the initial conditions chosen in (e)-(f), the LHC of SLCs along the second path is seen in the range of $\delta_b\in [0.00005,0.02379]$. Right outside this interval, the trajectories instead follow the first path.

{\bf $N > 2$ units.} The  choice of the second path in the hierarchical heteroclinic network is also seen for larger systems with $N > 2$ units. Results are obtained for a system of unidirectionally coupled units along a ring of size $N=16$ with $\delta=1.5$ for all links apart from the last one which  closes the chain to a ring. Since this coupling amounts to a back-coupling to the pacemaker, it is chosen less strong as $\delta_P=0.01$.
The time evolutions of the trajectories for the pacemaker and the driven units look similar to those in Fig.~\ref{fig:4}(c)-(f) for two sets of initial conditions, one choosing the first path and the other choosing the second path.
The amplitudes of entrained units are smaller than those of the pacemaker, both when we have to deal with a single driven unit (in case of $N=2$) and also as a function of the distance from the pacemaker (visible for larger $N$, not displayed).
The attractor corresponding to the second path seems  to have a smaller basin of attraction.
We ran the system of rings of heteroclinic units for $100$ initial conditions, of which  three  led to the path according to (b).
The corresponding spatiotemporal plots are shown in Fig.~\ref{fig:4}(g) and (h), respectively, where we plot the dominant item at each site as a function of time.

In view of possibly shared features with brain dynamics, Fig.~\ref{fig:4} (g) and (h) display synchronized units along the entire chain for both paths in the attractor space, and it is the temporal order of switches between subsequently visited saddles  that differs between both paths. This means that the corresponding sequence of dominant information items (coded in respective colors) is different, as visible in the spatiotemporal representation. Since the very temporal order of collective patterns is said to encode information \cite{odor1,odor2}, our results indicate possibly many more options for encoding information by choosing different paths in a given  heteroclinic network. In our implementation, different paths are selected as a result of the rates, the (back)-coupling and the initial conditions (in case of multistability).  In the brain, the counterpart of different paths  may be selected by external input.

In summary, depending on the selected path in the phase space of hierarchical heteroclinic motion, the encoded information string varies.  When it is transferred from the pacemaker along the chain, as a result of entrainment different hierarchical processes are synchronized and going on in parallel at different sites of the grid. This looks like a mathematical realization of the typical features of brain dynamics as conjectured in \cite{grossberg}.

\subsection{Noise acting as a control parameter}\label{subsec44}
As we know from a single heteroclinic unit with two hierarchy levels  \cite{max1}, noise can act similarly to the bifurcation parameter $\gamma$. It  acts as a control parameter to drive the dynamics from two levels of hierarchy to one level with an effective limit cycle between three 3-items saddle equilibria, and further to global coexistence of all items.
Similarly here, for example for a chain of $N=64$ units, with the pacemaker located at the edge position with $\gamma_P=1.05$ (from the 2L-regime) and all other units with $\gamma_D=1.47$ (from the CE-regime), we pursue the entrainment as a function of an increasing noise strength. Full entrainment to two levels of oscillatory motion is possible for low noise of the order of $\sigma=10^{-12}$; noise then prevents the slowing down of the pacemaker and supports entrainment to two levels  for driven units in the neighborhood of the pacemaker. However, at a more remote distance to the pacemaker, entrainment is achieved  only to one hierarchy level, as the influence of the pacemaker dampens with the distance. Stronger noise reduces the levels of the pacemaker from two to one and along with that, the levels of the driven units to a single one at all sites. A further increase of noise then drives the pacemaker to the resting state together with the driven units, whose parameters favor this state anyway. The corresponding figures and more details are presented in \cite{supp}.
\section{Maintenance of synchronization in the presence of resting-state units on a two-dimensional grid}\label{sec5}
Next we consider a two-dimensional spatial grid and assign to  each of its sites a hierarchical unit which in principle allows two hierarchy levels in time scale, depending on the choice of individual parameters. In earlier studies of one of us \cite{max3} on nested spirals in coupled heteroclinic networks, the parameters of individual units were chosen homogeneously from the HC-regime. Starting from random initial conditions, the synchronization to the heteroclinic motion of effectively a single unit was achieved by diffusive coupling, but the individual parameters were chosen to be the same and from the `right' regime in view of oscillations.
Here, for comparison with these results, the coupling is again diffusive (that is bidirectional rather than unidirectional as before), but the parameters of only a fraction $1-p_{d}$ of these units is chosen in the 2L-regime, while the fraction $p_d$  of  the remaining units is in the CE-regime. \footnote{The latter ones may be considered as defects in the sense that  they individually do not oscillate. In an analogous analysis of systems of coupled nonlinear oscillators \cite{daido1,daido2}, such defects are considered as a possible result of ageing, here we do not specify the origin of parameters from the CE-regime, they may correspond to default values in a kind of resting state, being disposed for getting entrained into oscillatory motion (rather than being aged).}

For a given set of parameters and initial conditions, we are interested in determining the fraction of defects that is tolerable in order not to destroy collective hierarchical heteroclinic motion with two time scales.
Thus we diffusively couple the units on the grid with the strength $\delta$. The dynamics of the system then follows the equation:
\begin{equation}
 \partial_t s_{i} = \delta \nabla_x^2 s_i + \rho s_i - \gamma_x s_i^2 -\sum_{j\neq i} A_{i,j} s_i s_j + \sigma |\xi_i(t)|,
\end{equation}
where $x$ indicates the location on the grid with $L\times L$ sites (here $L=128$), $\gamma_x$ is the decay rate at the spatial location $x$.
We again use $\gamma$ as the bifurcation parameter of an individual unit and fix the other  parameters at $r=1.25$, $e=0.2$, $f=0.4$, and $c=d=2$. For this choice of parameters and $\gamma_x<1.2$, the unit at location $x$ has two levels of hierarchy in the absence of coupling. For $\gamma_x>1.5$, the unit is in global coexistence of all nine items.  In the following we consider two assignments of defective units to the grid: The first one amounts to a uniform random distribution (Fig.~1(c)), the second one to all units being defective outside a little disc (Fig.~1(d)), to which the pacemakers are constrained in space.

\subsection{Randomly distributed defects}
For the first case we assign to each site a value from a uniform distribution in the interval $[0,1]$ resulting in defects at sites with values less than or equal to $p_d$. The second case will be discussed subsequently.
The decay rate of all possible pacemaker units is taken as $\gamma_x=0.8$, for the defective units we choose a uniform distribution from the interval $[1.5,1.6]$ of the CE-regime.
\begin{figure}
    \centering
    \subfloat[]{\includegraphics[width = 0.33 \textwidth]{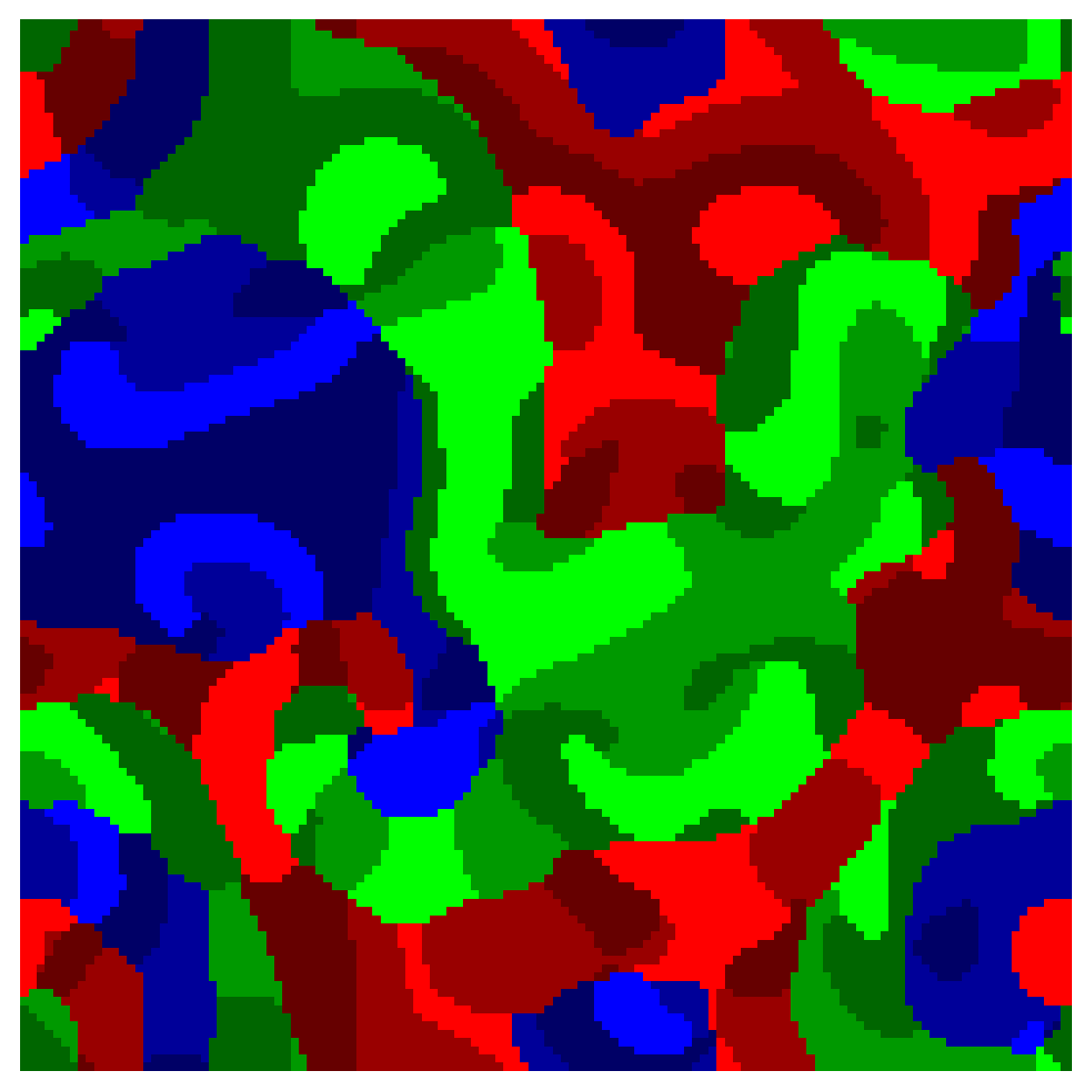}}
    \subfloat[]{\includegraphics[width = 0.33 \textwidth]{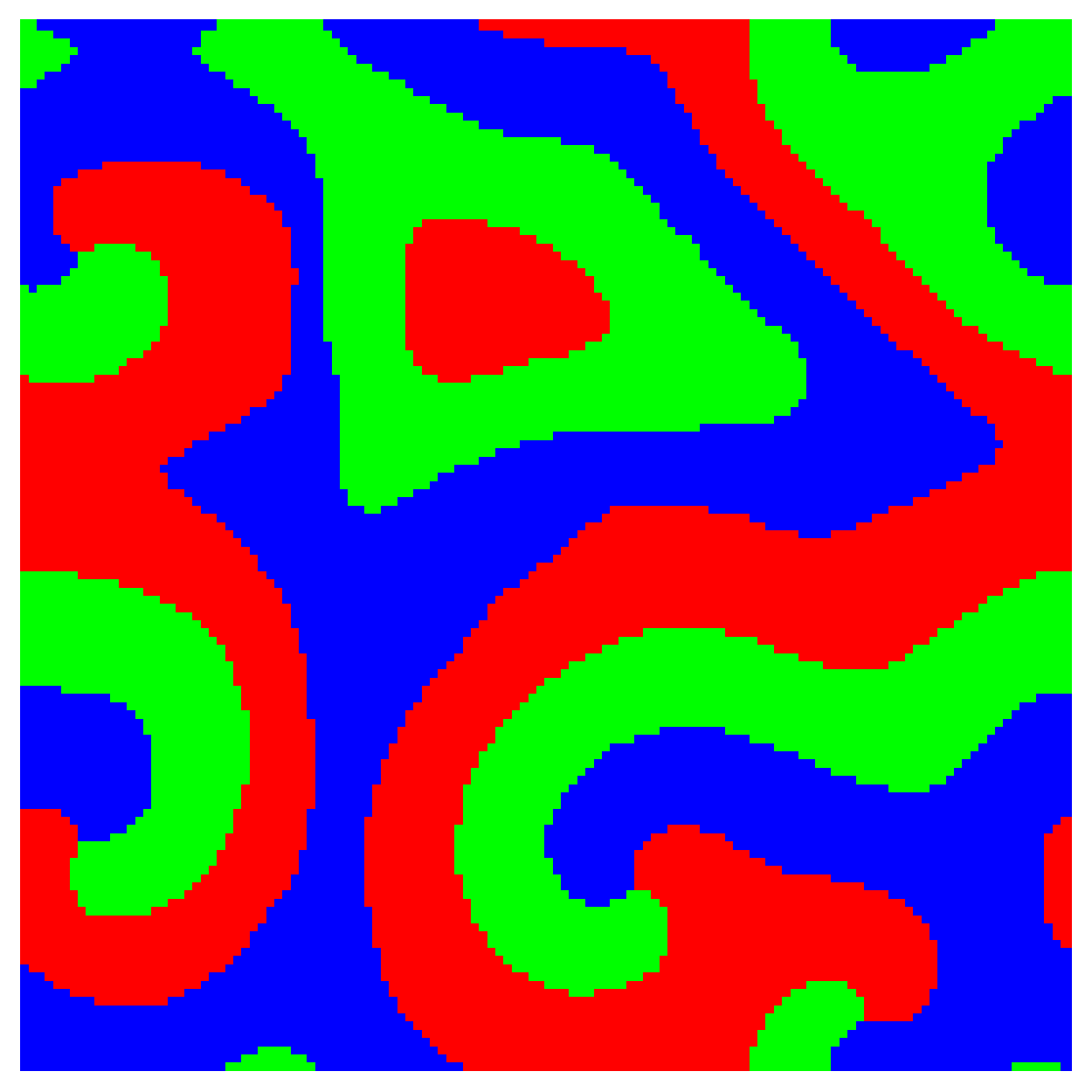}}
    \subfloat[]{\includegraphics[width = 0.33 \textwidth]{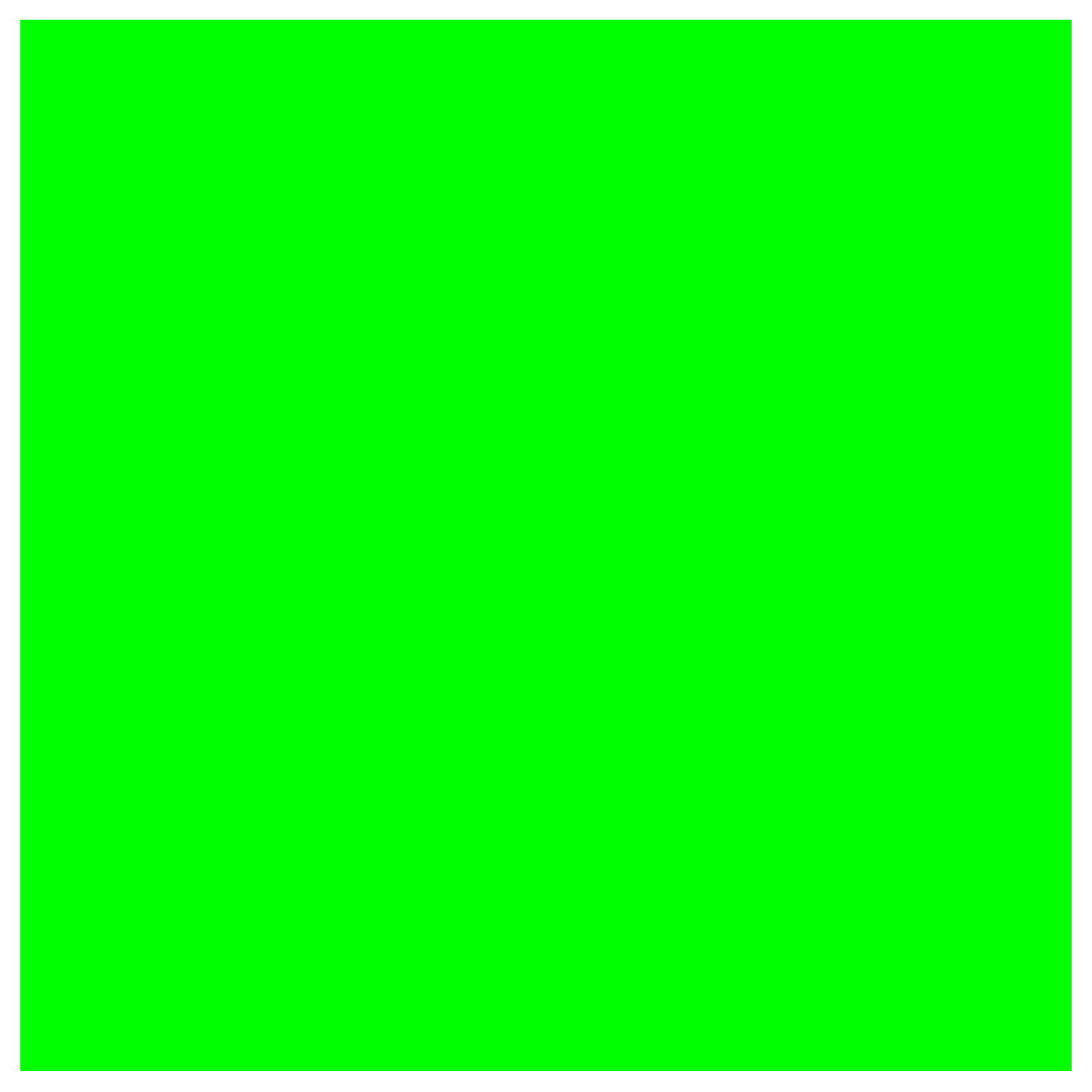}}\\
    \subfloat[]{\includegraphics[width = 0.3 \textwidth]{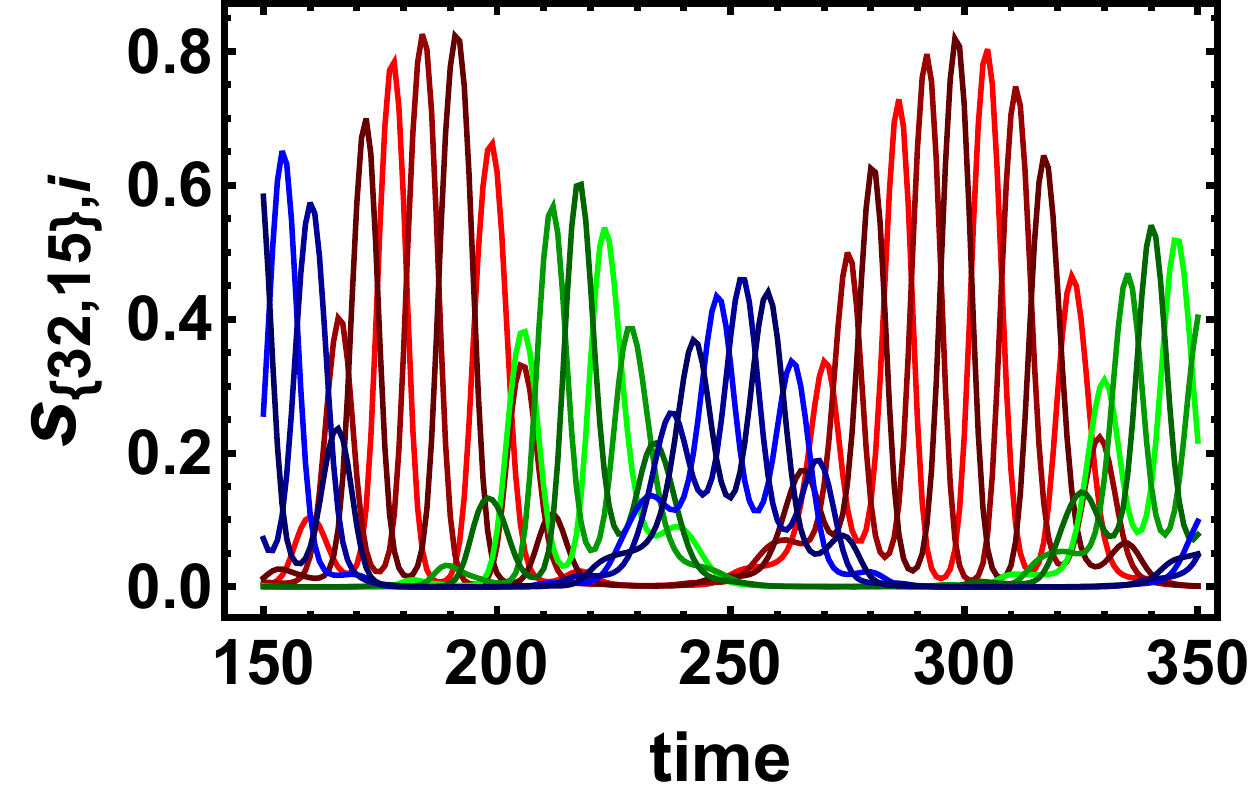}}
    \subfloat[]{\includegraphics[width = 0.3 \textwidth]{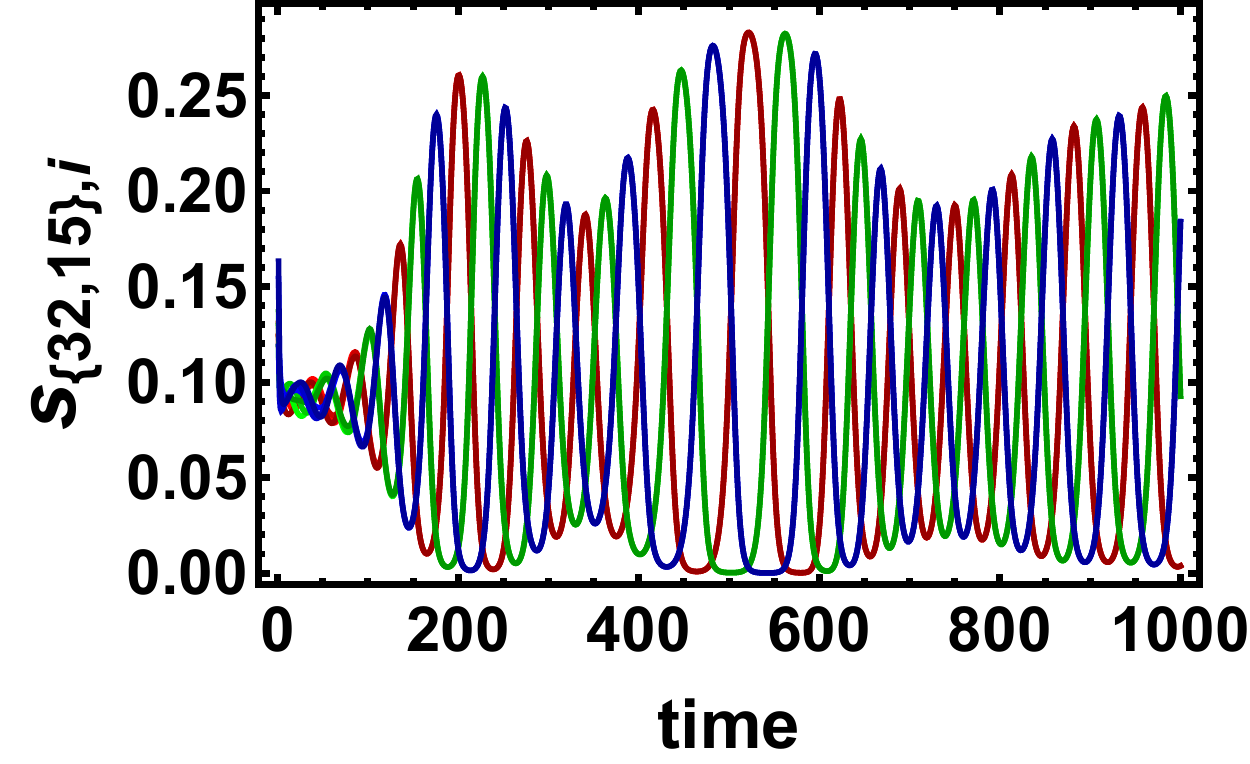}}
    \subfloat[]{\includegraphics[width = 0.3 \textwidth]{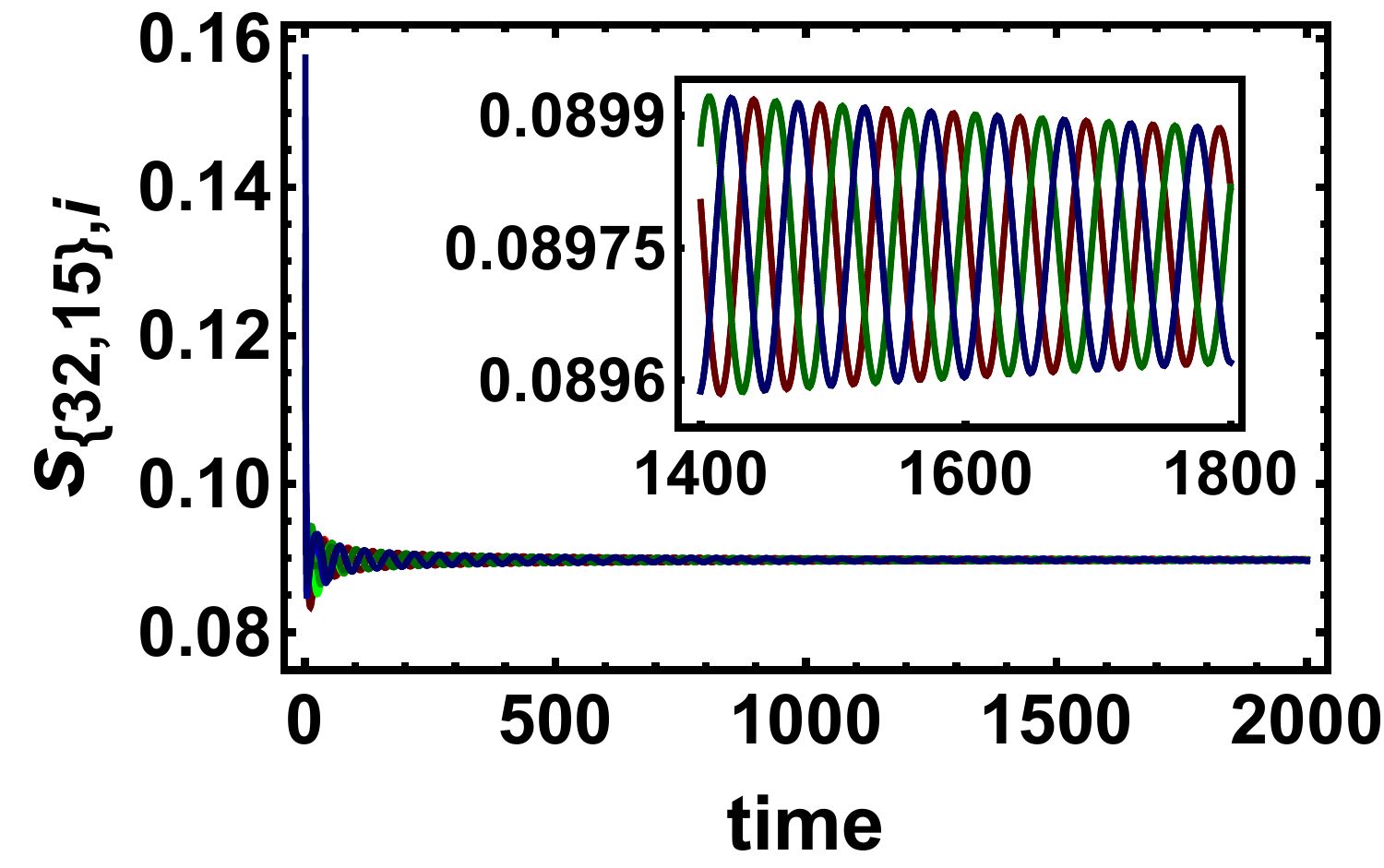}}
    {\includegraphics[width = 0.1 \textwidth]{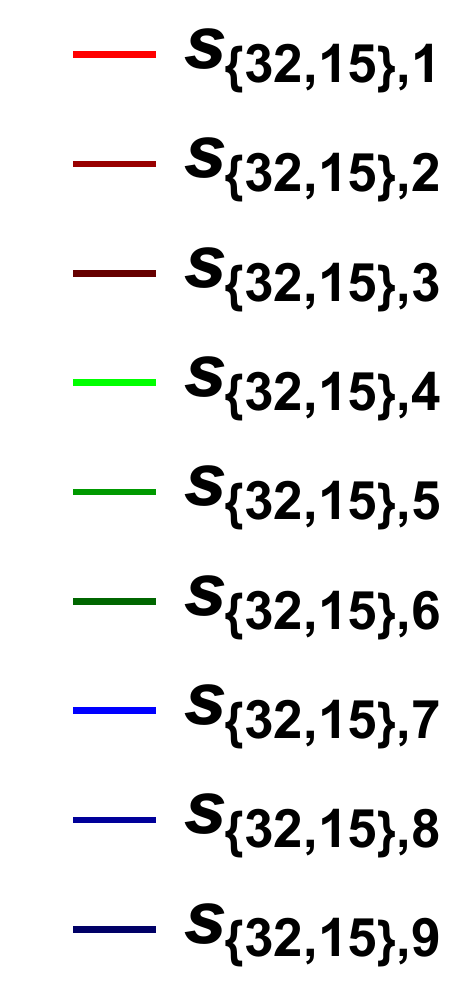}}
     \caption{Spatial entrainment of heteroclinic units on a two-dimensional grid of size $128\times 128$  for three values of the fraction of randomly distributed defects  $p_d$ . The parameter $\gamma_x$ of 2L-units is fixed at $0.8$, $\gamma_x$ of  CE-units are chosen from a uniform distribution in the interval $[1.5, 1.6]$. The  panels (a)-(c) show the respective pattern formation, the lower panels (d)-(f) the corresponding time series of the information items of a representative unit at $x=\{32,15\}$ with nested spirals for (a) $p_d=0.1$ at $\text{time}=1260$, LHC-spirals for (b) $p_d=0.5$ at $\text{time}=1850$, and global coexistence at each site for (c) $p_d=0.95$ at $\text{time}=1500$. A zoom in (f) shows remnants of the oscillations. Other parameters are $r=1.25$, $e=0.2$, $f=0.4$, $c=d=2$, and $\rho=1$.}
    \label{fig:5}
\end{figure}
For different values of the fraction of defects, we observe different patterns on the grid as shown in Fig.~\ref{fig:5}.
The subfigures (a)-(c) in the upper panel display the color of the dominant item at each location of the grid with species $1-3$, $4-6$, and $7-9$ in different shades of red, green, and blue, respectively.

For $p_d=0.1$,  most of the units are in the
2L-region. The pattern formation in Fig.~\ref{fig:5}(a) amounts to nested spirals on the scale of the LHCs, visible in spiral arms of different color, while different shades of a certain color represent the chasing inside an SHC. This means, at each location, two hierarchy levels are present. We show an example in Fig.~\ref{fig:5}(d) with the time series of all items of the unit at location $x=\{32,15\}$. 
The decay rate of this unit is $\gamma_{\{32,15\}}=1.531$, so it is a ``defect'', but gets entrained to join the collective hierarchical motion. We focus on such a particular unit as its decay rate does not change when we change the fraction $p_d$ in Figs.~\ref{fig:5}(b) and (c), because we have used the same random number seeds for the uniform distributions in each case.

Until around $p_d=0.4$, the dynamics on the grid is able to maintain two levels of hierarchy, visible in the corresponding time series of winnerless competition between nine items in Fig.~\ref{fig:5}(d), while for a larger fraction of defects  the lower level of hierarchy is eliminated. We observe LHC spirals in Fig.~\ref{fig:5} (b), that is, the spirals of clusters. The LHC connects the three 3-items saddle equilibria. As an example,  Fig.~\ref{fig:5} (e) gives the corresponding time series of all items of the unit at location $x=\{32,15\}$.

For an even larger number of defects, the second hierarchy level also collapses and the global coexistence state stabilizes. In Fig.~\ref{fig:5} (c) and (f), we present the results for $p_d=0.95$. Though the items of the unit shown in Fig.~\ref{fig:5} (f) are moving towards global coexistence, there is still some switching with a rather small amplitude between the three SHCs as shown in the inset.

\subsection{Pacemakers at the center of the grid}
\begin{figure}
    \begin{centering}
     \subfloat[]{\includegraphics[width = 0.33 \textwidth, height = 1.75in]{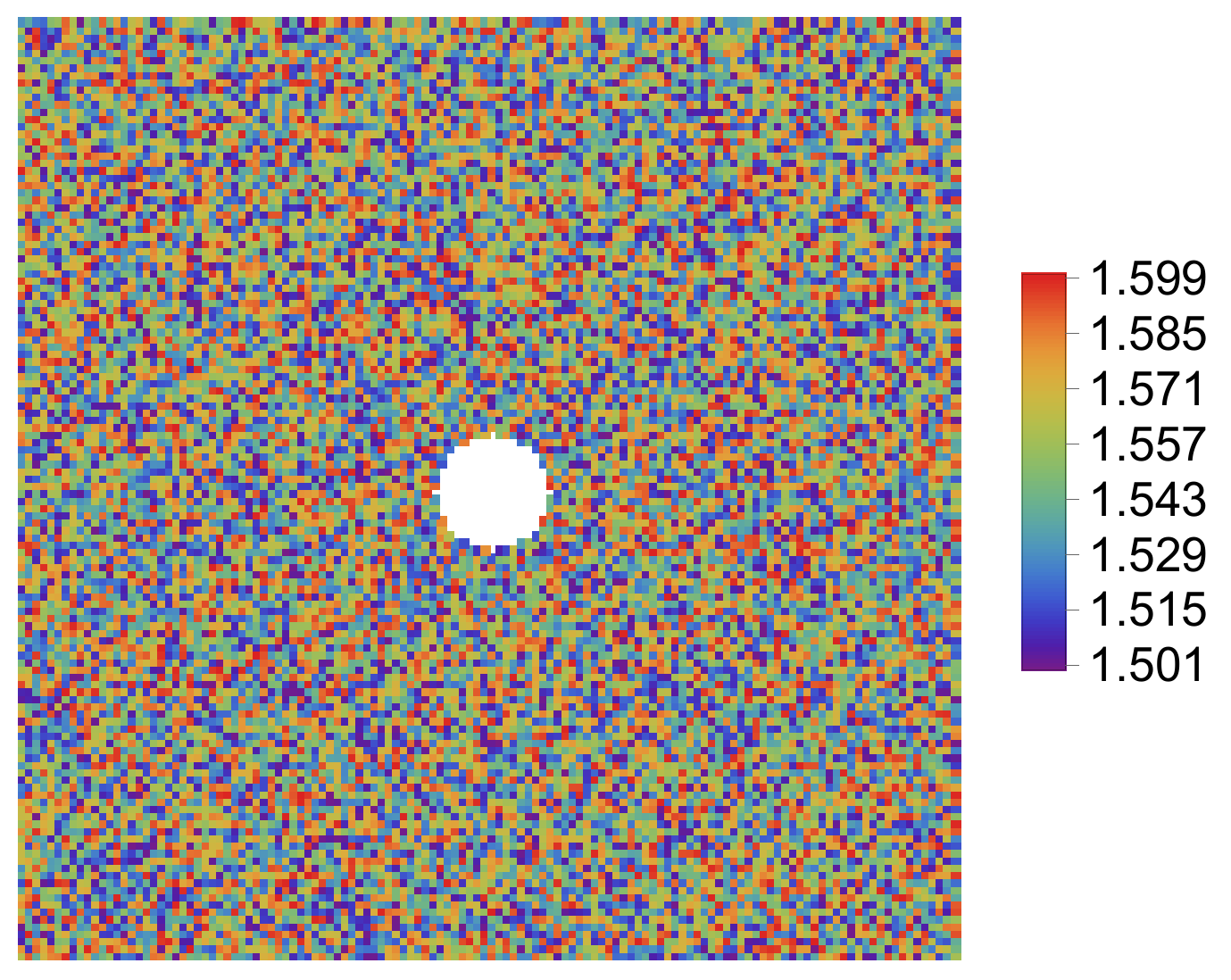}}
    \subfloat[]{\includegraphics[width = 0.3 \textwidth, height = 1.75in]{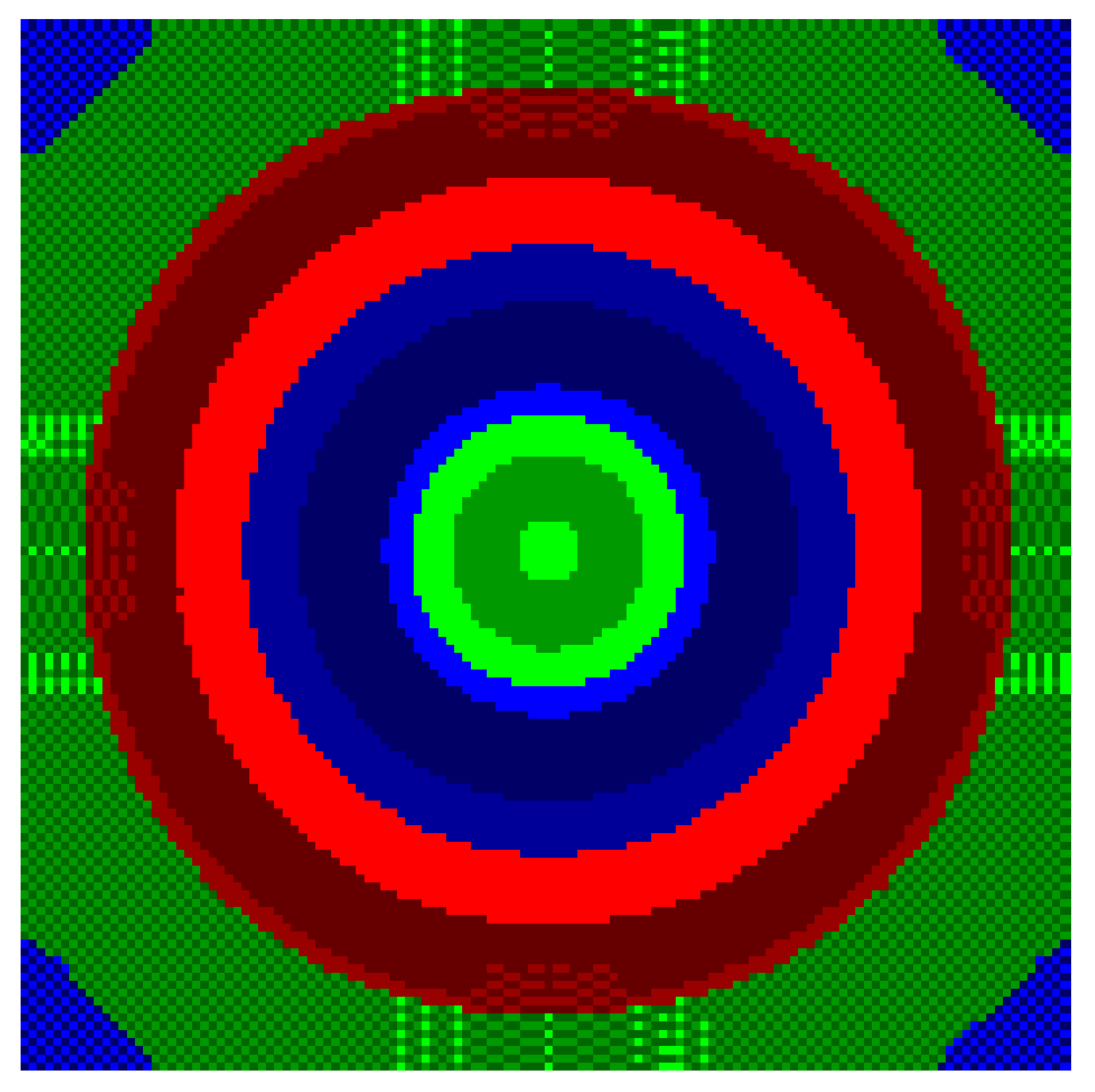}}
    \subfloat[]{\includegraphics[width = 0.33 \textwidth]{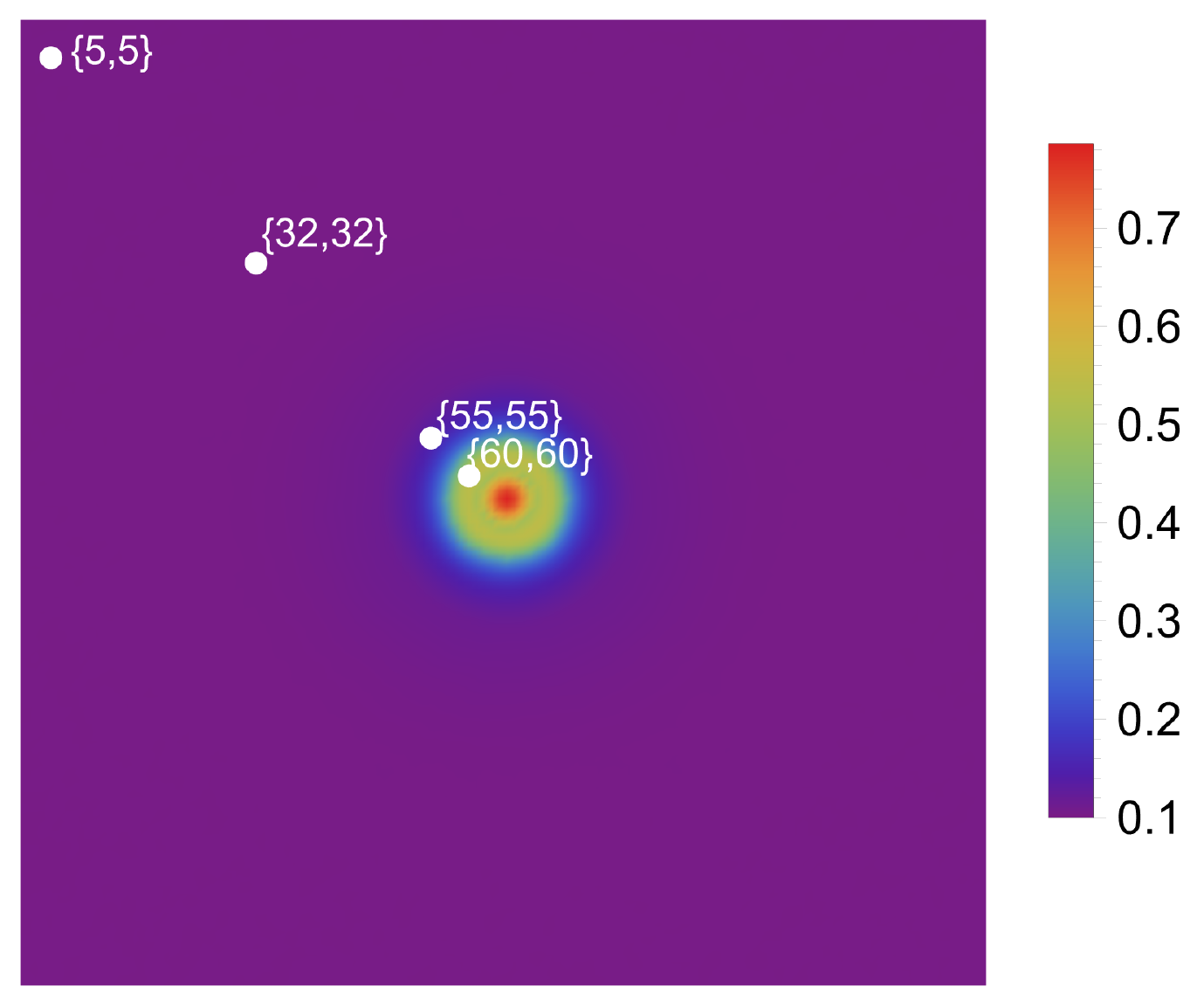}}\\
    \subfloat[]{\includegraphics[width = 0.44 \textwidth]{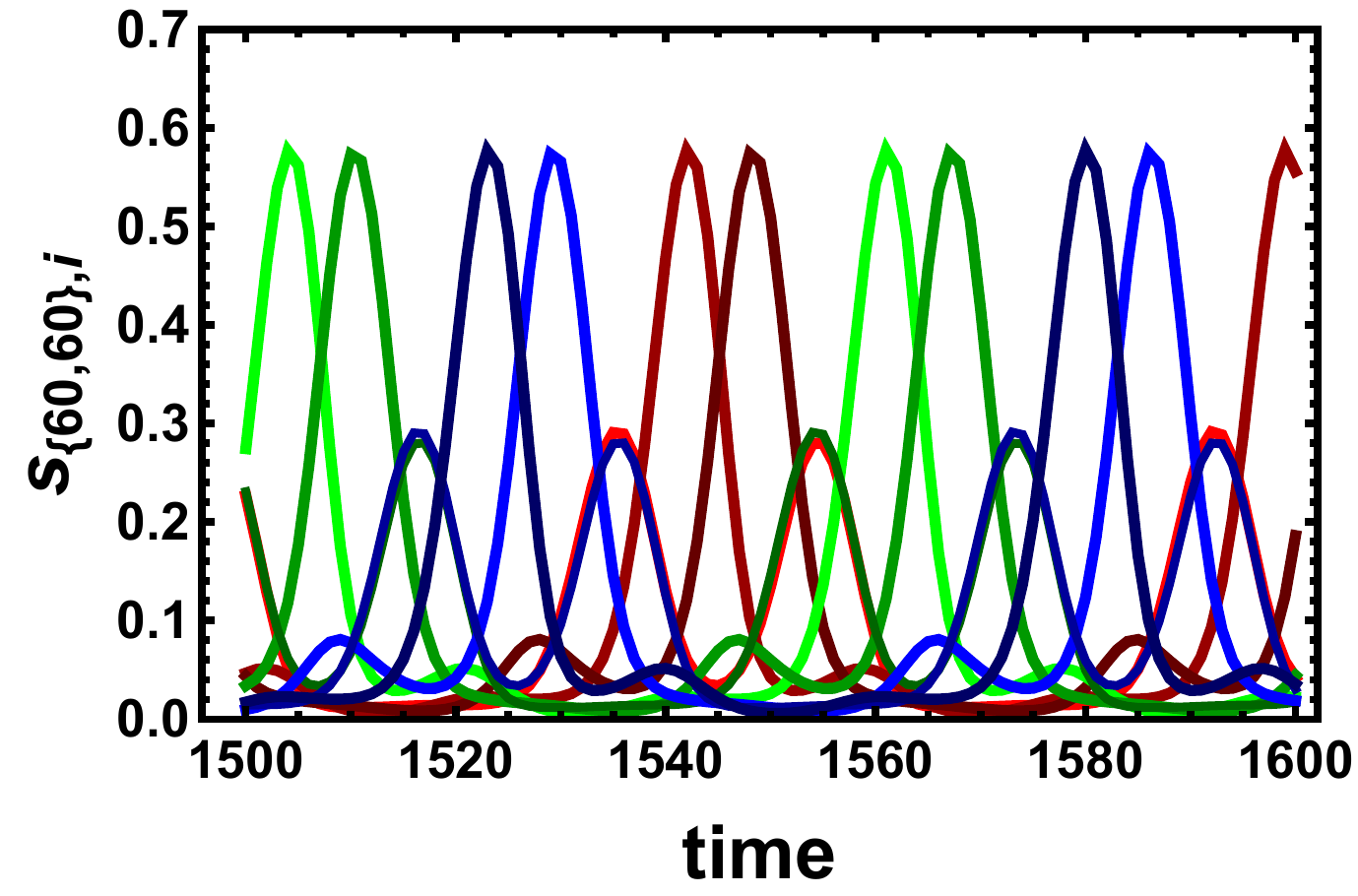}}
    \subfloat[]{\includegraphics[width = 0.44 \textwidth]{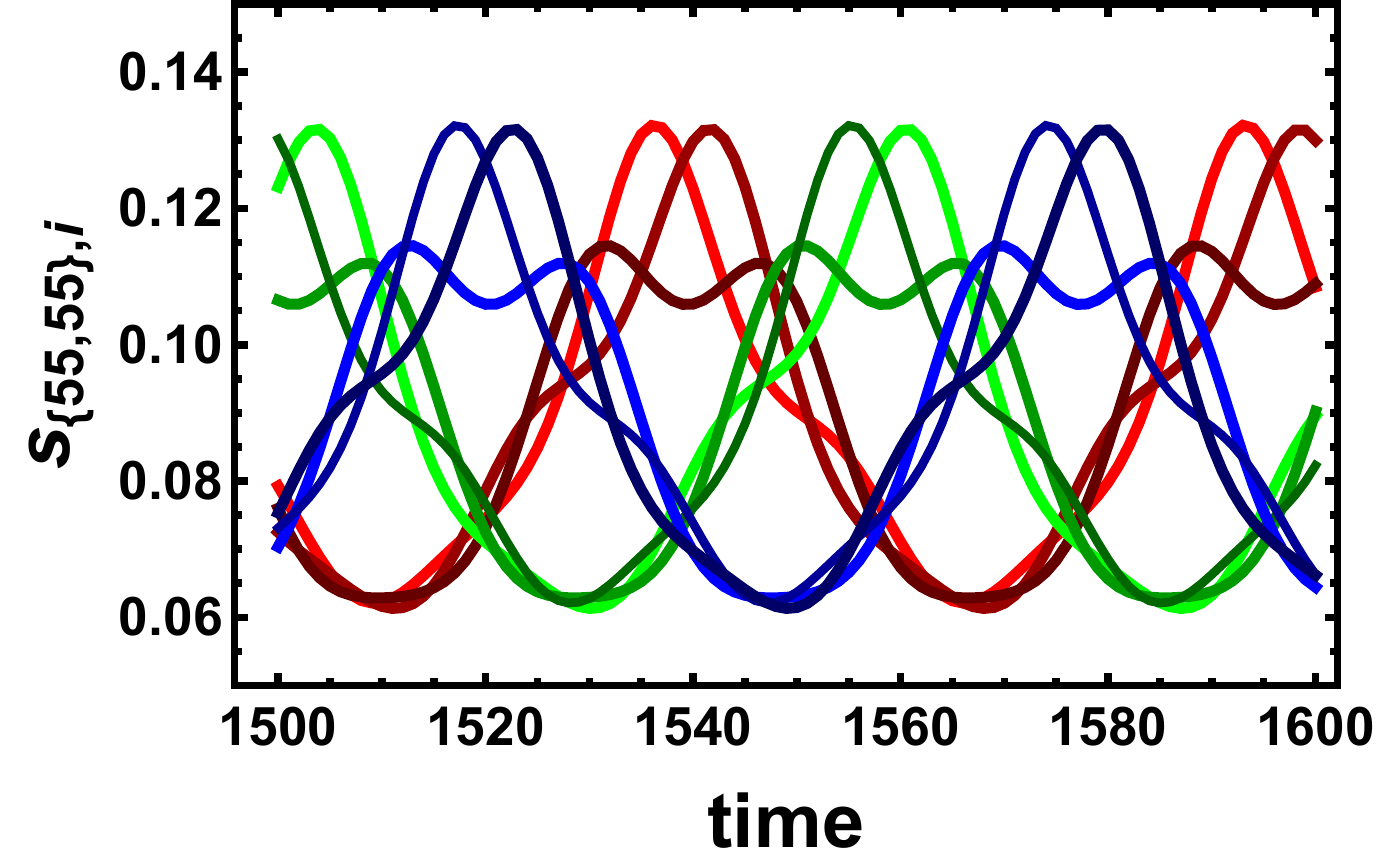}}\\
    \subfloat[]{\includegraphics[width = 0.44 \textwidth]{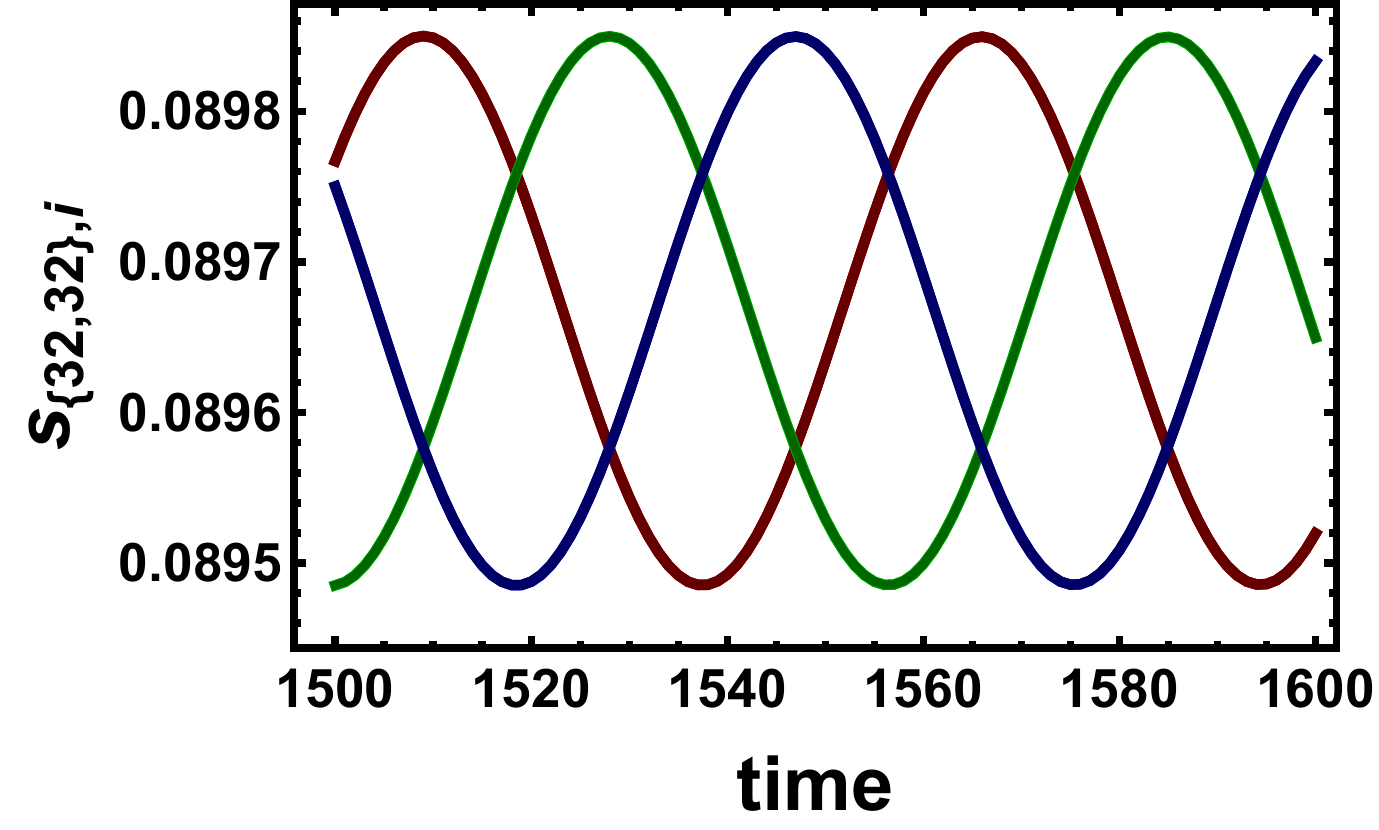}}
    \subfloat[]{\includegraphics[width = 0.44 \textwidth]{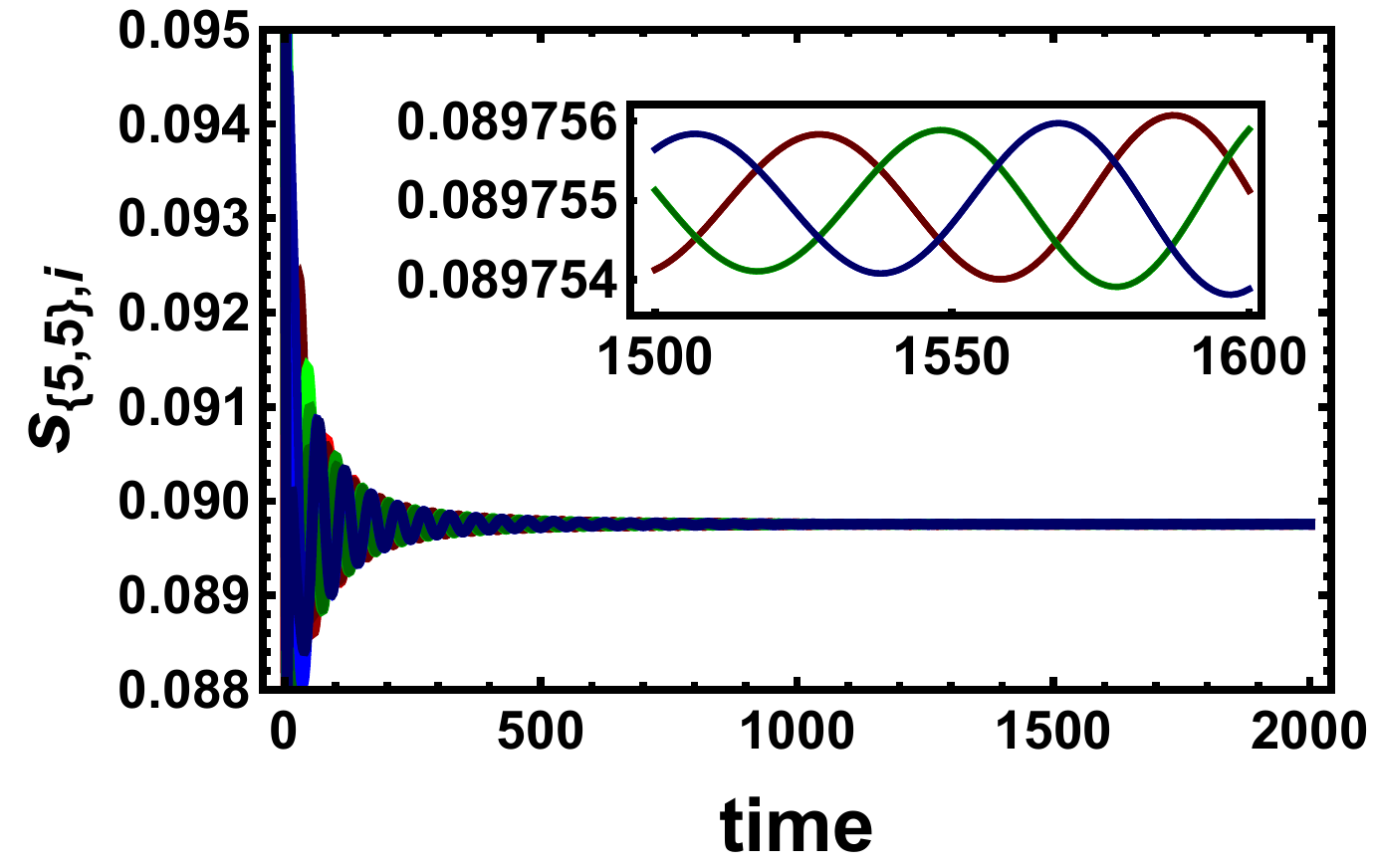}}
     \caption{Target waves emitted from a small set of pacemakers at the center of the grid. (a) Pacemakers chosen as 2L-units with $\gamma_x=0.8$ confined to a disc (white) of radius $R=8$ with color coded decay rates of CE-units with $\gamma_x \in [1.5, 1.6]$  outside the disc. (b) Snapshot  of target waves  taken at time instant $T=1565$. The different shades of the green, blue and red color indicate the maintained hierarchy of two time scales over a whole range outside the disc. (c): density plot of the concentration of the dominant item, plotted at each site at $T=1565$. Examples of time series are shown in (d) and (e) for the units at  location $\{60,60\}$ (one of the pacemakers), and $\{55,55\}$, close to the pacemakers, respectively. Units farther away (large red ring in (b)) apparently have the dynamics of a single hierarchy level ((f)), since the amplitudes for different shades of red are almost the same, and of a coexistence regime ((g)), although zooms in (f) and (g) display still remnants of the 2-level oscillations at the center. Other parameters are $r=1.25$, $e=0.2$, $f=0.4$, $c=d=2$, $\rho=1$, and $L=128$.}
    \label{fig:6}
    \end{centering}
\end{figure}

If the 2L-units as pacemakers  are confined to a circular region of radius $R$ on the grid, we observe target waves. The results are displayed  in Fig.~\ref{fig:6}. The size of the grid is again $128\times128$,
2L-units  have $\gamma_x=0.8$ inside a circular region of radius $R=8$ at the center. All units outside this region are ``defective'' (in a resting state) with  $\gamma_x$ chosen randomly from a uniform distribution in the interval $[1.5, 1.6]$. This distribution of the decay rates and therefore of the defects is color coded in Fig.~\ref{fig:6}(a).  It leads to the target wave pattern  of Fig.~\ref{fig:6}(b).
The wave fronts on the coarse scale of LHCs propagate inwards, whereas the fronts corresponding to the chasing inside an SHC propagate outwards, see the movie in \cite{supp}. In Fig.~\ref{fig:6}(c), we plot the concentration of the dominant-item (temporary winner) at each site.  Therefore, even though  the items at sites in the purple region have apparently the same concentration close to 0.09 (that is, close to their global coexistence equilibria), there are small oscillations as a remnant of the entrainment. For the units closer to the pacemaker disc, the two hierarchy levels are manifest in the switching between different shades of green, blue, and red in the inner rings. A zoom into the time evolution of the nine items at the inner site $(60,60)$ (Fig.~\ref{fig:6}(d)), or slightly more remote from the pacemakers at site $(55,55)$ (Fig.~\ref{fig:6}(e)), reveals the sequences of nine different temporary winners, partially differing in their amplitudes. Outside the red rings (outer red ring in (b)) with a time evolution according to (f), the dynamics seems to have collapsed to a 1L-unit with a heteroclinic cycle between coexistence equilibria, but a closer look into the motion shows still nine oscillating items with very small amplitudes and even smaller differences between their amplitudes. At the outmost site at (5,5) in the blue area of (b), the units seem to have approached a coexistence equilibrium (g), but also here remnants of a hierarchical heteroclinic motion are visible in a zoom. However, for the cases (f) and (g) the strong suppression of amplitudes and amplitude differences mean nevertheless a loss of information that was encoded in the winnerless competition sequence of information items closer to the disc of pacemakers.

\subsection{A quench in the bifurcation parameter}
Finally we would like to briefly mention another observation. So far it is preliminary and only based on numerical simulations. When we first let the target waves fully evolve and then quench the parameters of the pacemakers (2L-units) to parameters of CE-units, the whole grid of CE-units keeps on oscillating with hierarchical heteroclinic cycles for a rather long time. For example, if the quench is applied at time $t=1000$ time steps as in the movie of \cite{supp}, the waves, manifest in heteroclinic switching of items (colors), go on until at least $t=2500$. This means that the system effectively has a rather high ``inertia''. For future work it seems interesting  to pursue the very path in phase space that the system takes to the coexistence equilibrium after a quench in the bifurcation parameter. Moreover, our system should be complemented to account  for the energetic cost of maintaining synchronized heteroclinic motion in comparison to other modes of synchronization, and to consider the cost as a function of the system parameters and the type and  precision of the generated patterns.

\section{Summary and conclusions}\label{sec6}
We have analyzed heteroclinic units as pacemakers for entrainment to heteroclinic motion such as heteroclinic cycles or heteroclinic sequences (again cycles) with time scales of slow and fast oscillations. We studied such sets on chains, rings and two-dimensional grids. Entrainment to hierarchical heteroclinic motion is supported  by a location of the pacemaker at the edge of a chain for unidirectional coupling. A small back-coupling to the pacemaker  acts similarly to noise and reduces the slowing down of the heteroclinic motion, such that more units (themselves from the stable coexistence-regime) can be entrained. Moreover, the entrainment becomes easier when both the pacemaker and the driven units are closer to the bifurcation point, and the back-coupling  is not too strong; in the first case the vicinity to the bifurcation point reduces the slowing down of the pacemaker, in the second case the driven unit is closer to the regime  at which the coexistence equilibrium gets unstable.

The candidates for pacemakers were mainly distinguished  by their own parameters from the (hierarchical) heteroclinic dynamics as well as their location at the edge of a chain with strong forward and weak back-coupling to the pacemaker. On the two-dimensional grid, the coupling was bidirectional due to diffusion: there it was a small disc of pacemakers that generated an asymmetry in its neighborhood and led to target waves. Vice versa, when a unit in the resting state is at the edge of a chain with unidirectional coupling and all other units are in a mode of heteroclinic motion, these oscillations get stopped. This indicates the importance of the location of the pacemaker for directed coupling.

A small amount of noise facilitates the synchronization with the pacemaker as it prevents the driven units from exploring the substructure of their phase space. This may explain why the precision of stochastic patterns (in a brusselator model) is maximized for an intermediate thermodynamic cost \cite{barato} to suppress fluctuations. The substructure consists of  saddle equilibria which are generated under constant forcing and are seen as long as the pacemaker dwells on the respective saddles. The number of saddles increases with an increase in the distance of the driven unit to the pacemaker.  What the driven units experience is then a stepwise constant forcing from different dominant items of the respective saddles, interrupted by pulses during transits of the pacemaker between its saddle equilibria.

Moreover, as we have seen, initial conditions (in case of multistability) together with designed competition rates and a small amount of back-coupling allow the selection of different paths in the attractor space of the heteroclinic network. In view of brain dynamics this indicates a possibly rich repertoire of coding information via reproducible temporary sequences of entrained heteroclinic motion. The versatility would be generated  by taking different paths in the heteroclinic network. When assigned to a spatial grid, this setting  realizes hierarchical heteroclinic sequences, going on in parallel on different sites and getting synchronized due to the action of a pacemaker, without fine-tuning of the parameters of the driven units.

Furthermore, on a two-dimensional grid with diffusive coupling of pacemaker(s) and driven units, we have seen that a small disc of pacemakers in the center of the grid can entrain a whole grid to hierarchical heteroclinic motion. Loosely speaking, the target waves process the information that is encoded in the temporal sequence of temporarily dominant neural sub-populations (that is, the visited saddles along the path in the pacemaker's heteroclinic network). The entrainment works even if the main part outside the small disc is in a resting state  with a distribution of individual parameters from the coexistence regime. Also here, no fine-tuning of the parameters is needed for synchronization.
When the resting units are randomly assigned to the grid, for our concrete choice of parameters it is a fraction of around $60\%$ of working pacemakers that is sufficient to maintain the entrainment of the whole grid as a collective effect.

For future work it seems worthwhile to further explore how the selection of heteroclinic sequences in a heteroclinic network can be controlled as a function of the external input which we have not considered so far. A further extension will be an inclusion of the thermodynamic cost for this kind of intricate spatiotemporal pattern generation as well as a zoom into the transient dynamics after a quench in the pacemaker's parameters. When actually analyzing such patterns in brain activity, one may search for neuronal candidates that might play the role of  pacemakers.

\section*{Acknowledgment}
We thank the German Research Foundation (DFG) (grant number ME-1332/28-2) for financial support. One of us (B.T.) would like to thank Maximilian Voit (formerly at Jacobs University) for his help with the XMDS2 software.

\appendix
\section{}
\subsection{Effective model of a driven unit}
We consider the pacemaker to be very close to its saddle $(\rho/\gamma_P,0,0)$. Then we  substitute $s_{1,1}$ with $\rho/\gamma_P$ and $s_{1,2}=s_{1,3} = 0 $ in Eq.~\eqref{eq:4} and get the reduced system of equations
\begin{eqnarray}
d_ts_{2,1} &=& \rho s_{2,1} - \gamma_D s_{2,1}^2 - s_{2,1} (c s_{2,2} + e s_{2,3})+\delta  (\rho/\gamma_P - s_{2,1}),\nonumber\\
d_ts_{2,2} &=& \rho s_{2,2} - \gamma_D s_{2,2}^2 - s_{2,2} (c s_{2,3} + e s_{2,1})-\delta  s_{2,2},\nonumber\\
d_ts_{2,3} &=& \rho s_{2,3} - \gamma_D s_{2,3}^2 - s_{2,3} (c s_{2,1} + e s_{2,2})-\delta   s_{2,3}.\end{eqnarray}
This driven system has three equilibria of interest:
\begin{enumerate}
 \item A 1-item equilibrium: $(S_a,0,0)$ which is a stable node for forcing $\delta>\delta_{c_1}$.
In the original system, where the pacemaker is a heteroclinic cycle, the driven system has two additional equilibria, $(0,S_a,0), \text{and } (0,0,S_a)$, corresponding to the two other saddles of the pacemaker at $(0,\rho/\gamma_P,0), \text{and } (0,0,\rho/\gamma_P)$, respectively.
\item A 2-items equilibrium: $(S_b,S_c,0)$.
This equilibrium is a saddle equilibrium for $\delta > \delta_{c_1}$ and becomes a stable node for $\delta_{c_2}<\delta<\delta_{c_1}$.
Again, when the pacemaker performs a heteroclinic cycle, there are two additional two-items equilibria: $(0,S_b,S_c), \text{and } (S_c,0,S_b)$.
\item A 3-items equilibrium: $(S_d,S_e,S_f)$.
For $\delta>\delta_{c_2}$, this equilibrium is a saddle equilibrium,  for $\delta_{c_3}<\delta<\delta_{c_2}$, the equilibrium becomes a stable node. For $\delta_{c_4}<\delta<\delta_{c_3}$, the equilibrium is a stable focus-node and for $\delta < \delta_{c_4}$, it is an unstable focus node.  The unit driven by a pacemaker has two additional equilibria: $(S_f,S_d,S_e), \text{and }(S_e,S_f,S_d)$.
\end{enumerate}
\begin{figure}[ht!]
     \centering
 \includegraphics[width=0.75\textwidth]{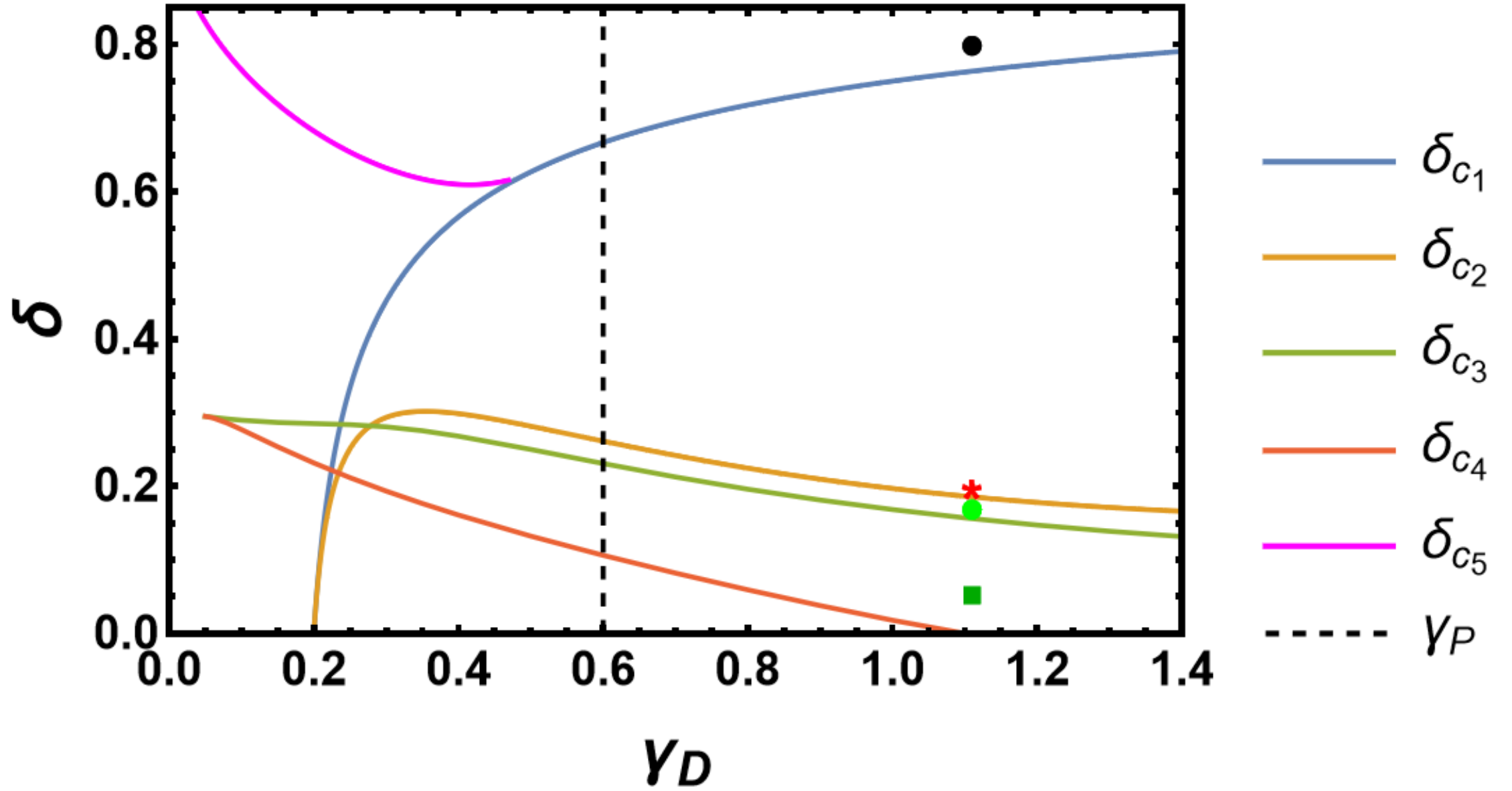}
\caption{Critical coupling strengths $\delta_c$ for the stability of different equilibria of a single  unit, driven under constant forcing $\delta$ to its first item, as a function of $\gamma_D$, the bifurcation parameter of the driven unit.   The other parameters are $\gamma_P = 0.6$, $\rho=1$, $c=2.0$, and $e=0.2$. In this work, $\gamma_D$ is chosen to the right of the dashed vertical line. The isolated symbols (green square, green dot, red star, black dot) correspond to the equilibria in Fig.~\ref{fig:A_2}. For further explanations see the text.}
 \label{fig:A1}
  \end{figure}
The expressions of $S_a,...,S_f$, $\delta_{c_1}$, and $\delta_{c_2}$ can be obtained analytically and are presented in \cite{supp}.
The expressions for $\delta_{c_3}$ and $\delta_{c_4}$ are only numerically accessible, thus we have plotted their numerical values along with $\delta_{c_1}$ and $\delta_{c_2}$ in Fig.~\ref{fig:A1} for a set of parameters.
In the following we restrict our discussion to the right of the vertical line in Fig.~\ref{fig:A1} with a unique correspondence of the equilibria between the critical lines, while multistability is observed to the left of this line that we do not further discuss. The vertical line corresponds to the choice of $\gamma_D=\gamma_P$.

\begin{figure}[ht!]
     \centering
 \includegraphics[width=\textwidth]{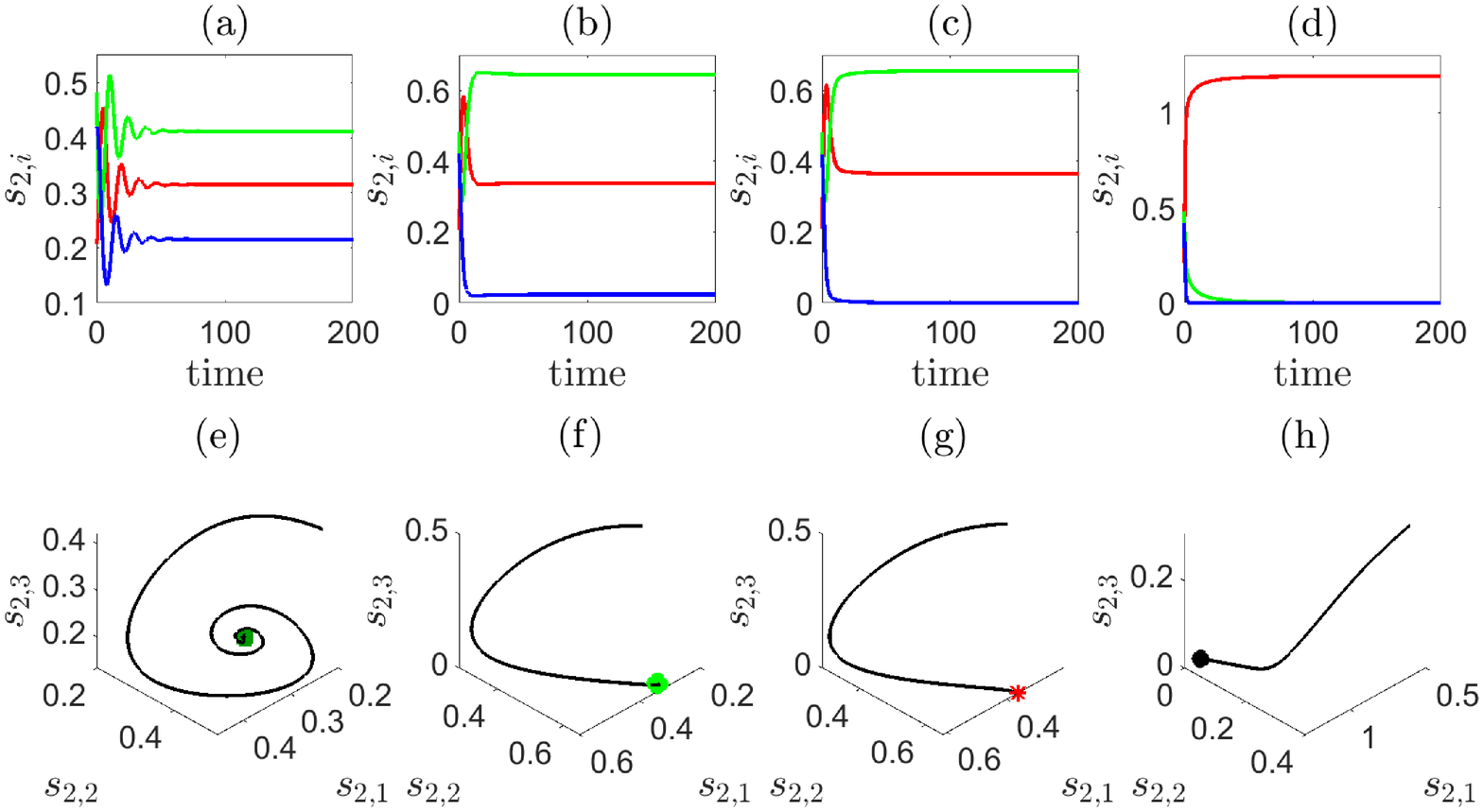}
\caption{Time series (upper panel)  and corresponding phase-space trajectories (lower panel) of the unit being driven by a constant forcing (a) $\delta=0.05$, (b) $\delta = 0.17$, (c) $\delta = 0.2$, and (d) $\delta = 0.8$. $\gamma_D=1.11$ chosen from the CE-regime. The 3S is shown by a green square (focus-node) or a green dot (stable node), the 2S by a red star, the 1S by a black dot. Other parameters are $\gamma_P = 0.6$, $\rho=1$, $c=2.0$, and $e=0.2$. The initial conditions are chosen randomly and uniformly from the interval [0,1].}
 \label{fig:A_2}
  \end{figure}

If we increase $\delta$ along a vertical line to the right of the dashed line in Fig.~\ref{fig:A1}, we observe the following features.
Below $\delta_{c_4}$, the driven unit has a stable limit cycle only if it is also in the regime of heteroclinic dynamics not considered further. If the driven unit is in the CE-regime, $\delta_{c_4}=0$. For $\delta_{c_4}<\delta<\delta_{c_3}$, the 3-items coexistence equilibrium (3S) is stable, it is a focus-node and therefore the trajectory spirals towards the equilibria (as shown in Fig.~\ref{fig:A_2}(a) for $\delta=0.05$ and $\gamma_D=1.11$). For $\delta_{c_3}<\delta<\delta_{c_2}$, the three-items equilibrium is still stable, but now it is a stable node (as shown in Fig.~\ref{fig:A_2}(b) for $\delta=0.17$). For $\delta_{c_2}<\delta<\delta_{c_1}$, the 2-items equilibrium (2S) becomes a stable node (as shown in Fig.~\ref{fig:A_2}(c) for $\delta=0.2$).
Finally, above $\delta_{c_1}$, the 1-item equilibrium (1S) becomes a stable node and the 2-items equilibrium loses stability (as shown in Fig.~\ref{fig:A_2}(d) for $\delta=0.8$). To the left of the vertical line, in the region bounded by $\delta_{c_1}$ from below and $\delta_{c_5}$ from above, 2-items equilibria coexist with 1-items equilibria. \\

\subsection{Back to the original system}
Now let us come back to our original system for which the pacemaker performs a heteroclinic cycle placed at the spatial index $k=1$. As the heteroclinic trajectory goes from saddle to saddle, the driven unit (at location $k=2$) goes from one of the above equilibria to the next equilibrium under the influence of the pacemaker.
Let us slowly increase the strength of the heteroclinic forcing, $\delta$, from zero on, and study its effect on the dynamics of the driven unit if $\gamma_D$ is chosen to the right of the vertical line in Fig.~\ref{fig:A1}.

\begin{figure}[ht!]
     \centering
 \includegraphics[width=\textwidth]{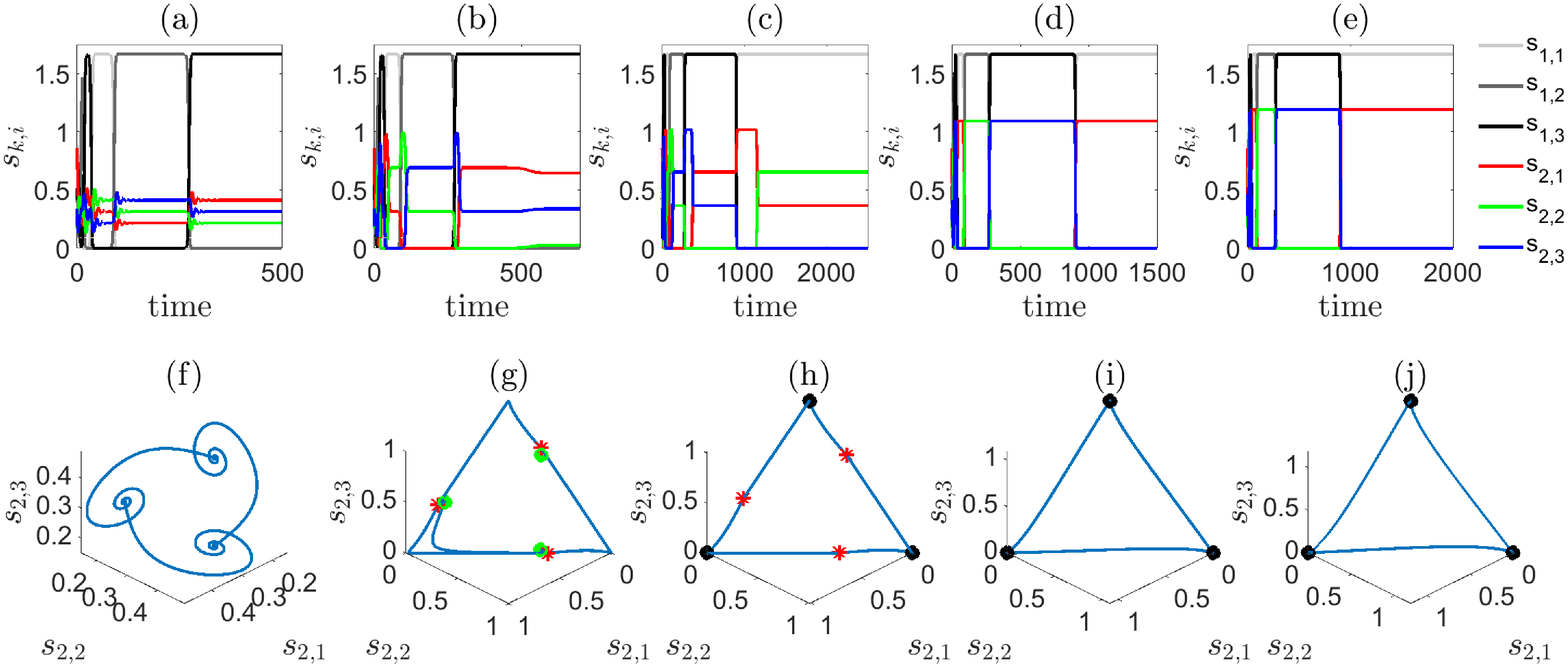}
\caption{Time series (upper panel (a)-(e))   and corresponding phase-space trajectories ((f)-(j))  of a unit at $k=2$ being driven by a heteroclinic cycle at $k=1$ for forcing (a) $\delta=0.05$, (b) $\delta = 0.17$, (c) $\delta = 0.2$, (d) $\delta = 0.4$, and (e) $\delta = 0.8$. $\gamma_P = 0.6$ from the HC-regime, $\gamma_D = 1.11$ from the CE-regime. In subfigures (a)-(e), the three items of the pacemaker are plotted in shades of gray, the items of the driven unit in red, green, blue.
In subfigures (f)-(j), the 1S are shown as black dots, the 2S
 as red stars, 3S  as green dots. The phase space trajectory of the driven unit is depicted  in blue. The 1S that the driven unit visits for $\delta = 0.4$ in (d) are not stable when adding a small amount of noise only to the driven unit. Other parameters are $\rho=1$, $c=2.0$, and $e=0.2$.}
 \label{fig:A3}
  \end{figure}

{\bf $\delta_{c_4}<\delta<\delta_{c_3}$ and $\gamma_D=1.11$.} For Fig.~\ref{fig:A3}(a)-(j), we consider the driven unit to be in the CE-regime  by taking $\gamma_D=1.11$.
When the forcing $\delta$ is between $\delta_{c_3}$ and $\delta_{c_4}$, the three 3S-equilibria are focus nodes. When the trajectory of the pacemaker goes to the vicinity of one of its saddles, the trajectory of the driven unit spirals to the corresponding focus-node. As the trajectory of the pacemaker goes to another saddle, the trajectory of the driven unit also moves to the corresponding focus node and stays there until the pacemaker moves to its next saddle as shown in Fig.~\ref{fig:A3}(a).

{\bf $\delta_{c_3}<\delta<\delta_{c_2}$.} For this range, the system with constant forcing goes to the 3S-equilibrium (Fig.~\ref{fig:A_2}(b)), but the trajectory of the unit driven by a heteroclinic cycle tries to first visit the 1S- and 2S-equilibria, pushed by the pacemaker. However, since the 1S and 2S-equilibria are unstable (under constant forcing), if the driven unit has enough time while the pacemaker dwells on the corresponding saddle, it relaxes finally  to one of the 3S-equilibria (which are stable as long as constant forcing is applied) (see Fig.~\ref{fig:A3}(b)).

{\bf $\delta_{c_2}<\delta<\delta_{c_1}$.} Similarly, the trajectory of the driven unit visits the 1S-equilibria in this range before relaxing to the stable  2S-state in Fig.~\ref{fig:A3}(c) (`stable'  again referring to constant forcing), or it may also keep on residing in the unstable  1S-equilibrium as shown in Fig.~\ref{fig:A3}(d) for $\delta=0.4$. However, if we add even a very small amount of noise to the dynamics of the driven unit, the trajectory relaxes to the 2S-equilibrium (not shown here).

{\bf $\delta>\delta_{c_1}$.} Finally, for $\delta>\delta_{c_1}$, the trajectory of the driven unit goes to one of its 1S-stable nodes when the pacemaker is in the vicinity of the corresponding saddle (see Fig.~\ref{fig:A3}(e)). So for large values of $\delta$ (above $\delta_{c_1}$), the pacemaker is able to keep the driven unit entrained to its dynamics.

In summary, based on the effective modeling of the driven unit under stepwise constant forcing, the driven unit directly coupled  to the pacemaker approaches a 1S, 2S, or 3S-saddle equilibrium.

If we introduce a back-coupling $\delta_b$ to the pacemaker,
a first effect  is  a slowdown of  the slowing down of the heteroclinic oscillations. Thus, increasing $\delta_b$ is similar to increasing $\gamma_P$, thus effectively shifting the pacemaker  closer to the coexistence regime with less pronounced slowing down. Beyond a certain strength of $\delta_b$, the heteroclinic dynamics loses stability to limit cycle oscillations, and after a further increase, the limit cycles lose stability in a Hopf bifurcation, and the dynamics of both units stabilize to their respective coexistence equilibria. Further details are presented in \cite{supp}.

\section{}
\subsection{Zooming into the proliferation of equilibria}
Consider a driven unit. Depending on whether the trajectory of its respective nearest-neighbor pacemaker (its driving unit) is in the vicinity of a 1S, 2S, or 3S-equilibrium, the driven unit experiences forcing to either one, two, or all three of its items. As discussed before for a single driven unit, if the driven unit experiences a constant force at one of its items it has three types of equilibria: a 1S, 2S, or 3S-equilibrium. Here, in addition, when the driven unit at $k\geq 3$ is forced at two  items at the same time, it has two equilibria: a 2S and a 3S. If it is forced at all three of its items, it has only one equilibrium of interest, a 3S.

\begin{figure}[ht!]
     \centering
 \includegraphics[width=\textwidth]{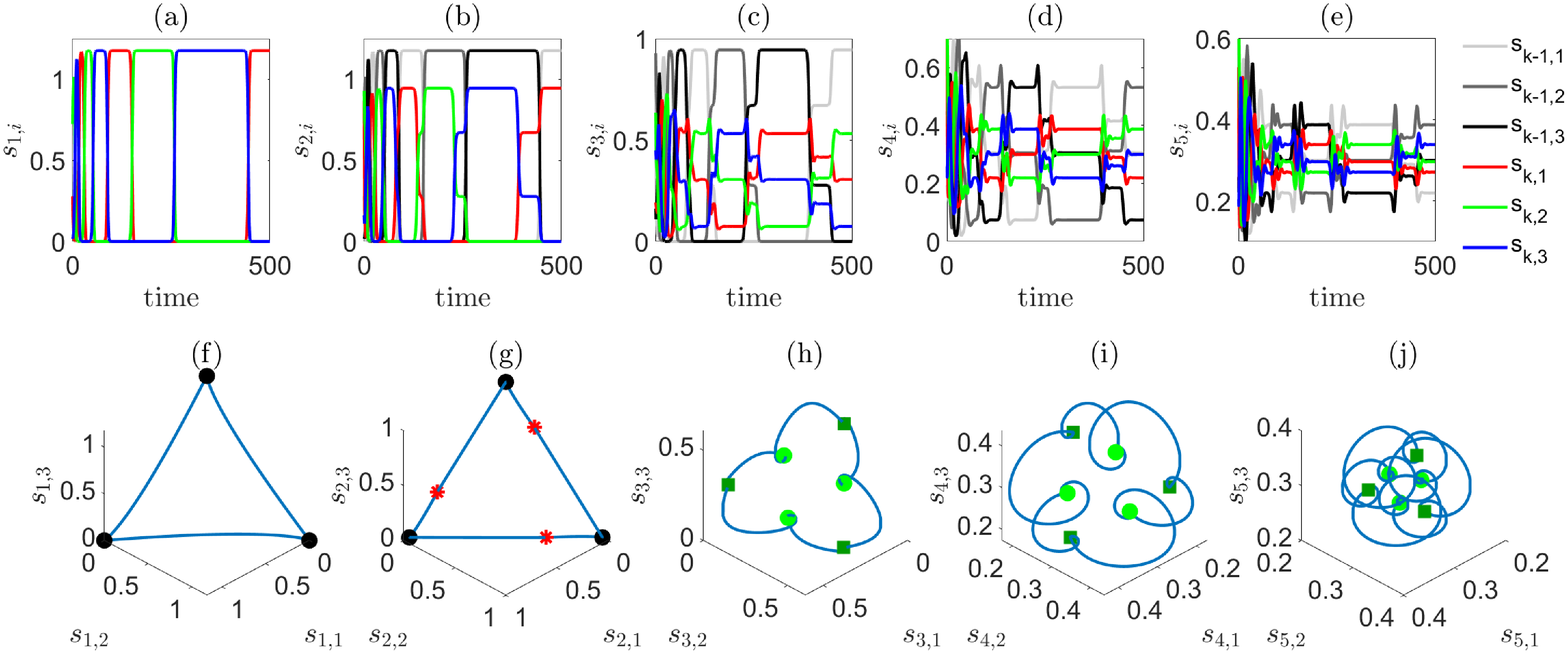}
\caption{Time series and phase-space trajectories of a system of five units unidirectionally coupled along a chain, from left to right at positions $k=1$ to $k=5$. The parameters are $\gamma_P=0.85$, $\gamma_D=1.11$. In (a)-(e), the time series of the three items of the $k^{th}$ unit are  color-coded in red, green, blue. The items of the $k-1^{th}$ unit, acting as a pacemaker to the $k^{th}$ unit, are shown in shades of gray. In the lower panel 1S in phase space are shown by black circles, 2S by red stars, the different 3S by dark green squares and light green circles. Other parameters are $\delta = 0.2$, $\rho=1$, $c=2.0$, and $e=0.2$.}
 \label{fig:B1}
  \end{figure}

As long as there is a pronounced slowing down of the heteroclinic motion, a driven unit is exposed for a long time interval to a constant forcing corresponding to the trajectory of its pacemaker approaching one of its equilibria and dwelling there, while for a comparatively very small time interval it is exposed to time-varying forcing when the pacemaker quickly transits from one of its equilibria to the next. Therefore the forcing experienced by the driven unit resembles a sequence of step functions with different constant forcing while different equilibria are visited by the driving unit,  interrupted by jumps, when the driving unit transits between different equilibria. The transitions are visible in spikes of the trajectories, the more pronounced, the farther the distance from the original pacemaker at $k=1$, as seen in Fig.~\ref{fig:B1}. 

This figure shows the time series and phase-space trajectories of a system of five units unidirectionally coupled along a chain at sites $k=1$ to $k=5$. The time series of the three items of the $k^{th}$ unit are  color-coded in red, green, blue. The items of the $k-1^{th}$ unit, acting as a pacemaker to the $k^{th}$ unit, are shown in shades of gray. The amplitudes decrease with the distance to the pacemaker, as well as the differences in the dominance of  different items. In the lower panel the phase space trajectories change from clear heteroclinic motion at $k=1$ and $k=2$ to increasingly more involved trajectories, exposed to more saddles and concentrated on a smaller part of phase space.

The possible stable or unstable equilibria are listed in Table.~\ref{table:2} of the main text.
The approaches to these equilibria is determined by the following rules:
\begin{subequations}\label{eqrules}
\begin{align}
1S\:\rightarrow\:1S \:\text{or}\: 2S \:\text{or}\: 3S \:\text{or}\: (1S\rightarrow 2S)\:\text{or}\:(1S\rightarrow 2S\rightarrow 3S)\: \text{or}\:(2S\rightarrow 3S),\\
2S\:\rightarrow\: 2S\:\text{or}\: 3S\:\text{or}\:(2S\rightarrow 3S),\\
3S\rightarrow \:3S.
\end{align}
\end{subequations}
For example, for unit-2  the rules should be read as follows: While unit-1 is at one of its three 1S-equilibria, unit-2 either dwells directly at one of its 1S, 2S, or 3S equilibria, or it moves to a more stable equilibrium according to the rules in the brackets of Eq.~\ref{eqrules}(a).
The options for unit-3  can be read off from Eq.~\ref{eqrules} in the same way, if unit-2 is in the vicinity of 1S, or 2S, or 3S.
Similarly for more remote distances from the pacemaker. Thus, the number of equilibria that are available to a unit to visit in its phase space increases, as the distance  from the  pacemaker at $k=1$ increases.

\subsection{Impact of  $\delta$ on the entrainment length}
\begin{figure}[ht!]
     \centering
 \includegraphics[width=\textwidth]{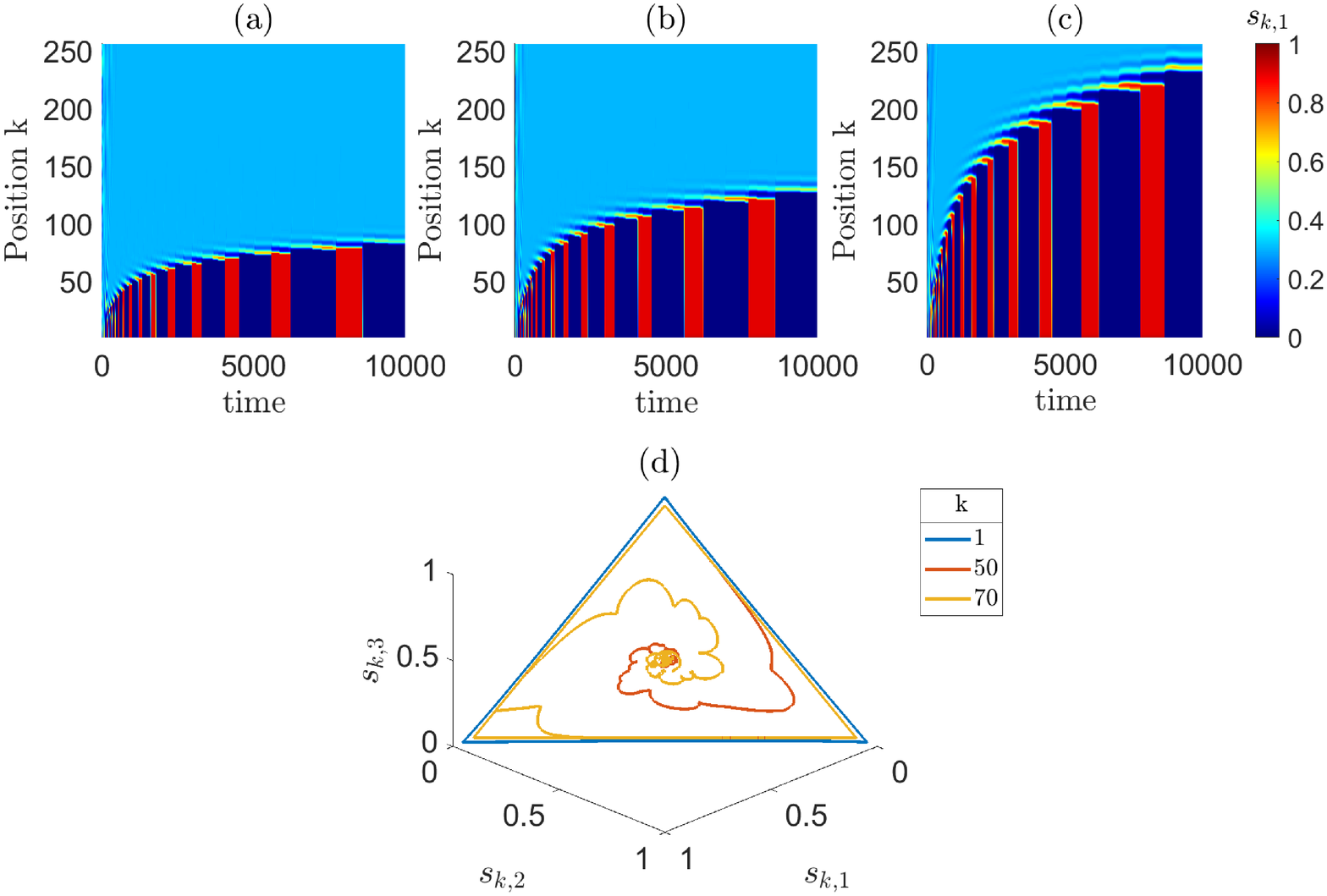}
\caption{Spatiotemporal plots of the concentration of the first item $s_{k,1}$ of each unit along the unidirectional chain of length $N=256$ for three values of the coupling strength (a) $\delta=0.5$, (b) $\delta=0.6$, and (c) $\delta=0.75$.  The phase-space trajectories of a few selected units are plotted in (d) for $\delta=0.5$, showing shaky trajectories due to the impact of many unstable equilibria. Other parameters are $\gamma_P=1.05$, $\gamma_D=1.11$, $c=2$, $e=0.2$, and $\rho = 1.0$.}
 \label{fig:B2}
  \end{figure}

The entrainment length corresponds to the number of units that get entrained to heteroclinic dynamics in a certain time window. For our numerical results, we consider the length of the chain to be $N=256$ units. The results of varying $\delta$  are shown in Fig.~\ref{fig:B2}. In panels Fig.~\ref{fig:B2}(a)-(c), we plot the spatiotemporal plots of the first items $s_{k,1}$ of each unit for three values of the coupling strength $\delta=0.5$, $\delta=0.6$, and $\delta=0.75$, respectively. In the red region of the plots, the concentration of items $s_{k,1}$ is dominant at the location $k$, in the dark blue region it approaches zero while either $s_{k,2}$ or $s_{k,3}$ are dominant (winnerless competition). The entrainment process takes some transient time that is not negligible.  The stronger the forcing $\delta$, the more units synchronize to the heteroclinic motion.
In Fig.~\ref{fig:B2}(d), we plot the  phase-space trajectories of a few units along the chain for $\delta=0.5$. While the pacemaker follows a heteroclinic cycle (blue line in (d)), a driven unit at position 50 wiggles on its approach of the heteroclinic cycle (red line) and even more so a unit at position 70 (yellow line). These features can be traced back to the proliferated saddle equilibria, felt by the respective trajectories. Here we also see that the effective modeling in terms of constant forcing to one, two, or three items loses its applicability, the more, the less pronounced the time scale separation between the dwell times in the vicinity of the saddles  and the transits to the subsequent saddles is.

\subsection{Impact of $\gamma_P$ and $\gamma_D$ on the entrainment length}
Another parameter that affects the entrainment length is the bifurcation parameter $\gamma_P$ of the pacemaker.
Here, the entrainment length increases with increasing $\gamma_P$, this means, the closer the pacemaker is to the bifurcation towards the CE-regime. At a first view this may seem counterintuitive, as one may expect that the deeper the pacemaker is in the regime of heteroclinic dynamics, the stronger is its impact on the driven units. However, what prevents an easy entrainment with a pacemaker far off the bifurcation point is the fact that the heteroclinic motion slows down much faster for parameters deeply in the heteroclinic region.

The bifurcation parameter of the driven units $\gamma_D$ also affects the entrainment length. Here, as expected, the entrainment length decreases with an increase in $\gamma_D$ for which the driven units move deeper into the CE-regime and are harder
for the pacemaker to entrain.

Moreover it should be mentioned that the entrainment length also depends on the initial conditions. The reason is that depending on the initial conditions, the transient time it takes the driving pacemaker to approach the vicinity of the heteroclinic trajectory differs and accordingly the time which is spent in the vicinity of the saddles. These differences propagate along the chain and get amplified, so that in a given time window the numbers of entrained units may slightly differ, even more so for a 2L-unit as a pacemaker, as there the initial conditions determine the number of revolutions in the SHCs.

\subsection{Ring configurations} If we close the unidirectionally coupled chain to a ring (Fig.~\ref{fig:1}(b)) such that only the $N^{th}$ unit applies a weak forcing $\delta_P$ to the pacemaker (weaker in comparison to the unidirectional forcing $\delta$ between all other units), the entrainment depends on the size of the ring. Larger rings may not be synchronized at all to heteroclinic dynamics. As $\delta_P$ is increased and other parameter are kept fixed, the system shows quasi-periodic dynamics coexisting with heteroclinic dynamics, for further increased $\delta_P$, limit cycles and  finally a stable coexistence equilibrium, details are described in \cite{supp}.

\subsection{Bidirectional coupling}
In \cite{supp} we discuss the role of bidirectional, but strongly asymmetric nearest neighbor coupling between a pacemaker and $N-1$ driven units. From a physical point of view some back-coupling seems more realistic than unidirectional coupling. The forward coupling $K_{k,l}=\delta$ if $k\in\{2,...,N\}$ and $l=k-1$, the back-coupling $K_{k,l}=\delta_b$ if $k\in\{1,...,N-1\}$ and $l=k+1$, $K_{k,l}=0$ otherwise. The main difference is that each unit is now sensitive to the dynamics on both sides along the chain, the entrainment length becomes sensitive to the size of the chain. Initially more and more units get entrained, but after reaching a peak value over time in the number of entrained units, some recently entrained ones start de-entraining when they receive the feedback from the back-coupling, and then the entrainment length settles to a value less than the peak. The stronger the back-coupling, the smaller the entrainment length.

\section*{References}
\providecommand{\newblock}{}


\begin{thebibliography}{10}
\expandafter\ifx\csname url\endcsname\relax
  \def\url#1{{\tt #1}}\fi
\expandafter\ifx\csname urlprefix\endcsname\relax\def\urlprefix{URL }\fi
\providecommand{\eprint}[2][]{\url{#2}}

\bibitem{EEG}
Michel C~M and Koenig T 2018 {\em Neuroimage\/} {\bf 180} 577--593

\bibitem{rabi3}
Rabinovich M, Volkovskii A, Lecanda P, Huerta R, Abarbanel H and Laurent G 2001
  {\em Physical review letters\/} {\bf 87} 068102

\bibitem{laurent}
Laurent G, Stopfer M, Friedrich R~W, Rabinovich M~I, Volkovskii A and Abarbanel
  H~D 2001 {\em Annual review of neuroscience\/} {\bf 24} 263--297

\bibitem{binding}
Afraimovich V, Gong X and Rabinovich M 2015 {\em Chaos: An Interdisciplinary
  Journal of Nonlinear Science\/} {\bf 25} 103118 (\textit{Preprint}
  \eprint{https://doi.org/10.1063/1.4932563})
  \urlprefix\url{https://doi.org/10.1063/1.4932563}

\bibitem{rabi1}
Rabinovich M~I, Varona P, Tristan I and Afraimovich V~S 2014 {\em Frontiers in
  computational neuroscience\/} {\bf 8} 22

\bibitem{valentin1}
Afraimovich V, Tristan I, Huerta R and Rabinovich M~I 2008 {\em Chaos: An
  Interdisciplinary Journal of Nonlinear Science\/} {\bf 18} 043103
  (\textit{Preprint} \eprint{https://doi.org/10.1063/1.2991108})
  \urlprefix\url{https://doi.org/10.1063/1.2991108}

\bibitem{nowak}
Nowak M~A and Sigmund K 2002 {\em Nature\/} {\bf 418} 138--139

\bibitem{szolnoki}
Szolnoki A, Mobilia M, Jiang L~L, Szczesny B, Rucklidge A~M and Perc M 2014
  {\em Journal of the Royal Society Interface\/} {\bf 11} 20140735

\bibitem{friston}
Friston K~J 2018 {\em private communication\/}

\bibitem{max1}
Voit M and Meyer-Ortmanns H 2018 {\em The European Physical Journal Special
  Topics\/} {\bf 227} 1101--1115

\bibitem{max2}
Voit M and Meyer-Ortmanns H 2019 {\em Applied Mathematics and Nonlinear
  Sciences\/} {\bf 4} 279--288

\bibitem{max3}
Voit M and Meyer-Ortmanns H 2019 {\em Physical Review Research\/} {\bf 1}
  023008

\bibitem{daido1}
Daido H and Nakanishi K 2007 {\em Physical Review E\/} {\bf 75} 056206

\bibitem{daido2}
Daido H and Nakanishi K 2004 {\em Physical review letters\/} {\bf 93} 104101

\bibitem{radicchi1}
Radicchi F and Meyer-Ortmanns H 2006 {\em Physical Review E\/} {\bf 73} 036218

\bibitem{radicchi2}
Radicchi F and Meyer-Ortmanns H 2006 {\em Physical Review E\/} {\bf 74} 026203

\bibitem{rabi2}
Rabinovich M~I, Huerta R and Varona P 2006 {\em Physical review letters\/} {\bf
  96} 014101

\bibitem{dawes1}
Dawes J and Tsai T~L 2006 {\em Physical Review E\/} {\bf 74} 055201

\bibitem{supp}
Supplementary material

\bibitem{townsend1}
Townsend R~G, Solomon S~S, Chen S~C, Pietersen A~N, Martin P~R, Solomon S~G and
  Gong P 2015 {\em Journal of Neuroscience\/} {\bf 35} 4657--4662

\bibitem{townsend2}
Townsend R~G and Gong P 2018 {\em PLoS computational biology\/} {\bf 14}
  e1006643

\bibitem{barato}
Rana S and Barato A~C 2020 {\em Physical Review E\/} {\bf 102} 032135

\bibitem{kirk}
Kirk V and Silber M 1994 {\em Nonlinearity\/} {\bf 7} 1605--1621

\bibitem{ashwin1999}
Ashwin P and Field M 1999 {\em Archive for Rational Mechanics and Analysis\/}
  {\bf 148} 107--143

\bibitem{ashwin2016}
Ashwin P and Postlethwaite C 2016 {\em Journal of Nonlinear Science\/} {\bf 26}
  345--364

\bibitem{grossberg}
Grossberg S 2000 {\em Trends in cognitive sciences\/} {\bf 4} 233--246

\bibitem{max4}
Voit M and Meyer-Ortmanns H 2020 {\em Physical Review Research\/} {\bf 2}
  043097

\bibitem{moelle}
M{\"o}lle M and Born J 2011 {\em Progress in brain research\/} {\bf 193}
  93--110

\bibitem{saleh}
Saleh M, Reimer J, Penn R, Ojakangas C~L and Hatsopoulos N~G 2010 {\em
  Neuron\/} {\bf 65} 461--471

\bibitem{ashwin2013}
Ashwin P and Postlethwaite C 2013 {\em Physica D: Nonlinear Phenomena\/} {\bf
  265} 26--39

\bibitem{kurabook}
Kuramoto Y 1984 {\em Chemical oscillations, waves, and turbulence\/} (Springer)

\bibitem{xmds2}
Dennis G~R, Hope J~J and Johnsson M~T 2013 {\em Computer Physics
  Communications\/} {\bf 184} 201--208

\bibitem{armbruster}
Armbruster D, Stone E and Kirk V 2003 {\em Chaos: An Interdisciplinary Journal
  of Nonlinear Science\/} {\bf 13} 71--79

\bibitem{bakhtin}
Bakhtin Y 2010 {\em Dynamical Systems\/} {\bf 25} 413--431

\bibitem{stone}
Stone E and Armbruster D 1999 {\em Chaos: An Interdisciplinary Journal of
  Nonlinear Science\/} {\bf 9} 499--506

\bibitem{odor1}
Laurent G, Stopfer M, Friedrich R~W, Rabinovich M~I, Volkovskii A and Abarbanel
  H~D 2001 {\em Annual review of neuroscience\/} {\bf 24} 263--297

\bibitem{odor2}
Nowotny T, Huerta R, Abarbanel H~D and Rabinovich M~I 2005 {\em Biological
  cybernetics\/} {\bf 93} 436--446

\end{thebibliography}
\end{document}


\title{Supplementary Material:\\Heteroclinic units acting as pacemakers:
Entrained dynamics for cognitive processes}
 \author{Bhumika Thakur and Hildegard Meyer-Ortmanns}
 \address{School of Science, Jacobs University Bremen, Campus Ring 1, 28759 Bremen, Germany}
\vspace{10pt}
\begin{indented}
\item[] July 2022
\end{indented}

\renewcommand{\thepage}{S\arabic{page}}
\renewcommand{\thesection}{S\arabic{section}}
\renewcommand{\thetable}{S\arabic{table}}
\renewcommand{\thefigure}{S\arabic{figure}}

\section{Movies}
\subsection{Target wave on a 2D grid}
The video shows target waves on a two-dimensional grid of 2L-units when the pacemakers are confined to a disc at the center as discussed in section 5 of the main text.
The video runs from $t = 1450$ to $t = 1600$ time points. The wave fronts on the coarse scale of LHCs (shown by different colors: red, green, and blue) propagate inwards whereas the fronts corresponding to the chasing inside an SHC (different shades of a colors) propagate outwards.

\subsection{Target wave with quenching}
The video shows the simulation when we first let the target waves fully evolve till $t<1000$ time points and apply a quench at $t = 1000$. The video runs from  $t = 850$ to $t = 2500$ time points. The units keep oscillating even until the end, shown by the switching of colors. The two hierarchy levels are visible in different shades of the same color until $t=1500$ time points and the switching between different colors goes on till the end of the simulation.

\section{Equilibria of a 1L-unit driven by  constant forcing}\label{S1}
As discussed in the Appendix A of the main text, a 1L-unit driven by  constant forcing has three equilibria of interest:
\begin{enumerate}
 \item A 1-item equilibrium: $(S_a,0,0)$, where
 \begin{equation}
  S_a = \frac{-\delta+\rho+\sqrt{4\delta\rho\gamma_D/\gamma_P+(\delta-\rho)^2}}{2\gamma_D}.
 \end{equation}
The eigenvalues of this equilibrium are $(-2 S_a \gamma_D - \delta+\rho, -S_a c - \delta+\rho, -S_a e - \delta+\rho)$. The first eigenvalue is along the radial direction, the second  along the contracting and the last one  along the expanding direction.
For our choice of parameters, radial and contracting eigenvalues are always real and negative.  The expanding eigenvalue is also real and negative for $\delta>\delta_{c_1}$, where
\begin{equation}
\delta_{c_1} = \rho \left[ 1- \frac{e}{2 \gamma_P (e-\gamma_D)}\left(e-\sqrt{e^2-4\gamma_P(e-\gamma_D)}\right) \right].
\end{equation}
Below $\delta_{c_1}$, the expanding eigenvalue becomes positive and therefore, the equilibrium becomes a saddle-equilibrium.

\item A 2-items equilibrium: $(S_b,S_c,0)$, such that
\begin{equation}
 S_b =  \frac{(c - \gamma_D)  (\delta - \rho)}{2 (-c e + \gamma_D^2)}+ G_1,
\end{equation}
\begin{equation}
 S_c = \frac{(c e + (e - 2 \gamma_D) \gamma_D) (\delta - \rho)}{2 \gamma_D (-c e + \gamma_D^2) } -\frac{e}{\gamma_D} G_1,
\end{equation}
where
\begin{equation}
 G_1 = \frac{\sqrt{(c -\gamma_D)^2 (\delta - \rho)^2 +
 4 (\gamma_D/\gamma_P)  (-c e +\gamma_D^2)\delta \rho}}{2 (-c e +\gamma_D^2) }.
\end{equation}
This equilibrium is a saddle equilibrium for $\delta > \delta_{c_1}$ and becomes a stable node for $\delta_{c_2}<\delta<\delta_{c_1}$, where
\begin{eqnarray}
 \delta_{c_2} &=& \rho\left(1 - \frac{1 - \sqrt{1 - 2 G_2}}{G_2}\right),\nonumber\\
 G_2 &=& \frac{2 (e - \gamma_D) (c^2 + e^2 + \gamma_D^2 - e \gamma_D  -
   c e - c \gamma_D) \gamma_P}{(e^2 - c \gamma_D)^2}.
\end{eqnarray}
%
\item A 3-items equilibrium: $(S_d,S_e,S_f)$, such that
\begin{eqnarray}
S_d &=&  \frac{(ce-\gamma_D^2)(-m_1 + \sqrt{m_1^2-2\delta\rho m_2})}{m_2}\\
S_e &=& -\frac{m_1 ((c^2-e\gamma_D)-2(\gamma_D^2-ce)) +(c^2-e\gamma_D)\sqrt{m_1^2-2\delta\rho m_2}}{m_2}\\
S_f &=& -\frac{m_1 ((e^2-c\gamma_D)-2(\gamma_D^2-ce)) +(e^2-c\gamma_D)\sqrt{m_1^2-2\delta\rho m_2}}{m_2},
\end{eqnarray}
where $m_1 = (c^2 + e^2 + \gamma_D^2 - e \gamma_D -
   c e - c \gamma_D)(\delta-\rho) \gamma_P$, and $m_2 = 2\gamma_P (c^3+e^3+\gamma_D^3-3 c e\gamma_D) (ce-\gamma_D^2)$.
\end{enumerate}

\section{Bidirectional coupling of two 1L-units and non-vanishing $\delta_b$}
\begin{figure}[ht!]
     \centering
  \subfloat[]{\includegraphics[width = 0.19 \textwidth]{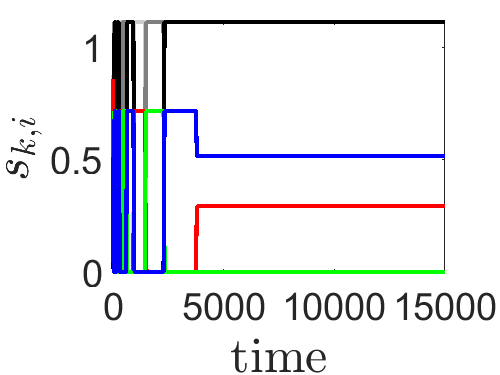}}
    \subfloat[]{\includegraphics[width = 0.19 \textwidth]{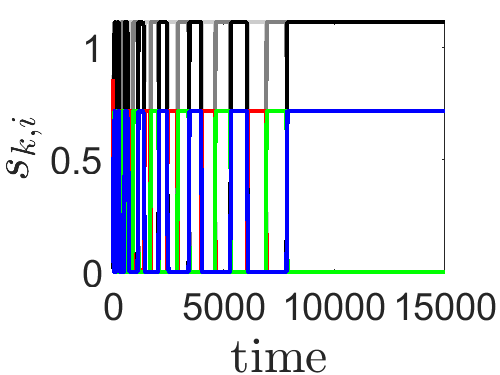}}
    \subfloat[]{\includegraphics[width = 0.19 \textwidth]{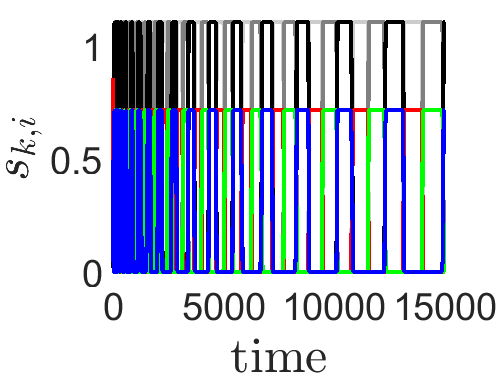}}
    \subfloat[]{\includegraphics[width = 0.19 \textwidth]{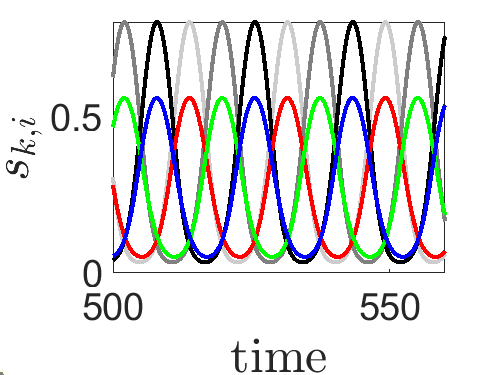}}
    \subfloat[]{\includegraphics[width = 0.19 \textwidth]{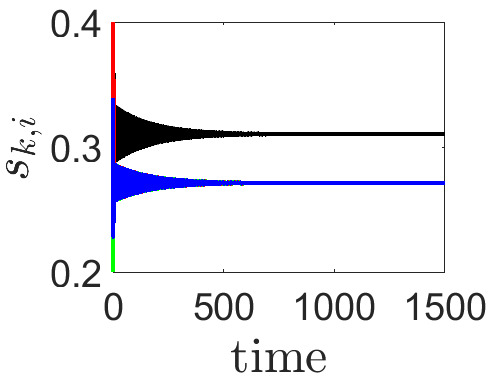}}
    {\includegraphics[width = 0.05 \textwidth]{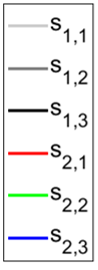}}\\
    \subfloat[]{\includegraphics[width = 0.25 \textwidth]{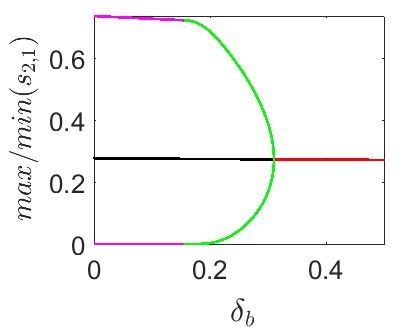}}
    \subfloat[]{\includegraphics[width = 0.25 \textwidth]{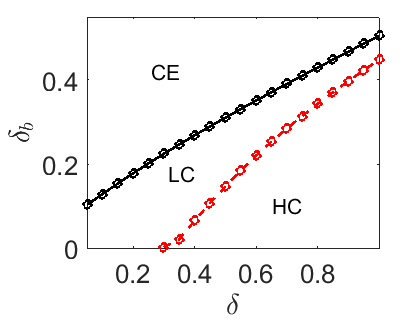}}
    \subfloat[]{\includegraphics[width = 0.25 \textwidth]{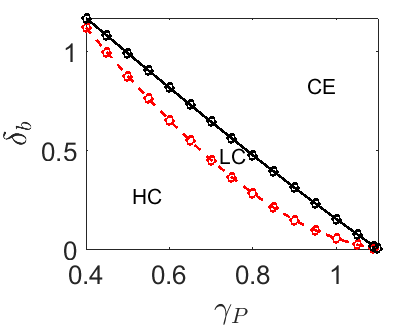}}
    \subfloat[]{\includegraphics[width = 0.25 \textwidth]{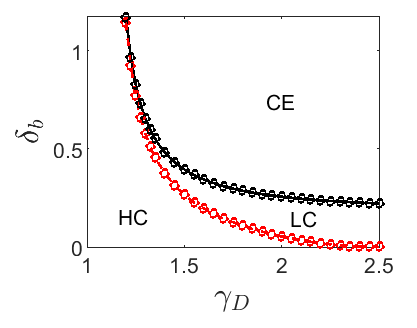}}
\caption{Effect of introducing a back-coupling ($\delta_b$) to the pacemaker on the dynamics of two 1L-units. For (a)-(e),  $\delta=0.4$. In (a) $\delta_b=0$,  (b) $\delta_b=10^{-8}$,  (c) $10^{-2}$ ,(d)  $\delta_b=0.2$, (e) $\delta_b=0.3$.  (f): bifurcation diagram where the maximum and minimum of the first item $s_{2,1}$ of the driven unit are plotted against $\delta_b$ for fixed value of $\delta=0.5$. For (a)-(g), $\gamma_P=0.9$,  $\gamma_D = 1.7$. The regime of stable heteroclinic dynamics (HC), limit cycles (LC), and coexistence equilibria (CE) are shown by magenta, green, and red lines, respectively. Before the Hopf bifurcation point, the coexistence equilibrium is unstable (black line). In (g)-(i), three dynamical regimes in the parameter space of $\delta-\delta_b$, $\gamma_P-\delta_b$, and $\gamma_D-\delta_b$, respectively. In (h), $\delta = 0.5$ and $\gamma_D=1.7$. In (i), $\delta = 0.5$ and $\gamma_P=0.9$.
Other parameters are $c=2.0$, $e=0.2$, and $\rho=1$.
}
 \label{fig:back_coup}
  \end{figure}

If the back-coupling ($\delta_b$) from the driven unit to the pacemaker is also present, the coexistence equilibrium state which is otherwise unstable for the pacemaker  becomes stable for both units in a certain parameter range.
Depending on the parameters ($\gamma_P$, $\gamma_D$, $\delta$, and $\delta_b$), the bi-directional coupling can either lead to heteroclinic dynamics or limit cycle behavior or stabilization of the coexistence equilibrium state.
In Fig.~\ref{fig:back_coup}, we show the effect of $\delta_b$ on the dynamics of the system.  In Fig.~\ref{fig:back_coup}(a), with $\delta_b=0$, we choose the other parameters such that the trajectory of the driven unit eventually settles to a 2S-equilibrium when the dynamics of the pacemaker slows down considerably. If we  now introduce a small amount of back-coupling (see Fig.~\ref{fig:back_coup}(b)), the first effect we observe is that the dynamics of the pacemaker slows down much more slowly than  for $\delta_b=0$. Also, it is able to keep the driven unit entrained even after its dynamics slows down considerably. As we keep increasing $\delta_b$, the dynamics of the pacemaker takes more and more time to slow down as shown in Fig.~\ref{fig:back_coup}(c). Thus, the effect is similar to increasing the value of $\gamma_P$. If we increase $\gamma_P$ in the HC-regime, the dynamics of the unit also takes longer and longer to slow down. Thus,  by increasing the back-coupling, the driven unit effectively pulls the pacemaker closer to the CE-regime. However, beyond a certain value of $\delta_b$, the heteroclinic dynamics loses stability to limit cycle oscillations as shown in Fig.~\ref{fig:back_coup}(d). On further increasing $\delta_b$ beyond another critical value corresponding to a Hopf bifurcation point, the limit cycle loses stability and the dynamics of both units stabilize to their respective coexistence equilibrium as shown in Fig.~\ref{fig:back_coup}(e).

In Fig.~\ref{fig:back_coup}(f), we plot the corresponding bifurcation diagram with $\delta_b$ as the bifurcation parameter for a fixed value of $\delta=0.5$. Until around $\delta_b<0.146$, both units are entrained to heteroclinic dynamics shown by a magenta colored line. Beyond $\delta_b>0.146$, the HC loses stability and gives way to limit cycle oscillations which remain stable till the system undergoes a Hopf bifurcation at $\delta_b=0.31$. The regime of limit cycle oscillations is shown in green. At the Hopf bifurcation, the limit cycle loses stability to coexistence equilibria.  The stable equilibria are shown in red. Before the Hopf bifurcation point, the coexistence equilibrium is unstable as shown by the black line.
The death rates $\gamma_P$ and $\gamma_D$ are chosen as $0.9$ and $1.7$, respectively.

If we just focus on what supports entrainment to heteroclinic motion, based on the phase diagrams in Figs.~\ref{fig:back_coup}(g)-(i), it is a stronger (forward) coupling from the pacemaker to the driven units than the back-coupling in the reverse direction, and the stronger the forward-coupling is, the larger the back-coupling may be. The closer the pacemaker is to the bifurcation point, the less back-coupling is compatible with such entrainment. The deeper the driven unit is in the CE-regime, the smaller the back-coupling to the pacemaker should be.

\section{Rings of N units}
In Appendix B.3, we have summarized the results of a ring configuration with directed coupling
obtained on closing a chain of unidirectionally coupled units. Here we expand on the discussion and give supporting numerical results. The schematic of this configuration is shown in Fig.~1(b) of the main text.
The elements of the coupling matrix $K$ of such a configuration are

\begin{figure}[ht!]
     \centering
 \includegraphics[width=\textwidth]{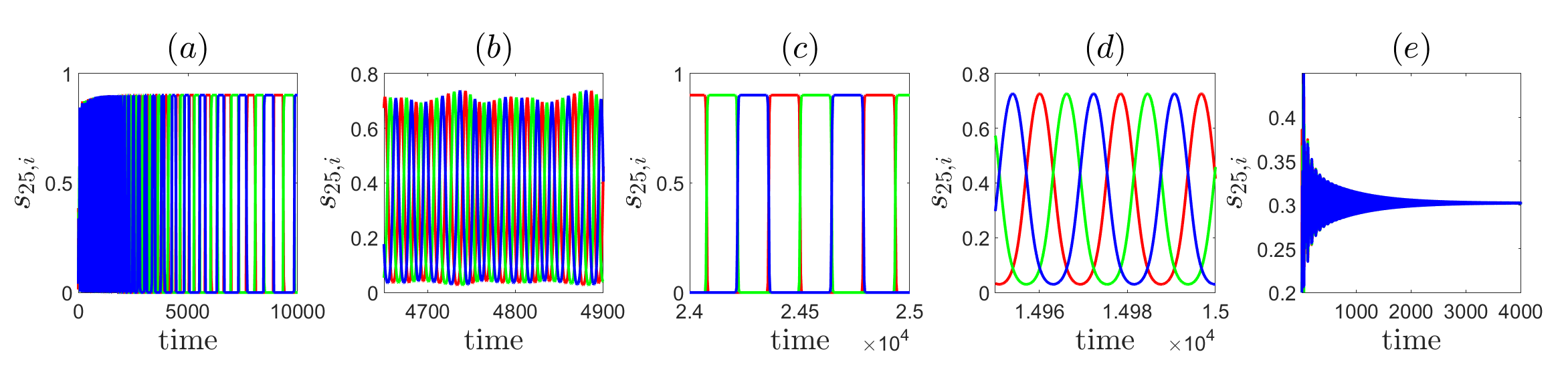}
\caption{Three items of the unit at position $k=25$ in a unidirectionally coupled ring of size $N=50$ for different values of $\delta_P$. (a) For $\delta_P=0.001$ all units are entrained to heteroclinic dynamics. Quasi-periodic dynamics in (b) for $\delta_P=0.03$, but it coexists with the heteroclinic dynamics as shown for another initial condition in (c). Limit cycle behavior  for $\delta_P=0.07$ (d) which loses stability via a Hopf bifurcation and gives way to the coexistence equilibrium as shown in (e) for $\delta_P=0.15$.
Other parameters are $\delta=0.75$, $\gamma_P=1.05$, $\gamma_D=1.11$, $c=2$, $e=0.2$, and $\rho = 1.0$.
}
 \label{fig:ring}
  \end{figure}

$K_{k,l} =
\begin{cases}
    \delta,& \text{if } k\in \{k_p+1,k_p+2,...,k_p+N-1 \} \text{ and } l=k-1\\
    \delta_P,& \text{if } k=k_p \text{ and } l=N+k_p-1\\
        0,              & \text{otherwise},
\end{cases}
$

\noindent where $k_p$ is the location of the pacemaker on the ring, chosen here as $k_p=1$.
The indices $k$ and $l$ are in terms of the elements of quotient rings of integer numbers.
Fig.~\ref{fig:ring} shows how the dynamics of this system changes as we tune $\delta_P$ for a fixed value of $\delta$ and other parameters for $N=50$.  For low values of $\delta_P$, all units are entrained to heteroclinic dynamics. An example is shown in Fig.~\ref{fig:ring}(a) for $\delta_P=0.001$ where we show the time series of the unit at the location $k=25$. The entrainment also depends on the size of the ring. Larger systems may not be synchronized to heteroclinic dynamics if we keep all the other parameters the same. As $\delta_P$ is increased, the system shows quasi-periodic dynamics (see Fig.~\ref{fig:ring}(b) for $\delta_P=0.03$) (also visible in the frequency spectrum, not displayed), but it coexists with the heteroclinic dynamics as shown for another initial condition in Fig.~\ref{fig:ring}(c). Upon further increasing $\delta_P$, the units show limit cycle behavior (an example is shown in Fig.~\ref{fig:ring}(d) for $\delta_P=0.07$) which loses stability via a Hopf bifurcation and gives way to a stable coexistence equilibrium as shown in Fig.~\ref{fig:ring}(e) for $\delta_P=0.15$.

\section{Bidirectional coupling along a chain of N units}
\begin{figure}[ht!]
     \centering
 \includegraphics[width=\textwidth]{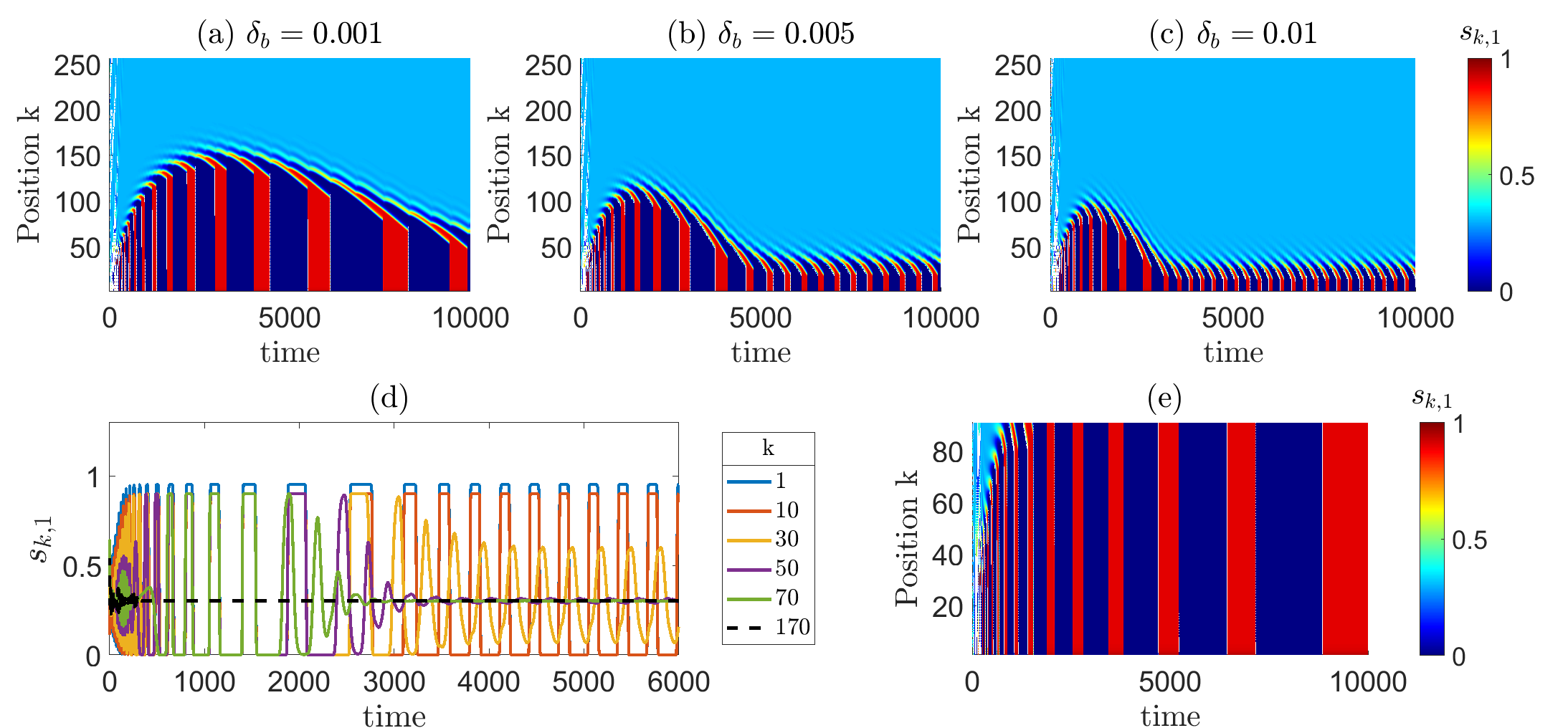}
\caption{Spatiotemporal plots of the first items $s_{k,1}$ of each unit of the bidirectionally coupled chain of length $N=256$ for three values of  $\delta_b$.
 When the size of the chain is reduced to $N=91$ units and $\delta_b=0.01$, all units get entrained (d).
Other parameters are $\delta=0.75$, $\gamma_P=1.05$, $\gamma_D=1.11$, $c=2$, $e=0.2$, and $\rho = 1.0$.}
 \label{fig:bidir_chain}
  \end{figure}

Here we discuss the results of Appendix B.4 with $N$ 1L-units coupled bidirectionally along a chain with a much weaker backward coupling compared to the forward coupling.
The results are shown in Fig.~\ref{fig:bidir_chain} for $N=256$. In Fig.~\ref{fig:bidir_chain}(a)-(c), we show the spatiotemporal plots of the first items $s_{k,1}$ of each unit for three values of the back coupling $\delta_b$.
When $\delta_b$ is finite, each unit is also affected by the dynamics of the units to its right.  This way, the entrainment length becomes sensitive to the size of the chain, differently from the unidirectional coupling.
Initially more and more units  get entrained along the chain, but after reaching a peak value over time in the number of entrained units, the recently entrained ones start de-entraining when the feedback from the units to the right reaches them. Finally the entrainment length settles to a smaller value. The ``peak'' of the entrainment decreases with an increase in the back-coupling, as for stronger back-coupling,  less units can be temporarily entrained.

If we reduce the size of the chain to $N=91$ units for $\delta_b=0.01$, all units can be entrained to heteroclinic dynamics as shown in Fig.~\ref{fig:bidir_chain}(d). An estimate of the size of the chain that can sustain the synchronized heteroclinic dynamics throughout the chain for a given $\delta$ and $\delta_b$ is indicated by the value of the peak. For example, if maximally $N=140$ units can be entrained over time from a larger set of units, reducing the chain to this number guarantees full entrainment. It should be noticed that the entrainment length also depends on the initial conditions similar to what we have seen for the unidirectional chain.

\section{Distance-dependent couplings between the pacemaker and the driven units}\label{subsec33}
\begin{figure}[ht!]
     \centering
  \includegraphics[width=0.45\textwidth]{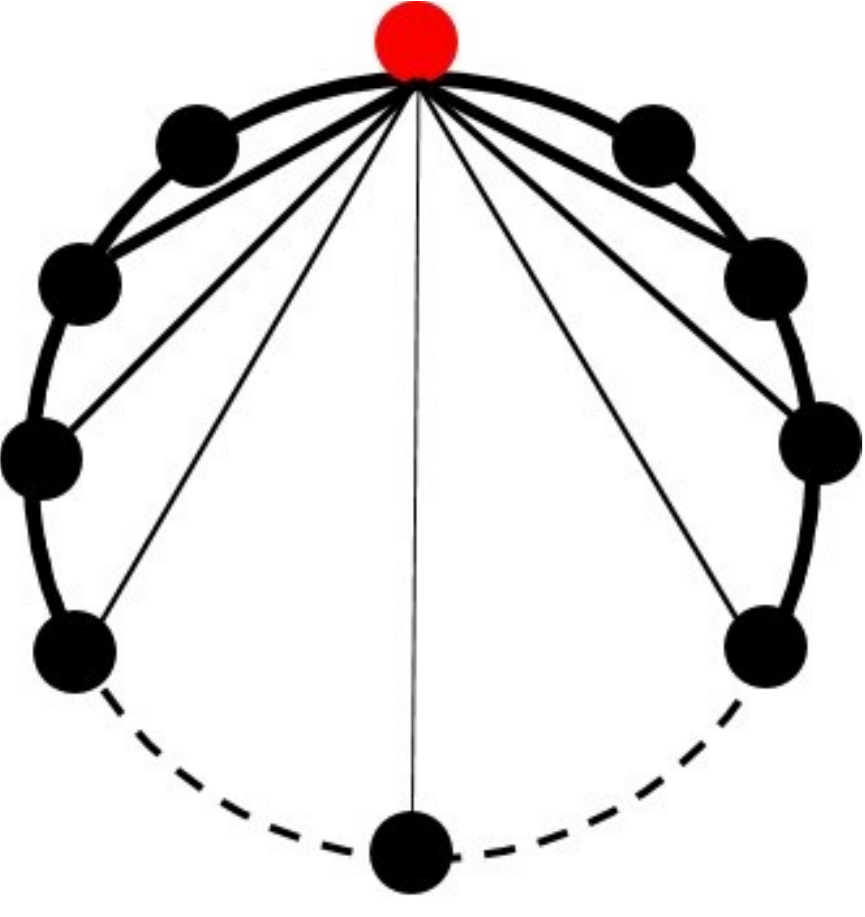}
\caption{A ring configuration with a distance-dependent forcing from the pacemaker to the other units. There is also a small distance-dependent back-coupling from each unit to the pacemaker. The units closer to the pacemaker are coupled more strongly to it, represented by thicker lines in comparison to the connections with units farther away. The interaction between the driven units is confined to only nearest neighbors.}
 \label{fig:schematic_nonlocal}
  \end{figure}

We have also studied a system of $N$ 1L-units in a ring configuration with distance-dependent coupling between the pacemaker and the other $N-1$ driven units. The schematic of this configuration is shown in Fig.~\ref{fig:schematic_nonlocal}. For this implementation, the effects of the back-coupling  on the phase space can be more versatile and lead to four different regimes that we describe in this section.
The coupling from the pacemaker to the driven units ($\delta$) as well as the back-coupling from the driven units to the pacemaker decay with  distance according to $1/r$. Since it is a ring configuration, the location of the pacemaker on the ring does not matter, again it is labelled as $k=1$.
The elements of the coupling matrix $K$ thus are

\begin{figure}[ht!]
     \centering
 \includegraphics[width=\textwidth]{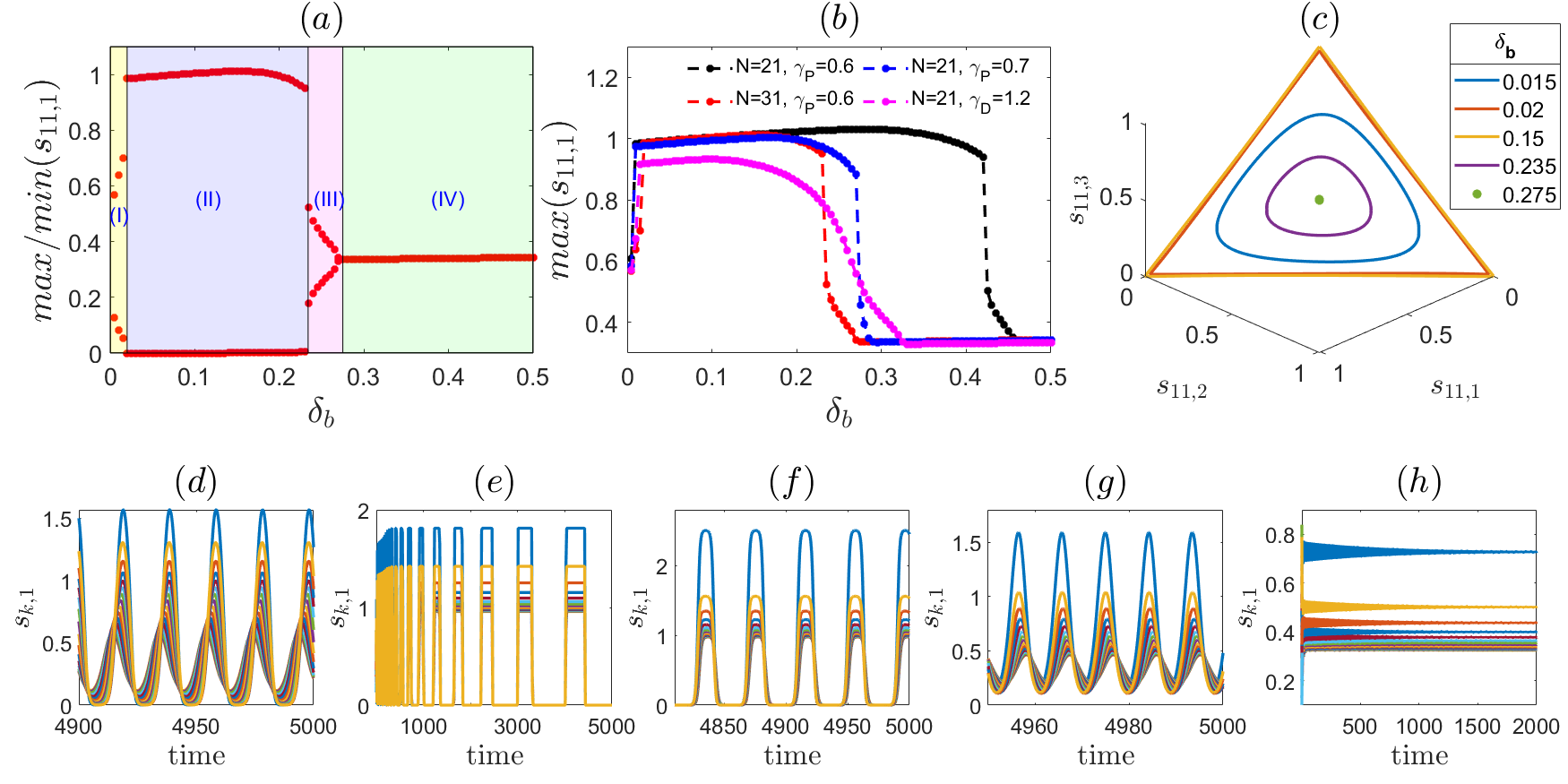}
\caption{Results for the ring configuration with the distance-dependent coupling. In (a), we take $N=31$ and plot the maximum and the minimum value of the first item of $11^{th}$ unit ($s_{11,1}$) against the back-coupling $\delta_b$. There are four regions with distinct dynamical behavior. Region-(I) has limit cycle behavior as shown in (d) for $\delta_b=0.015$. In region-(II), units show synchronized heteroclinic dynamics as shown in (e) for $\delta_b = 0.02$ and (f) for $\delta_b = 0.15$. In region-(III), the trajectories converge to limit cycles  far from the heteroclinic contour. An example is shown for $\delta_b=0.235$ in (g). In region-(IV), the trajectories of all units settle to their corresponding coexistence equilibrium as shown in (h) for $\delta_b = 0.275$. In (b), we compare the results for two different sizes of the ring $N=\{21,31\}$, and two different values of $\gamma_P=\{0.6,0.7\}$ and $\gamma_D=\{1.11,1.2\}$. The phase-space trajectory of the unit at position $k=11$ is shown in (c) for five values of $\delta_b$ corresponding to (d)-(h). Other parameters are $\delta=0.5$, $\gamma_P=0.6$, $\gamma_D=1.11$, $c=2$, $e=0.2$, and $\rho = 1.0$.
}
 \label{fig:nonlocal_fig}
  \end{figure}

$K_{k,l} =
\begin{cases}
    \frac{\delta}{r_{k,l}},& \text{if } k\in \{2,...,N \} \text{ and } l=1\\
    \frac{\delta_b}{r_{k,l}},& \text{if } k=1 \text{ and } l\in \{2,...,N \}\\
    \delta ,& \text{if } l=\{k-1,k+1\} \text{ and } k,l\in \{2,...,N \}\\
    0,              & \text{otherwise,}
\end{cases}
$
where  $r_{k,l}= \text{min}(|k-l|,N-|k-l|)$.
We fix the coupling strength at $\delta=0.5$ and slowly increase the back coupling strength $\delta_b$ from zero.
The results are shown in Fig.~\ref{fig:nonlocal_fig}. We distinguish four regions with distinct dynamical behavior. We show them in Fig.~\ref{fig:nonlocal_fig}(a), where we plot the maximum and the minimum value of the first item of the $11^{th}$ unit ($s_{11,1}$) against the back coupling $\delta_b$. The choice of the item $s_{11,1}$ is purely arbitrary as all the units have similar dynamical behavior in each region and only their amplitudes differ.
For very small values of $\delta_b$, i.e. in region-(I), all units entrain to limit cycle behavior. The frequencies of all  units are synchronized while their phases are not. An example is shown in Fig.~\ref{fig:nonlocal_fig}(d) for $\delta_b=0.015$ where we plot the time series of the first item of units at all locations. As the back coupling exceeds a certain threshold and $\delta_b$ enters region-(II), there is a gap, all  units synchronize to heteroclinic dynamics. An example is shown for $\delta_b = 0.02$ in Fig.~\ref{fig:nonlocal_fig}(e). As we further increase $\delta_b$ inside region-(II), the characteristic slowing down of heteroclinic dynamics becomes less and less pronounced, and the trajectories seem to converge to a limit cycle in the vicinity of the heteroclinic contour. An example is shown in Fig.~\ref{fig:nonlocal_fig}(f) for $\delta_b = 0.15$. As we increase $\delta_b$ further and enter region-(III), there is another gap: the values of the amplitudes drop. The trajectories still converge to limit cycles but far from the heteroclinic contour. An example is shown for $\delta_b=0.235$ in Fig.~\ref{fig:nonlocal_fig}(g). On increasing $\delta_b$ further, the limit cycles lose stability via a Hopf bifurcation. The trajectories of all units settle to their corresponding coexistence equilibrium in region-(IV) as shown in Fig.~\ref{fig:nonlocal_fig}(h) for $\delta_b = 0.275$. The phase-space trajectory of the unit at position $k=11$ is shown in Fig.~\ref{fig:nonlocal_fig}(c) for five values of $\delta_b$ corresponding to Fig.~\ref{fig:nonlocal_fig}(d)-(h).

In Fig.~\ref{fig:nonlocal_fig}(b), we compare the results for two different sizes of the ring $N=\{21,31\}$, and different values of $\gamma_P=\{0.6,0.7\}$ and $\gamma_D=\{1.11,1.2\}$. The entrainment length (that is the size of region-(II)) decreases with an increase in the size of the ring and an increase in $\gamma_P$. On increasing $\gamma_D$, the transition from region-(II) to region-(III) becomes gradual.

\section{Noise acting as a control parameter}
\begin{figure}
    \centering
 {\includegraphics[width = 0.194 \textwidth,height=1in]{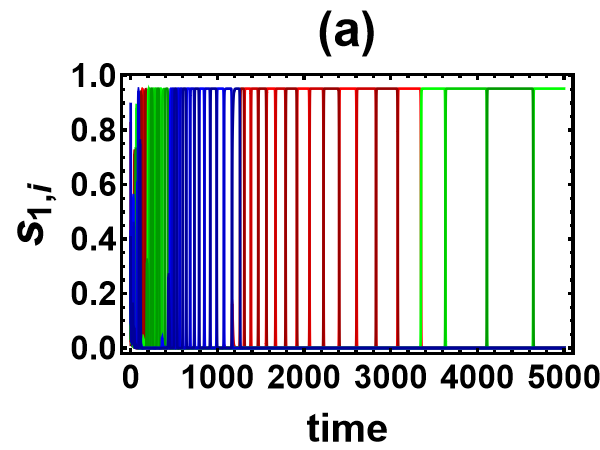}} \hspace{-2mm}
    {\includegraphics[width = 0.194 \textwidth,height=1in]{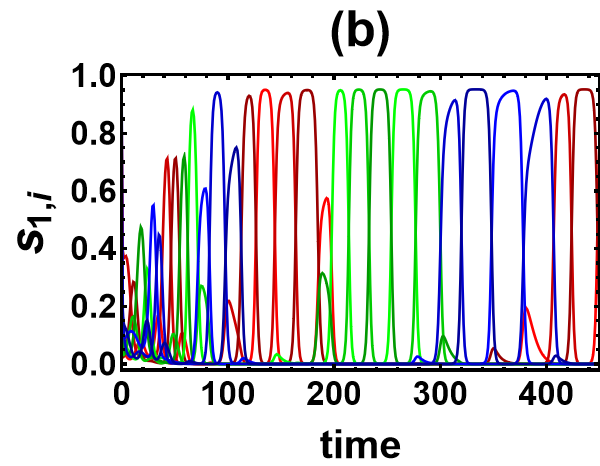}} \hspace{-2mm}
    {\includegraphics[width = 0.194 \textwidth,height=1in]{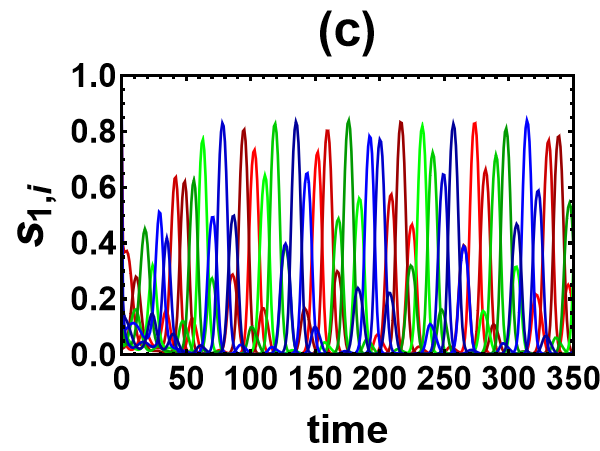}}   \hspace{-2mm}
    {\includegraphics[width = 0.194 \textwidth,height=1in]{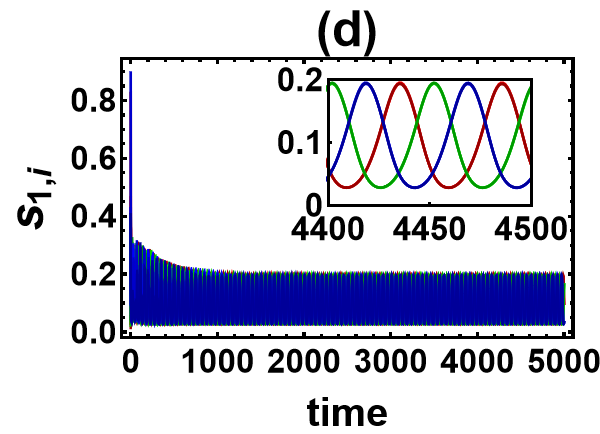}}\hspace{-2mm}
    {\includegraphics[width = 0.194 \textwidth,height=1in]{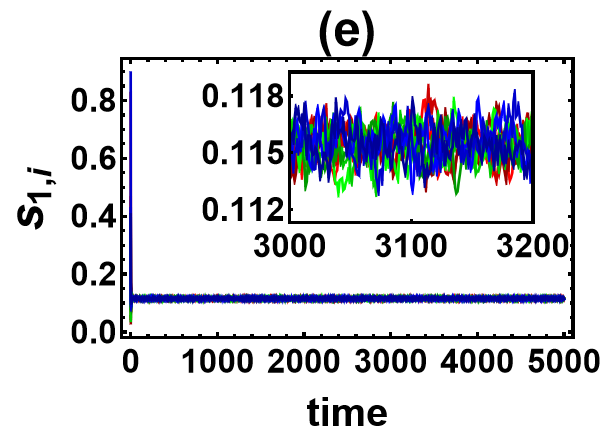}}\\
 {\includegraphics[width = 0.194 \textwidth,height=1in]{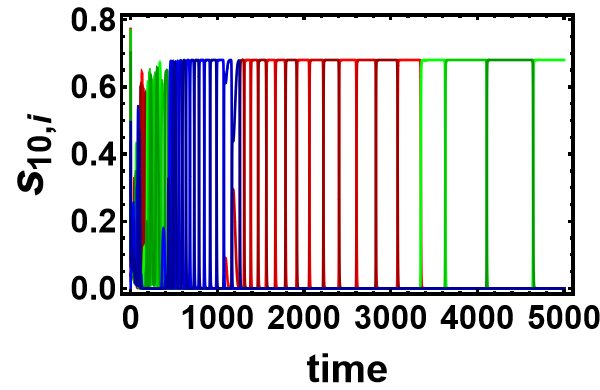}}  \hspace{-2mm}
    {\includegraphics[width = 0.194 \textwidth,height=1in]{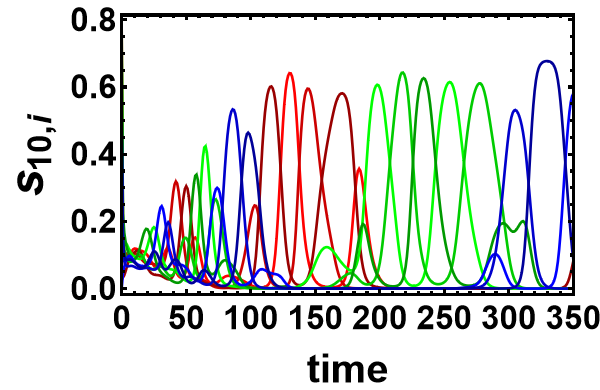}}  \hspace{-2mm}
    {\includegraphics[width = 0.194 \textwidth,height=1in]{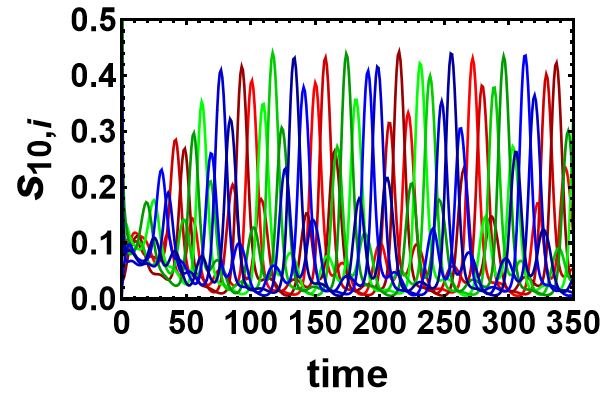}} \hspace{-2mm}
    {\includegraphics[width = 0.194 \textwidth,height=1in]{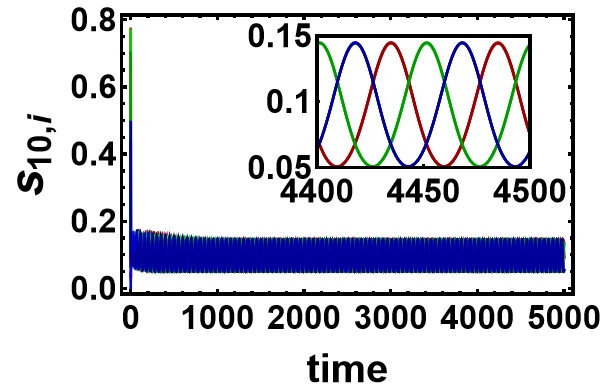}} \hspace{-2mm}
    {\includegraphics[width = 0.194 \textwidth,height=1in]{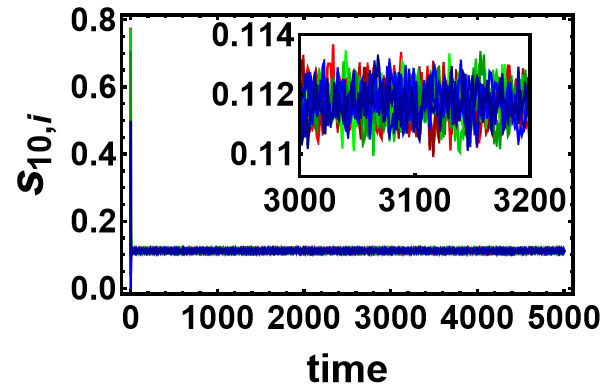}}\\
     {\includegraphics[width = 0.194 \textwidth,height=1in]{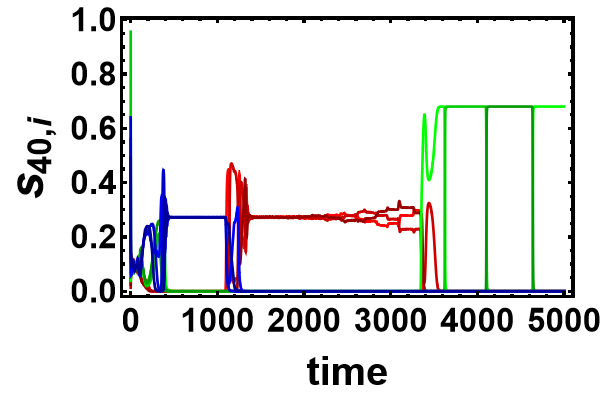}}\hspace{-2mm}
    {\includegraphics[width = 0.194 \textwidth,height=1in]{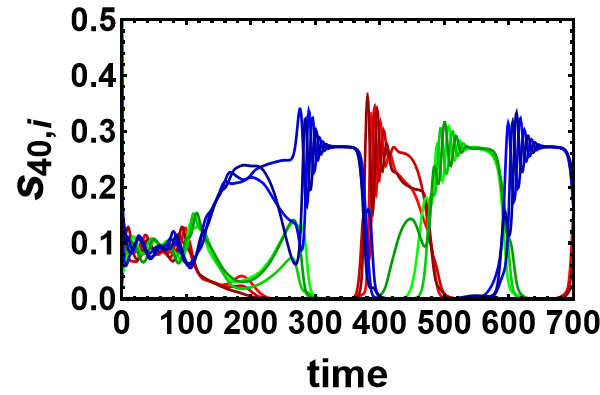}} \hspace{-2mm}
    {\includegraphics[width = 0.194 \textwidth,height=1in]{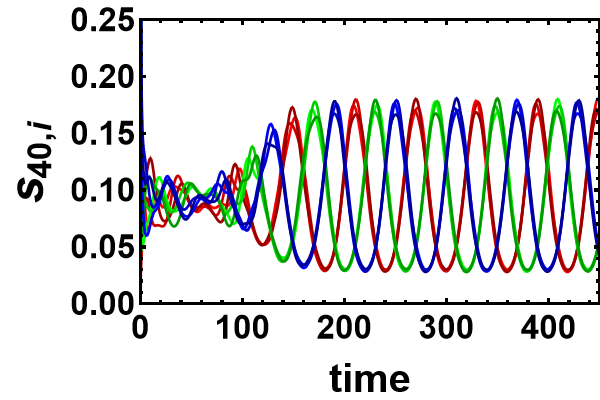}}\hspace{-2mm}
    {\includegraphics[width = 0.194 \textwidth,height=1in]{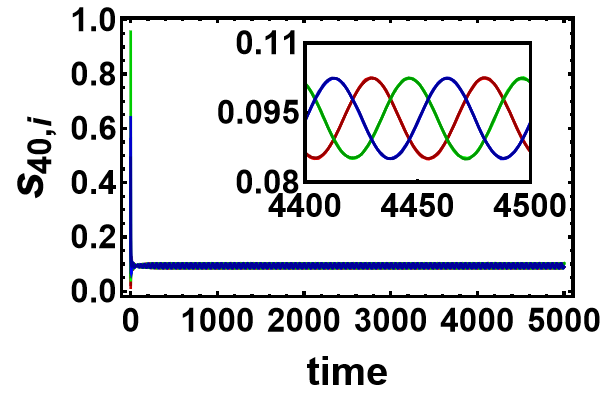}}\hspace{-2mm}
    {\includegraphics[width = 0.194 \textwidth,height=1in]{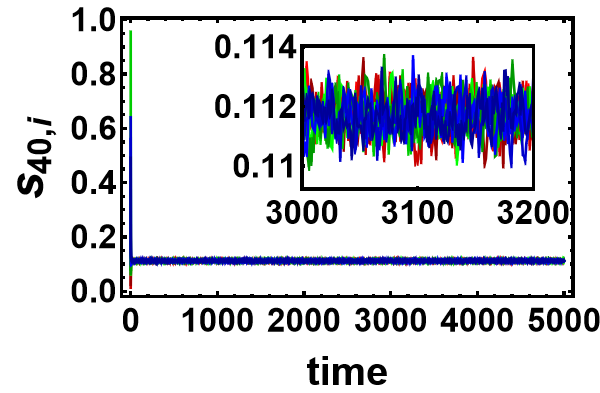}}\\
     {\includegraphics[width = 0.194 \textwidth,height=1in]{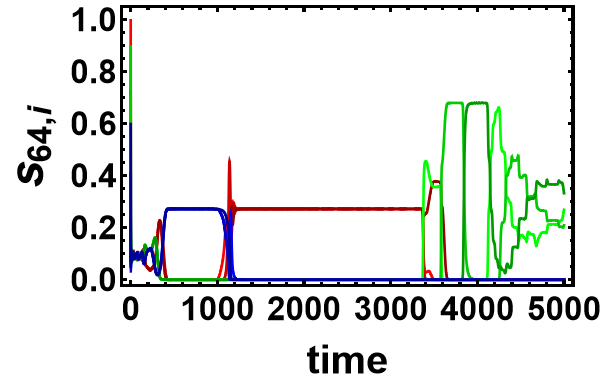}}\hspace{-2mm}
    {\includegraphics[width = 0.194 \textwidth,height=1in]{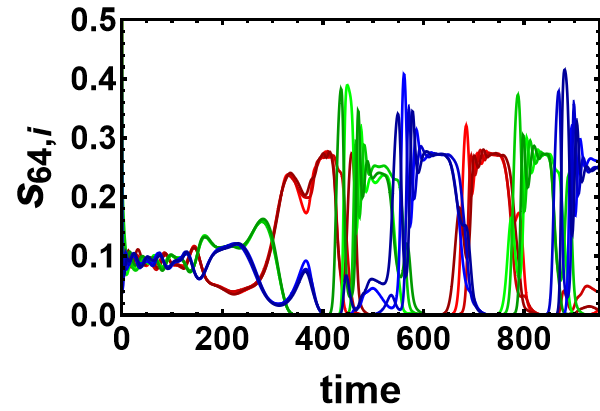}} \hspace{-2mm}
    {\includegraphics[width = 0.194 \textwidth,height=1in]{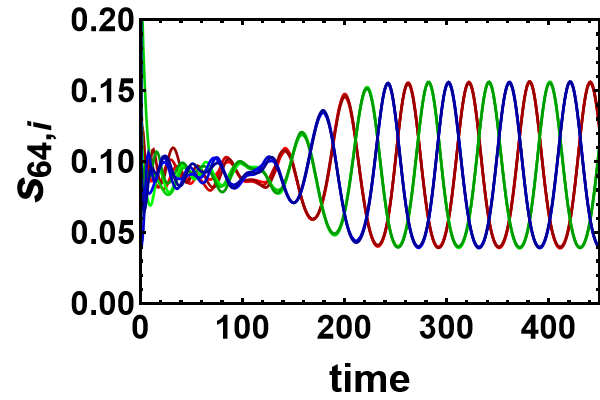}}\hspace{-2mm}
    {\includegraphics[width = 0.194 \textwidth,height=1in]{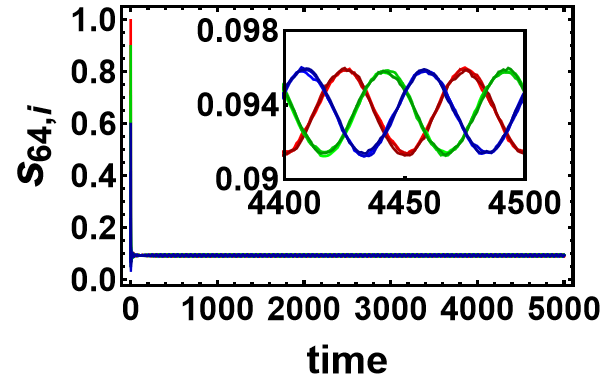}}\hspace{-2mm}
    {\includegraphics[width = 0.194 \textwidth,height=1in]{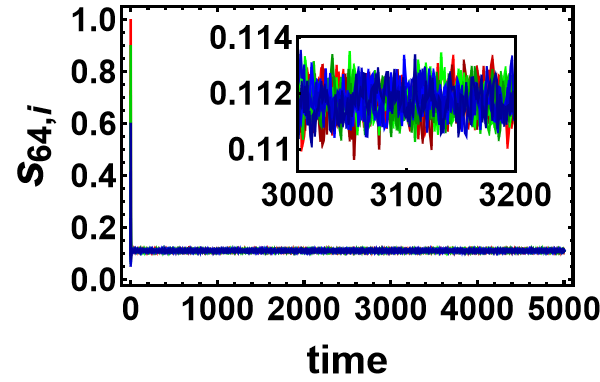}}\\
    {\includegraphics[width = 0.4\textwidth,height=1in]{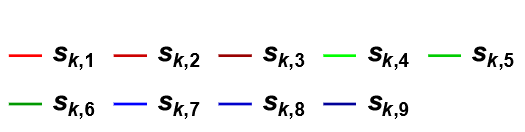}}
     \caption{The figure shows the reduction of the hierarchy levels in a unidirectionally coupled chain of $N=64$ 2L-units as the noise strength is increased from (a) $\sigma=0$, to
     (b) $\sigma=10^{-12}$, (c) $\sigma=10^{-5}$, (d) $\sigma=10^{-4}$, and (e) $\sigma=10^{-3}$. The first row shows the time series of all nine items of the pacemaker (located at $k=1$) for each noise level. Second, third, and fourth rows show the time series of the units at spatial locations $k=10, 40$, and $64$ respectively. In the inset plots of (d) and (e), we zoom into the time series in a small time window. The hierarchy level of the pacemaker collapses from two levels (for $\sigma=10^{-12},10^{-5}$) to one level (for $\sigma=10^{-4}$) and finally to zero level (for $\sigma=10^{-3}$). The units closer to the pacemaker are entrained to the same hierarchy level as the pacemaker, those farther away could be entrained to a lower level. For example, see the time series of $k=40,64$ in the third and fourth rows of (b) and (c). Other parameters are $\gamma_P=1.05$, $\gamma_D=1.47$, $\delta=0.5$, $r=1.25$, $e=0.2$, $f=0.3$, $c=d=2$, and $\rho=1$. }
    \label{fig:hier_noise}
\end{figure}

In  section 4.4 of the main text we discuss the role of noise as a control parameter to steer the dynamics of the system between different levels of hierarchy. Noise, here applied to all units, drives the pacemaker from a 2L-unit to a 1L-unit to a CE-unit and along with that the entrained units assigned to the chain.
Here we discuss these results in more detail and plot them in Fig.~\ref{fig:hier_noise}.
The units are unidirectionally coupled along a chain of length $N=64$, the pacemaker is again positioned at $k=1$.  The pacemaker is a 2L-unit with $\gamma_P=1.05$ and all the other units are in the CE-regime with $\gamma_D=1.47$. The first row of Fig.~\ref{fig:hier_noise} shows the time series of the pacemaker (at $k=1$) for five values of noise strength $\sigma$. Second, third, and fourth rows show the time series of units located at position $k=10,40,64$ respectively. In column (a) we show the results when there is no noise present in the system. The pacemaker is able to entrain the units to hierarchical heteroclinic dynamics for sufficiently strong coupling strength. However, when the dynamics of the pacemaker slows down so much that it gives sufficient time to the trajectories of the driven units to travel to more stable states (equilibria with multiple items), the driven units get de-entrained. If we add a very small amount of noise to the system, the slowing down is prevented and the pacemaker still has two levels of hierarchy for low  noise levels ($\sigma=10^{-12}$ in Fig.~\ref{fig:hier_noise}(b) and $\sigma=10^{-5}$ in Fig.~\ref{fig:hier_noise}(c)). The units closer to the pacemaker are also entrained to the dynamics with two levels of hierarchy (see the second row of Fig.~\ref{fig:hier_noise}(b) and Fig.~\ref{fig:hier_noise}(c)), while the units far away only have a single level of hierarchy as the influence of the pacemaker with distance decreases (see the third and last row of Fig.~\ref{fig:hier_noise}(b) and Fig.~\ref{fig:hier_noise}(c)). When the noise is increased to $\sigma=10^{-4}$ in Fig.~\ref{fig:hier_noise}(d), the hierarchy level of the pacemaker collapses to a single level, thus all the driven units are entrained to the dynamics with a single level of hierarchy. The three SHCs are replaced by three 3-items coexistence fixed points connected by an effective limit cycle. On further increasing the noise to $\sigma=10^{-3}$ in Fig.~\ref{fig:hier_noise}(e), the noise drives the pacemaker to the vicinity of its global coexistence equilibrium. The driven units, which are already in the CE-regime, stay near their coexistence equilibrium.
